\documentclass[11pt]{article}

\usepackage{enumerate,amsmath,amssymb,amsthm,color}

\makeindex             

\usepackage{enumerate,amsmath,amssymb,amsthm,color}

\usepackage{hyperref}

\usepackage{mathptmx}       
\usepackage{helvet}         
\usepackage{courier}        
\usepackage{type1cm}        
\usepackage{makeidx}         
\usepackage{graphicx}        
\usepackage{subfigure}

\usepackage{multicol}        
\usepackage[bottom]{footmisc}

\usepackage{authblk}
\usepackage{caption}

\DeclareMathOperator{\sech}{sech}

\def\Na{Na$^{+}$}
\def\Ca{Ca$^{2+}$}
\def\Cl{Cl$^{-}$}
\def\K{K$^{+}$}

\def\vecn{\vec{n}}
\def\vecx{\vec{x}}
\def\vecr{\vec{r}}
\def\dd{\displaystyle}

\def\8{\infty}
\def\f{\frac}
\def\p{\partial}

\makeindex             

\begin{document}

\title{Variational Methods for Biomolecular Modeling} 
\author[1]{Guo-Wei Wei \thanks{wei@math.msu.edu}}
\author[2]{Y. C. Zhou \thanks{yzhou@math.colostate.edu}}
\affil[1]{Department of Mathematics, Michigan State University, East Lansing, MI 48824}
\affil[2]{Department of Mathematics, Colorado State University, Fort Collins, CO 80523}
\maketitle

\abstract{Structure, function and dynamics of many biomolecular systems can be characterized by the energetic variational principle 
and the corresponding systems of partial differential equations (PDEs). This principle allows us to focus on the identification of 
essential energetic components, the optimal parametrization of energies, and the efficient computational implementation 
of energy variation or minimization. Given the fact that complex biomolecular systems are structurally non-uniform 
and their interactions occur through contact interfaces, their free energies are associated with various interfaces 
as well, such as solute-solvent interface, molecular binding interface, lipid domain interface, and membrane 
surfaces. This fact motivates the inclusion of interface geometry, particular its curvatures, to the parametrization of  
free energies. Applications of such interface geometry based energetic variational principles are illustrated through three concrete topics: 
the multiscale modeling of biomolecular electrostatics and solvation that includes the curvature energy of the molecular 
surface, the formation of microdomains on lipid membrane due to the geometric and molecular mechanics at 
the lipid interface, and  the mean curvature driven protein localization on membrane surfaces. By 
further implicitly representing the interface using a phase field function over the entire 
domain, one can simulate the dynamics of the interface and the corresponding energy variation by evolving the phase field 
function, achieving significant reduction of the number of degrees of freedom and computational complexity.  Strategies for improving 
the efficiency of computational implementations and for extending applications to coarse-graining or multiscale molecular simulations are outlined.
}

\newpage

{\setcounter{tocdepth}{4} \tableofcontents}
\newpage

\section{Introduction} \label{sec:intro}

Living biological systems require a constantly supply of energy to generate and maintain certain biological orders 
that keep the systems alive. This warrants the biophysical models that quantify the management and 
balance of energy in biological systems, i.e., the energy budget of metabolism. Taking cells 
- the building blocks of life - as an example, energy is derived from the chemical bond energy in food 
molecules, passed through a sequence of biochemical reactions, and is used in cells to produce activated 
energy carrier molecules (i.e., ATPs) for powering almost every activity of the cells, including muscle contraction,  generation of 
electricity in nerves, and DNA replication \cite{Molecularbiologyofthecell}. For  solvated  biomolecular systems \footnote{Water constitutes a large 
percentage of cellular mass and therefore biomolecules are mostly living in an aqueous environment where various types of ions such as 
sodium (\Na), potassium (\K), calcium (\Ca), and  chloride (\Cl) present at different concentrations.} discussed in this chapter, 
including solvated proteins,  bilayer membranes, or their complexes, one can make similar energy budgets too. 
Various types of energies can be identified for   biomolecular systems, such as
\begin{enumerate}
\item kinetic energies of atoms or molecules in motion;
\item potential energies for  bonded atoms: potential energies characterizing the stretching, bending, torsion of the
covalent bonds between atoms;
\item potential energies for unbounded atoms: electrostatic energy and van der Waals energy; and 
\item kinetic and potential energy interconversions in enzymatic processes and chemical reactions.
\end{enumerate}
The first three energy terms constitute the basis for molecular dynamics (MD) simulations of non-reactive solvated biomolecular systems. Using 
the spatial coordinates of individual atoms as parameters, MD simulations trace the motion of each atom by using the Newton second 
law, where the force applied to each atom is computed as the variational of the total energy with respect to the atom's spatial 
coordinates \cite{CHARMM,CaseD2002,GROMACS,NAMD}. Additional forces that models  temperature-dependent thermal fluctuations can be 
added, giving rise to Langevin dynamics simulations \cite{Schlickbook}. In this regard,  MD 
simulation is  indeed a classical application of the variational principle. 

The large amount of solvent molecules in a molecular dynamics simulation of solvated biomolecular system can make the  
simulation daunting and expensive. This deficiency motivates the development of various continuum or multiscale models for part 
of or the entire solvated biomolecular system \cite{TangY2006a,DzubiellaJ2006PRL,Grabe2008a,Grabe2012a,ChenX2010a,WeiG2010a,ZhouY2010b,ChenZ2012a,FeigM2004,SkynerR2015a}. Notably among these simplifications are  implicit solvent models, which manage to replace the atomic degrees of 
freedom of solvent molecules with a continuum description of averaged behavior of solvent molecules while retain an 
atomistic description of the solute molecule \cite{FeigM2004,SkynerR2015a}. Accordingly, the solvent-solute interface must be identified as the boundary 
between the continuum solvent region and the discrete biomolecular domain. This interface is of particular importance 
because it is related to a range of solvent-solution interactions such as hydrogen bonding, ion-ion, ion-dipole, dipole-dipole 
and multipole interactions, and Debye attractions \cite{DillK_MolecularDrivingForce}. Thus the parametrization of the total energy 
of the system must include the geometry of this interface. Mean and Gaussian curvatures are generally involved in such parametrization because  
they measure the variability or non-flatness of a biomolecular surface and characterize respectively the extrinsic and intrinsic 
measures of the surface \cite{Wolfgang_DiffGeo}.  In these multiscale models of solvated biomolecules systems the motion of 
the atoms still follows the Newton's law  where the force is given as the variational of the total energy with respect to the atoms' spatial coordinates, the electrostatic potential, and the 
interface \cite{Gilson:1993,Wagoner06,WeiG2010a,Geng:2011,ZhouS2014a}. The change in the solvent-solute interface induces 
variation in curvatures, whose energies might be treated as a part of the total energy functional. 
These curvature based or differential geometry based biomolecular models offer a manifest of mathematical analysis and 
computational methodologies for  the dynamics of the solvent-solute interface and the equilibrium energy landscape of solvated 
biomolecules. In other words, one can derive dynamic partial differential equations to evolve the interface morphology, 
and this evolution can be mapped to the path toward the global or local minimum on the landscape of the total energy. Here 
in this chapter we shall present three representative applications of interface geometry based variational principles to the modeling 
of biomolecular interactions: (i) biomolecular electrostatics and solvation, (ii) surface microdomain formation in 
bilayer membranes, and (iii) curvature driven protein localization in bilayer membranes.

In the first application we consider the long-range electrostatic interactions among partially charged static atoms in the solute and  the aqueous solvent with mobile ions. These interactions strongly depend on the position of solvent-solute boundary, also  referred to as the molecular surface in this context, where a rapid transition of dielectric permittivity is observed. Inclusion of this interface, albeit implicitly, in the formulation of the total energy of the system facilitates the coupling of polar and nonpolar solvent-solute interactions, as well as  the nonlinear solvent response, in the form of interface energy functional of surface curvature energy, electrostatic energy and van der Waals potential. Such a coupling finally gives rise to a novel variational multiscale solvation model \cite{DzubiellaJ2006PRL,Dzubiella06a,WeiG2010a,ChenZ2010a,ChenZ2011a}. In a more elaborated model, the solute molecule can be described in further detail by using the quantum density functional theory (DFT) in an iterative manner, which allows a more accurate account of solvent-solute interaction and response \cite{ZhanChen:2011a}. Differential geometry based solvation models have been shown to deliver superb predictions of solvation free energies for hundreds of molecules \cite{ChenZ2012a,BaoWang:2015a}.    This variational principle based solvation model can be further extended to describe essential biological transportation such as transmembrane ion or proton flows that depend critically on the geometry of the associated protein channels. By including the chemical potential and entropy of the diffusive ion species into the total energy functional one can obtain simultaneously the optimized channel protein surfaces as well as the corresponding I-V (current-voltage) curve \cite{ZhengQ2011a,ChenZ2012a,WeiG2012a}.

Curvature is believed to play an important role in many biological processes, such as protein-DNA and protein-membrane interactions, including membrane curvature sensing. Classical phase field modeling of surface pattern formation in bilayer membranes contains a curvature term in its definition of the total energy \cite{ChenX1994a,DuQ2005c,DuQ2005d,ParthasaR2006a,CampeloF2007a,GaoL2009a}. However, when modeling the surface pattern formation in our second application here, we  show that  it is the geodesic curvature rather than the curvature of pattern interfaces that plays an essential role in modulating the interface energy. Noting that this geodesic curvature is defined on a general differentiable manifold, and thus the classical  phase field modeling of phase separation with specified intrinsic curvature can be regarded as a special case of this geodesic curvature model in the Euclidean spaces. By providing various intrinsic geodesic curvatures that model the geometry of the contact of different species of lipids, we are able to simulate the generation of lipid rafts as the formation and equalization of localized surface domains.

In contrast to most amphiphilic lipids whose relatively long and geometrically regular hydrophobic tails allow they to pack  together, membrane proteins usually do not present in large distinct domains in membrane surfaces, although small amount of membrane proteins can compound together forming functional complexes such as ion channels or membrane transporters. Most membrane proteins have amphipathic helices, which contain both hydrophobic and hydrophilic groups, complementing to   amphiphilic lipids.  Therefore, the localization of these membrane proteins in general can not be modeled using the geodesic curvature based phase separation model as described  in our second application. Many membrane proteins, however, do prefer bilayer membranes with particular curvature, in the sense that they can induce particular curvature in the bilayer membrane and they tend to be localized in regions with specific curvature. Therefore, one can imagine that membrane curvature can provide a driving force for the distribution of membrane proteins in the bilayer, and thus an appropriate energy functional that represents the membrane curvature must be added to the classical electrochemical potential and entropy to describe the localization of membrane proteins. 

These three applications of variational principles in biomolecular modeling are by no means exhaustive, even in the context  of 
solvation analysis and membrane-protein interactions. There are inspiring studies of ion and water transport in membrane channels using 
energetic variational approaches, where the effects of surface charge density and non-uniform particle sizes can be readily included in investigations 
thanks to the flexibility of variational approaches 
\cite{WeiG2010a,HorngT2012a,XuS2014a,HyonY2014a,HyonY2011a,HyonY2012b,ZhouY2011a,LiB2009b,WeiG2012a}. Similar flexibility also enables the extension of the application of variational principles from the standard phase field modeling of bilayer membrane deformation and morphology \cite{DuQ2006a,DuQ2005c,DuQ2005d} to multi-components membranes \cite{LiS2012a,ZhaoY2011a}, pore formation \cite{RyhamR2013a,CohenF2012a}, and double layer \cite{DaiS2013a,GavishN2011a}. Some of models, particular those for bilayer membranes, share various degree of similarity to the models used for self-assembly or phase separation of polymers or co-polymers. It is this wide diversity of lipid structures and the complicated interactions between proteins and lipid bilayers in solution that makes the energetic variational modeling of bilayer membranes unique and challenging. As we shall present below, most of our efforts are concentrated on the formulation of potential energy functional of these interactions so that the variational principle can be applied and numerical solutions can be found by solving the corresponding systems of nonlinear partial differential equations (PDEs).

\section{Variational Multiscale Methods for Biomolecular Electrostatics and Solvation} \label{sec:solvation}

By definition, the solvation energy of biomolecules is the cost of free energy required to transfer the biomolecules from the vacuum to the solvent environment. It is therefore an essential quantitative characterization of the solute-solvent interactions.  Electrostatic free energy, also called polar solvation free energy,  is an important component of the solvation free energy since most biomolecules are charged and there are always mobile ions in the solvent under physiological conditions.  Various critical applications of the electrostatic and solvation free energies can be found in chemistry, biophysics, and medicine. We refer the reader  to \cite{Davis90a,Madura96,Luzhkov92,Warshel94,Warshel97,FeigM2004,WarshelA2006a,Wagoner06,JaramilloA2005a,LeeM2013a,RenP2012a,EdelsbrunnerH2005a,ChipotC2006a,ShirtsM2010a} for theoretical underpinning of these applications and the determination of the electrostatics and solvation free energies.  Apart from electrostatic effects, the solvation free energy also involves the nonpolar energy, namely, the energy cost for creating a suitable cavity in the continuum solvent to allow the transferring of the biomolecules and for the dispersive interactions between the solvent and the biomolecule on the surface of this cavity. Implicit solvent models are particularly appearing for computing the solvation free energy since the number of solvent degrees of freedom can be dramatically reduced by a well fitted bulk dielectric permittivity while the atomistic representations of solute biomolecules can be  retained to maintain a detailed  modeling of the solute.  The framework of    implicit solvent models allows  the solvation free energy to  be decomposed into two components, polar solvation and nonpolar solvation \cite{LeeM2013a, Wagoner06,LevyR2003a}. In this approach, the electrostatic contribution can be readily computed from the solution of the Poisson-Boltzmann equation, or  the Poisson equation if there is no explicit ion in the solvent \cite{Lu07e,Grochowski2008a,Wagoner04,Park03,Grant01a,Azuara06,Baker05a}.  The solution of these equation depends on the contrast of dielectric permittivity in vacuum and the solvent environments, and this contrast is concentrated at the boundary between the biomolecule and the solvent. Likewise, the calculation of nonpolar solvation free energy  depends on the geometry of the biomolecular surface.  The fact that both polar and nonpolar components are determined by the solvent-solute interface warrants the importance of a biophysically justifiable, mathematically well-posed, and computational feasible definition of the molecular surface or dielectric interface. In fact, the decoupling of polar and nonpolar components makes implicit solvent models conceptually convenient and computationally simple.

However, there are many structural imperfections  associated with implicit solvent models. First, intrinsic thermodynamical and kinetic  coupling makes it impossible to completely separate the electrostatic component from the non-electrostatic components in the solvation modeling. Additionally, a pre-prescribed solvent-solute interface, such as solvent excluded surface and van der Waals surface,  decouples polar and nonpolar components.   As a result, the solvation induced solute polarization and solvent response are not appropriately accounted in implicit solvent models. Moreover, implicit solvent models neglect potential solvation induced surface reconstruction and possible conformational changes. Finally, thermodynamically, the change in the Gibbs free energy of solvation can be formally decomposed into the change in internal energy, work, and entropy effect. There is no guarantee that all of these components are fully accounted  in implicit solvent models.  In addition to the aforementioned  structural or organizational imperfections, the performance of implicit solvent models is subject to a wide range of implementation deficiencies, such as  the modeling of nonpolar component, the treatment of the electrostatic component, the exclusion of high-order polarization, the exclusion of curvature,  the geometric singularity of solvent-solute interface, the stability of numerical schemes and algorithms, the grid convergence of the solvation free energy, to mention only a few.    

Some of the aforementioned problems have been the subjects of intensive study in the past few decades.   One approach starts from improving the surface definitions, so that earlier van der Waals surface, solvent accessible surface \cite{Lee:1971}, and molecular surface (MS) \cite{Richards:1977}  are replaced by smooth surface expressions  \cite{Grant01a,Grant:1995,Grant:2007,QZheng:2012,MXChen:2012}.  Geometric analysis, which combines differential  geometry (DG) and differential equations, is a powerful mathematical tool  for signal and image processing, data 
analysis, and surface construction \cite{SOsher:1988,Wei:1999a,YWang:2011c,YWang:2012a,YWang:2012b}. Geometric PDEs and DG theories of surfaces provide a natural and simple description for a solvent-solute interface. The first curvature-controlled PDEs  for molecular surface construction and solvation analysis was introduced  in 2005 \cite{Wei:2005}. A variational solvent-solute interface, namely a minimal molecular surface (MMS), was proposed for molecular surface generation in 2006  \cite{Bates:2006,Bates:2008}. In this work, the minimization of surface free energy is equivalent to the minimization of surface area, which can be implemented via the mean curvature flow, or the Laplace-Beltrami flow, and gives rise to the MMS. The MMS approach has been used in implicit solvent models \cite{ChenZ2012a, Bates:2008}. Potential-driven geometric flows, which admit potential driven terms, have also been proposed for biomolecular surface construction \cite{Bates:2009}. This 
approach was adopted by many researchers \cite{Cheng:2007e,CheJ2008,ChengL2009a,Yu:2008g,SZhao:2011a,SZhao:2014a} for biomolecular surface 
identification and electrostatics/solvation modeling.  

It is natural to extend DG based variational theory of the solvent-solute interface  into a full solvation model by incorporating a variational formulation of the PB theory \cite{Sharp:1990, Gilson:1993, WeiG2010a,ChenZ2012a} following the spirit of a similar approach by McCammon and coworkers \cite{Dzubiella06a,DzubiellaJ2006PRL}. However, the formalism of McCammon and coworkers does not involve geometric flow and has a Gaussian curvature term that might lead to jumps in the energy when there are topological changes. Our DG based variational  model addresses  many of the aforementioned imperfections of implicit solvent models. For example, by parametrizing both polar and nonpolar components of the solvation energy using the geometry of the interface, these two components can be coupled naturally in a single free energy functional. Application of the variational principle and the equilibrium solution of the associated Laplace-Beltrami flow  gives rise to an optimal biomolecular surface along with an optimized solvation energy.


\subsection{Polar solvation free energy}

We start with the definition of polar solvation energy, which is associated with the energy difference for charging  biomolecules in the vacuum and 
the solvent environment. Variational formulation of Poisson-Boltzmann equation was discussed in the lietrature  \cite{Sharp:1990,Gilson:1993}. Here we recast this formulation in our DG based formalism. Considering a solvated biomolecular system occupying a three-dimensional (3D) domain $\Omega \in \mathbb{R}^3$, one can relate the polar solvation energy
of the biomolecule to the electrostatic potential $\Phi({\bf r}): \mathbb{R}^3\rightarrow \mathbb{R} $ by the formulation \cite{WeiG2010a,ChenZ2010a}
\begin{equation} \label{eqn:Gp}
G_p = \int_{\Omega} \left \{ S\left [ \rho_m \Phi - \frac{1}{2} \epsilon_m | \nabla \Phi|^2 \right ] -  
(1-S) \left [ \frac{1}{2} \epsilon_s | \nabla \Phi|^2 + k_BT \sum_{i=1}^{N_c} c_i (e^{-q_i \Phi/K_BT} - 1)   \right ]  \right \} d  \vecr,
\end{equation}
where $S(\vecr)$ and $1-S(\vecr)$ are respectively the domain indicators for the solute and the solvent domains. We set    $0\leq S({\bf r})\leq 1$, which is  related to   the widely used  phase-field function $|\bar{\phi} ({\bf r})|\leq 1$ by  
\begin{equation} \label{eqn:SmSs}
S = \frac{1+\bar{\phi}}{2}, \quad  1-S= \frac{1-\bar{\phi}}{2}.
\end{equation}
Here $S$ and $1-S$ are introduced to distinguish the contributions to the total free energy from the solute region $\Omega_m$ and solvent region $\Omega_s$. The 
dielectric permittivity in these two complementary subdomains of $\Omega$ are given by $\epsilon_m$ and $\epsilon_s$, respectively. The  fixed charge density $\rho_m$ of biomolecule consists of a summation of  partial charges ($Q_j$) from atoms
\begin{equation} \label{eqn:fixed_charge_atom}
\rho_m(\vecr) = \sum_{j}Q_j \delta(\vecr - \vecr_j),
\end{equation}
where  $\vecr_j\in \mathbb{R}^3$ is the position of $j$th charged atom. In Eq. (\ref{eqn:Gp}),  $q_i$  and $c_i$ are respectively the  charge and bulk concentration of the $i$th ion species, $N_c$ is the number of ions species in the solvent,  $k_B$ is the Boltzmann constant, and $T$ is the temperature. 

\begin{figure}[!ht]
\begin{center}
\includegraphics[width=5.5cm]{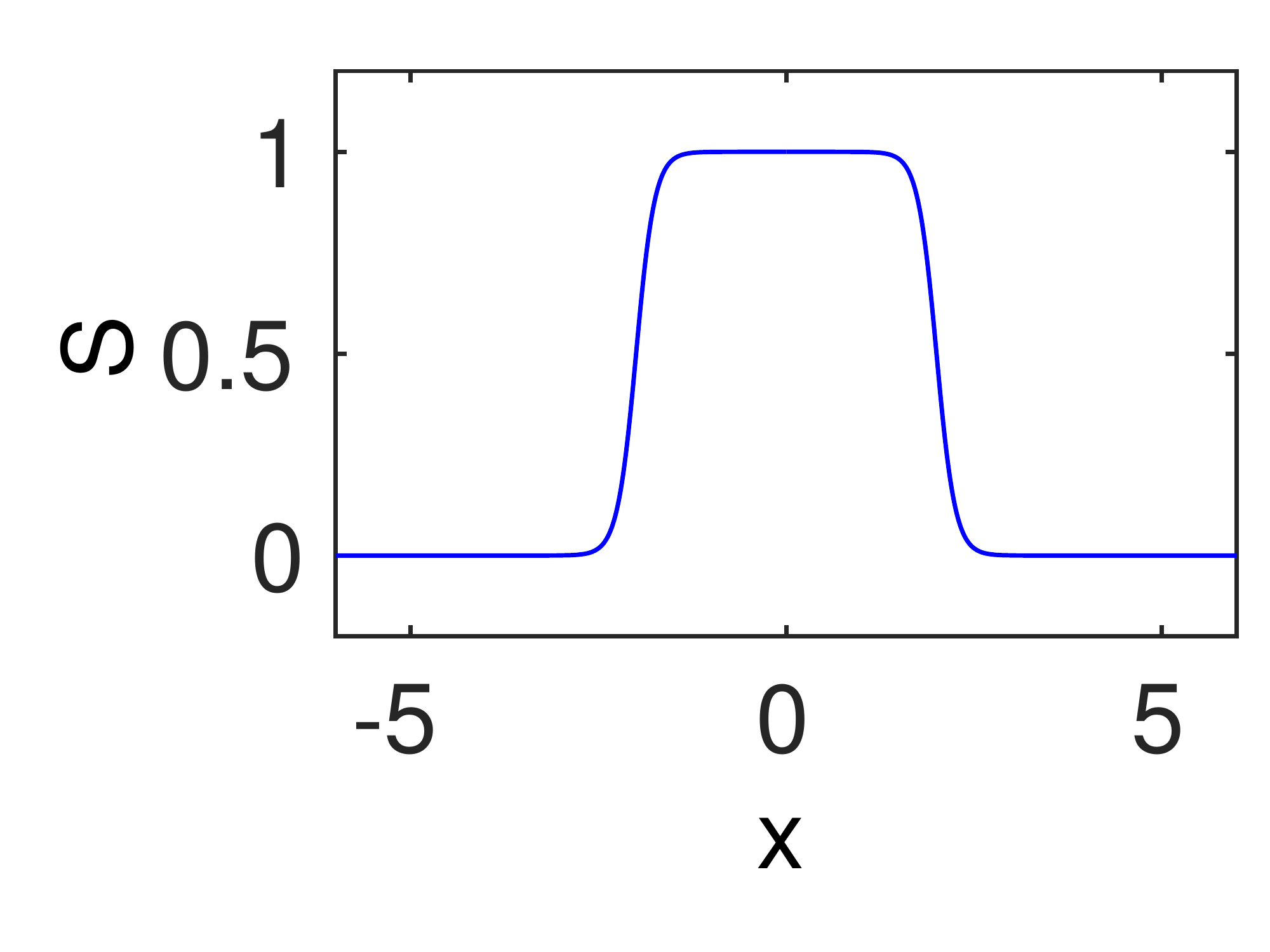}
\includegraphics[width=5.5cm]{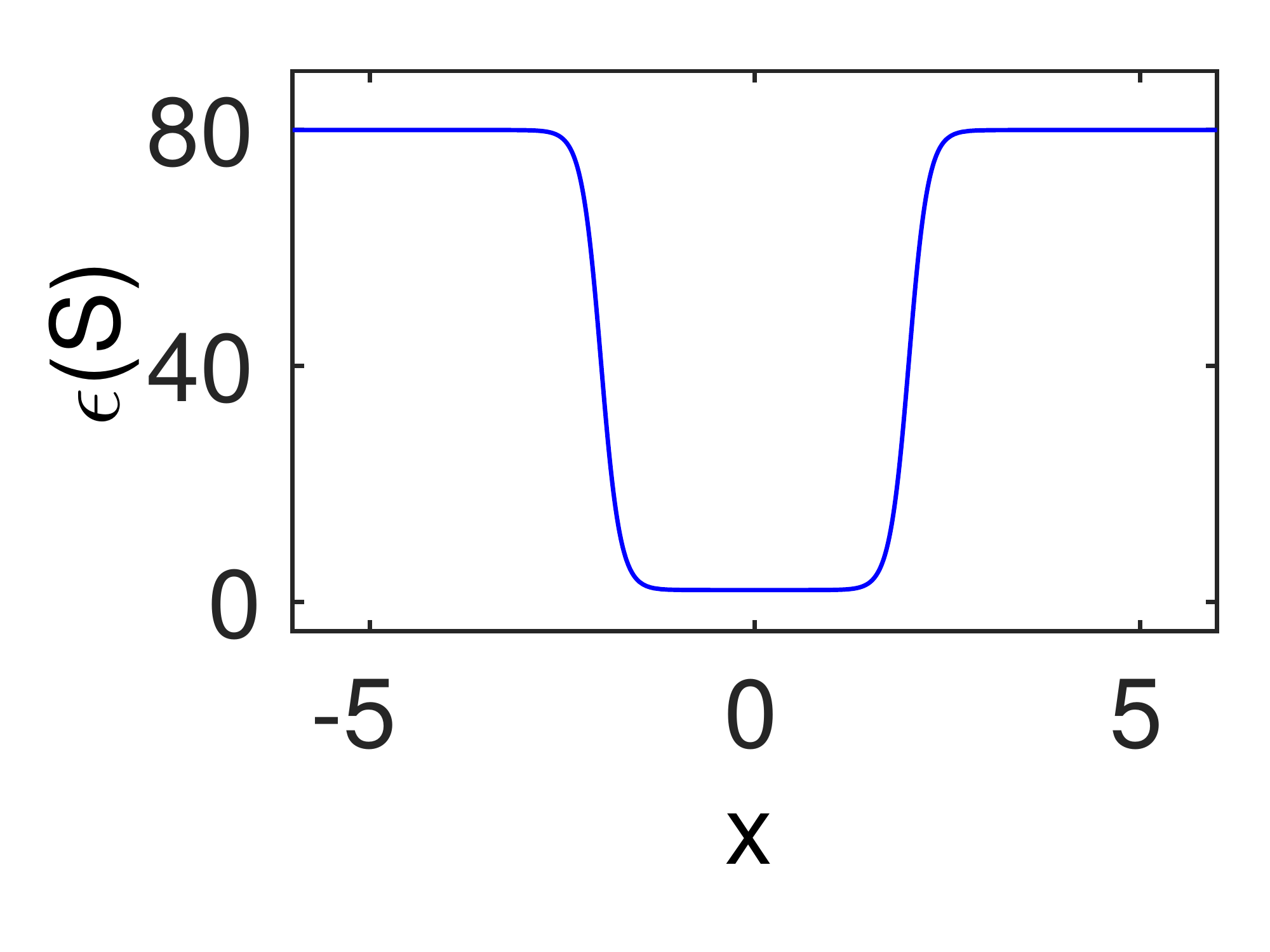}
\end{center}
\captionsetup{width=0.85\textwidth}
\caption{Left: A typical phase field function $S$ changes smoothly from its value of $-1$ in the solvent domain to
the value of $1$ in the solute domain. Right: The dielectric constant $\epsilon(S)$ depends on the phase field function 
and changes smoothly from a value of $78$ (or 80) in the solvent domain to a value of $2$ (or 1) in the solute domain.}
\label{fig:phi_and_S}
\end{figure}

The surface function $S(\vecr)$ can be chosen initially as a smooth function to simplify the numerical implementation, as seen in the left chart of
Fig.~\ref{fig:phi_and_S}. We show below the classical Poisson-Boltzmann equation can be reproduced by using this energy functional 
when a sharp solvent-solute interface is adopted, i.e., when $S$ becomes a Heaviside function. In the sequel we shall work on a generalized 
Poisson-Boltzmann equation in the sense that the transition from the solvent region to the solute region is smooth rather than 
discontinuous.

\subsection{Nonpolar solvation free energy}

The nonpolar solvation energy involves a number of terms. The scaled-particle theory (SPT) for nonpolar solutes in aqueous solutions \cite{Stillinger:1973,Pierotti:1976} utilize a solvent-accessible surface area  term \cite{Swanson:2004,massova2000combined}.  Solvent-accessible volume was shown to be relevant in large length scale regimes \cite{Lum:1999,huang2000temperature}. It was pointed out that van der Waals (vdW) interactions near solvent-solute interface are important as well \cite{Gallicchio:2002,Gallicchio:2004,choudhury2005mechanism,Wagoner06}. 
Dzubiella {\it et al}  convert these terms into a nonpolar energy functional, which, however involves a Gaussian curvature term \cite{DzubiellaJ2006PRL}.  We modify this  functional in spirit of our MMS  \cite{Bates:2006,Bates:2008} to give the following nonpolar term \cite{WeiG2010a,ChenZ2010a}
\begin{equation}  \label{eqn:Gnp}
G_{np} = \gamma A_m + p V_m + \rho_0 \int_{\Omega_s} U^{\rm att} d \vecr.
\end{equation}
Here the first term on the right is the surface energy given by the surface tension $\gamma$ and the biomolecule's surface area $A_m$.
This term measures the disruption of inter- and intra-molecular noncovalent bonds of solvent molecules when an internal surface is created. In our approach, the surface tension  $\gamma$ does not depend on Gaussian curvature so that the first term in Eq. (\ref{eqn:Gnp})
avoids possible energy jumps suggested by the Gauss-Bonnet theorem. Additionally, such a term  follows our  minimum surface energy functional formulation \cite{Bates:2006,Bates:2008}.
The   second term represents the mechanical work for expanding a volume of $V_m$ in the solvent against a hydrostatic pressure $p$. The last term
quantifies the attractive dispersion effects near the solvent-solute interface, determined by the solvent bulk density $\rho_0$ 
and the attractive portion of the van der Waals potential $U^{\rm att}$ at position $\vecr$. Since the biomolecular surface is not
explicitly known in the present modeling, we relate the surface area and its enclosed volume to the surface function $S$ through 
\begin{equation}
V_m = \int_{\Omega_m} d \vecr = \int_{\Omega} S d \vecr
\end{equation}
and the coarea formula \cite{GeometricIntegrationTheory,WeiG2010a}
\begin{equation}
A_m = \int_{\Omega} | \nabla S | d \vecr.
\end{equation}
With these relations we can assemble the polar and nonpolar contributions to give the formulation of the total solvation 
free energy functional for biomolecules at equilibrium   \cite{WeiG2010a,ChenZ2010a}
\begin{eqnarray} \label{eqn:Gtot}
G_{\rm tot} & = & \int_{\Omega}  \left \{  \gamma | \nabla  S | + p S + (1-S)\rho_0 U^{\rm att} + S \left [ (\rho_m \Phi) - 
\frac{1}{2} \epsilon_m | \nabla \Phi|^2 \right ] + \right  . \nonumber \\
& & \left . (1-S) \left [ -\frac{1}{2} \epsilon_s | \nabla \Phi|^2 - k_BT \sum_{i=1}^{N_c} c_i (e^{-q_i \Phi/K_BT} - 1)   \right ]  \right \} d  \vecr.
\end{eqnarray}
There are a variety of definitions of nonpolar free energies alternative to that in Eq. (\ref{eqn:Gnp}), but most of them are determined by the surface area, its enclosed volume and ver der Waals term in a similar way \cite{LevyR2003a,Wagoner06,LeeM2013a}. The present formulation and the variational principle introduced here are applicable to these alternative nonpolar solvation models as well.

\subsection{Governing equations}

We search for the critical point of the free energy functional to obtain the optimal free energy of the biomolecular systems. By construction, the 
free energy functional is determined by the surface function $S$ and the potential $\Phi$. The latter indeed depends on the position of 
dielectric interface hence on the surface function $S$ as well. Since the electrostatic potential follows the Poisson equation, it is
theoretically possible to replace the electrostatic potential using the convolution of the Green's function with the change density. However,
the dependence of this Green's function on the surface function $S$ does not have an explicit representation. Consequently,  it is practically 
impossible to represent the total energy as the functional of the surface function only and compute its variation. In our investigations we shall
compute the critical point by evolving the gradient flow of the free energy functional   to a steady state; while the electrostatic
potential defined by the vanishing variation $\dd{\frac{\delta G_{\rm tot}}{\delta \Phi}}$ is used as a constraint during the evolution. 
These two variations are
\begin{eqnarray}
\frac{\delta G_{\rm tot}}{\delta \Phi} & = & S \rho_m + \nabla \cdot ( (1-S) \epsilon_s + S \epsilon_m) \nabla \Phi) + (1-S) 
\sum_{i=1}^{N_c} c_i q_i e^{-q_i \Phi/K_BT},  \label{eqn:variation_phi} \\
\frac{\delta G_{\rm tot}}{\delta S} & = & -\nabla \cdot \left ( \gamma \frac{\nabla S}{|\nabla S|}\right) + p - \rho_0 U^{\rm att} + \rho_m \Phi 
+\frac{1}{2} (\epsilon_s-\epsilon_m) |\nabla \Phi|^2  \nonumber \\
& & + k_BT \sum_{i=1}^{N_c} c_i (e^{-q_i \Phi/K_BT} - 1).  \label{eqn:variation_S}
\end{eqnarray}
The vanishing variation in Eq. (\ref{eqn:variation_phi}) gives rise to a generalized Poisson-Boltzmann equation (GPBE) \cite{WeiG2010a,ChenZ2010a}
\begin{equation} \label{eqn:GPBE}
-\nabla \cdot (\epsilon(S) \nabla \Phi) = S \rho_m + (1-S) \sum_{i=1}^{N_c} c_i q_i e^{-q_i \Phi/K_BT}.
\end{equation}
where the dielectric function
\begin{eqnarray} \label{eqn:epsilon_S}
 \epsilon(S) = (1-S) \epsilon_s + S \epsilon_m,
\end{eqnarray}
is also plotted in the right chart in Fig.~\ref{fig:phi_and_S}. The gradient flow for the surface function $S$ follows 
the following a generalized Laplace-Beltrami equation \cite{WeiG2010a,ChenZ2010a}
\begin{equation} \label{eqn:SE}
\f{\p S}{\p t} =-|\nabla S| \frac{\delta G_{\rm tot}}{\delta S} = |\nabla S| \left[ \nabla \cdot \left( \gamma \f{\nabla S}{|\nabla S|}  \right ) + V\right],
\end{equation}
where a generalized potential function $V$ collects the relevant terms in Eq. (\ref{eqn:variation_S}) as
\begin{equation} \label{eqn:V}
V =-p + \rho_0 U^{\rm att} - \rho_m \Phi +\frac{1}{2} (\epsilon_m -\epsilon_s)|\nabla \Phi|^2 - k_BT \sum_{i=1}^{N_c} c_i (e^{-q_i \Phi/K_BT} - 1),
\end{equation}
and $|\nabla S|$ is added to the front of the variation to introduce the local curvature of the molecular surface to 
adjust the rate at which the surface function evolves toward its steady configuration. In this sense Eq. (\ref{eqn:SE}) is
a generalized geometric flow equation. Note that the time in Eq. (\ref{eqn:SE}) is artificial.  

We expect that the GPBE with smooth $S$ converges to its sharp interface limit when $S$ becomes a Heaviside function with a discontinuity located at 
the dielectric interface $\Gamma$, in that case the GPBE can be written as the following two elliptic equations 
\begin{eqnarray}
\quad &  -\epsilon_m \nabla^2 \Phi_m = \rho_m, &  \vecr \in \Omega_m,  \\
\quad & -\epsilon_s \nabla^2 \Phi_s = \dd{ \sum_{i=1}^{N_c} c_i q_i e^{-q_i \Phi_s/K_BT}}, &  \vecr \in \Omega_s. 
\end{eqnarray}
These two equations are coupled through the interface conditions on $\Gamma$. In this case, to make the above two equations well posed, one has to introduce two interface jump conditions,  
\begin{eqnarray}
\Phi_s = \Phi_m, \quad \epsilon_m \nabla \Phi_m \cdot \bar{\vec{n}} = \epsilon_s \nabla \Phi_s \cdot  \bar{\vec{n}},\quad  \vecr \in \Gamma
\end{eqnarray}
where $\Phi_m, \Phi_s$ are the limit values of the electrostatic potential from  solution domains $\Omega_m$ and $\Omega_s$, respectively, and $\bar{\vec{n}}(\vecr)$ is the unit normal vector on $\Gamma$.  

\subsection{Computational simulations and summary}

\begin{figure}[!ht]
\begin{center}
\includegraphics[width=11cm]{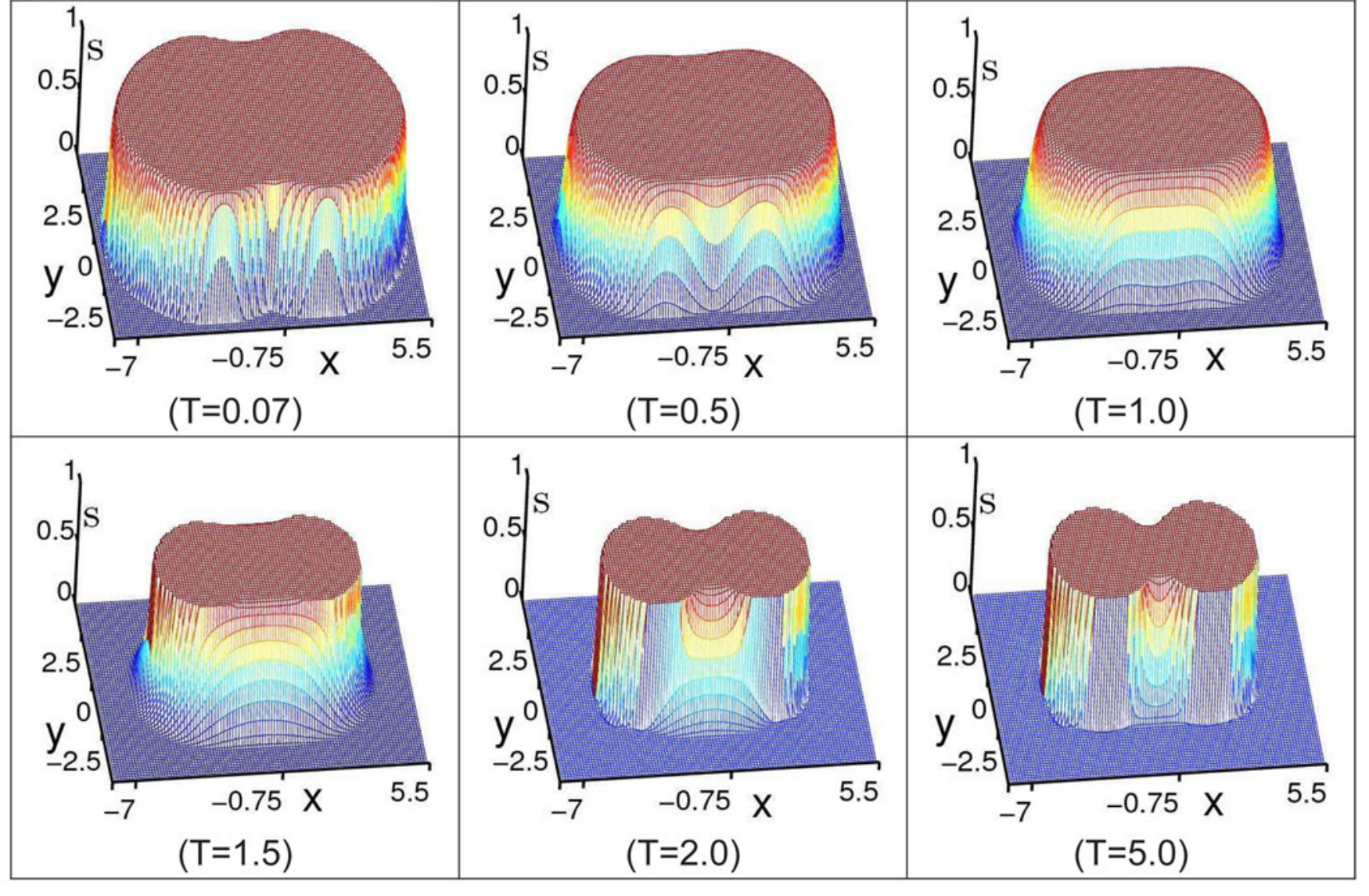}
\end{center}
\captionsetup{width=0.85\textwidth}
\caption{The phase field function evolves from its initial configuration to the final state where the surface $S=0.0$ 
fits the molecular surface for a diatomic system. Here we show only the profiles of $S$ at the cross 
section $(x,y,0.05)$ sampled at six moments during the evolution.}
\label{fig:solvation_diatomic}
\end{figure}

\begin{figure}[!ht]
\begin{center}
\includegraphics[width=11cm]{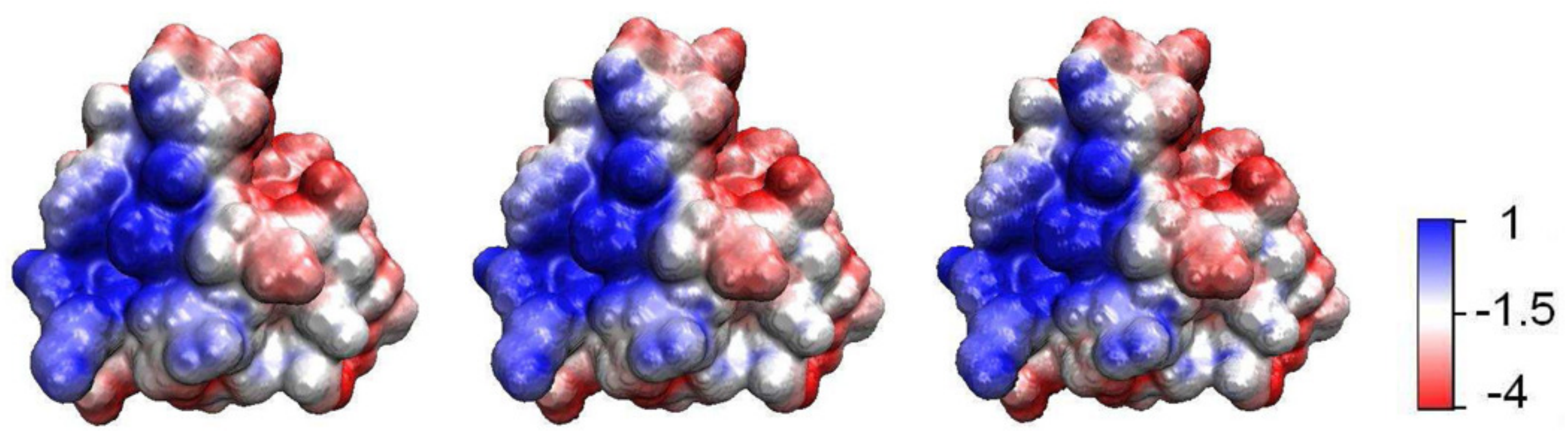}
\end{center}
\captionsetup{width=0.85\textwidth}
\caption{Electrostatic potential on molecular surfaces with different values of $S$. Left: $S=0.25$; Middle: $S=0.5$; Right: $S=0.75$.} 
\label{fig:solvation_1frd}
\end{figure}

A second-order finite difference scheme was designed to solve the coupled generalized Poisson-Boltzmann equation (\ref{eqn:GPBE}) 
and the Laplace-Beltrami equation (\ref{eqn:SE}). Most of physical parameters involved in Eq. (\ref{eqn:SE}) are taken from the
references \cite{LevyR2003a,NichollsA2008a} and the CHARMM force field. A constant surface tension $\gamma$ is chosen in our investigation whose value shall 
vary for different molecular surfaces \cite{LevyR2003a,NichollsA2008a}. In particular, $\gamma$ is implemented as a fitting parameter so that the 
optimized solvation free energy $\Delta G$ from our computational studies can match the experimental measurements. By definition,
\begin{equation} \label{eqn:dG}
\Delta G = G_{\rm tot} - G_0, 
\end{equation} 
where $G_{\rm tot}$ is defined in Eq. (\ref{eqn:Gtot}) and $G_0$ is the total energy of the solvent molecules in vacuum 
with $\epsilon_s = \epsilon_m=1$ and without nonpolar energy. To facilitate the fitting of $\gamma$ we rewrite Eq. (\ref{eqn:SE}) as
\begin{equation}
\f{\p S}{\p t} = \gamma |\nabla S| \left[ \nabla \cdot \left( \f{\nabla S}{|\nabla S|}  \right ) + \frac{V}{\gamma} \right].
\end{equation}
More details of the numerical methods for solving the coupled partial differential equations can be found in \cite{ChenZ2010a}. In Fig.~\ref{fig:solvation_diatomic} we show a simulation where the initial surface function is set such that the target  diatomic system is well contained in the region $S=1$. The surface function evolves from the initial profile toward  the final configuration that fits the molecular surface of a diatomic system, reaching a state where the total solvation energy is optimized. A more realistic simulation on the protein (PDB ID: 1frd) is shown in Fig.~\ref{fig:solvation_1frd}, where isosurfaces defined by different $S$ are plotted along with the electrostatic potential $\Phi$ on the surface. While $S=\frac{1}{2}$ is usually chosen as the molecular surface, the three surfaces are very close due to the high resolution of the numerical method. The availability of the surface position and surface potential could significantly facilitate the analysis of binding affinity  of protein-protein or protein-ligand systems, of which the electrostatic potential is an important  component \cite{Simonson02a,Archontis01,Grochowski2008a,Chu07,LinH2006a,ReddyM2014a}.

Numerically, this model can be computed by using  both  the Eulerian formulation, in which the solute boundary is embedded in the   3D  Euclidean space so evaluation of the electrostatic potential can be carried out directly \cite{ChenZ2010a}, and the Lagrangian formulation, wherein the solvent-solute interface is extracted as a sharp surface and subsequently used in solving the GPB equation for the electrostatic potential \cite{ChenZ2011a}. 
Lagrangian formulation requires direct tracking of the sampling points on the  molecular surface, which is convenient for the surface visualization, the mapping of the surface electrostatic potential field, and the enforcement of the van der Waals radii in constraint. However, it suffers  from the development of singularities while evolving molecular surface and the difficulty of handling the change of topology. In contrast, the Eulerian representation gets around of the explicit tracking of sampling points by modeling the solvent-solute interface either a smooth 3D density profile or as a specific level set of  the smooth profile.  The dynamics of the solvent-solute interface can be obtained by evolving this 3D density profile following the Laplace-Beltrami flow of the energy functional. The Eulerian representation is therefore capable of reproducing complicated dynamics of surface topology. As we shall introduce below, it also greatly  facilitates the computation of a number of geometric quantities that are otherwise difficult to compute in the Lagrangian representation, such as  the area of entire surface and  surface enclosed volume.

The parametrization of solvation energy using the surface function $S$ allows one to  track the molecular surface by following the isosurface $S=0.5$ during the evolution of $S$. This formulation is referred to the Eulerian formulation. Alternatively, one can explicitly define a molecular surface $\Gamma$ to separate the solvent and solute domains, and to use this surface to parametrize the solvation energy. Denote such an energy functional as $G_{\rm tot}(\Gamma)$. Similar to the optimization procedure presented above, the total energy is optimized by evolving $\Gamma$ following the gradient flow of the energy, and in this case, the energy variation is with respect to the spatial coordinates of this explicitly defined surface $\Gamma$. Numerically, this can be achieved by discretizing $\Gamma$ into a collection of surface elements or surface vectors $\{ \hat{S}_j \} $, each element parametrized by a local coordinate system $(x_1, x_2)$, and thus $G_{\rm tot}(\Gamma)$ becomes $G_{\rm tot}( \hat{S}_j)$. Furthermore, we can constrain the motion of $\Gamma$ to the normal direction $\vec{n}(x_1,x_2)$ only, for that a tangential displacement of $\Gamma$ does not change the surface configuration and the solvation energy. A scalar displacement field $\psi(x_1,x_2)$ in the normal direction can be defined through
\begin{equation}
 \hat{S}_j^{\sigma}(x_1,x_2)  =  \hat{S}(x_1,x_2) + \sigma \psi(x_1,x_2) \vec{n}(x_1,x_2),
\end{equation} 
which states that the surface element $ \hat{S}_j$ is updated from its original position by $\sigma \psi(x_1,x_2)$ along the normal direction to the new position $ \hat{S}_j^{\sigma}$, where $\sigma$ is a number to scale the normal displacement field $\psi(x_1,x_2)$. The optimization of the total energy at a particular molecular surface $\Gamma$ means that any normal displacement will violate the nature of optimum at this point, indicating
\begin{equation}
\left . \frac{\p  \hat{S}_j^{\sigma}}{\p \sigma} \right |_{\sigma = 0} = 0.
\end{equation}
Now we can observe the transition of the independent variables in calculating the energy variation:
\begin{equation}
\f{\delta }{\delta \Gamma} \rightarrow  \f{\p }{\p  \hat{S}_j^{\sigma}} \rightarrow \f{\p }{\p \sigma},
\end{equation}
as a result of replacing the motion of the explicit surface $\Gamma$ using the scaled normal motion of a collection of surface elements. The readers
are referred to \cite{ChenZ2011a} for the detailed calculation of the energy variation, the derivation of the equation governing the gradient flow, and the 
numerical techniques for solving the equation. This investigation also shows that the optimized solvation energy and molecular surface are well matching 
those generated by the Eulerian formulation if there is no topological change in $\Gamma$ during its evolution. Notice that a single point on $\hat{S}_j$ may evolves to two distinct points, or two distinct points in two different surface elements may converge to a single point when there is a topological change during the 
evolution of $\Gamma$. This intrinsic singularity in handling the topological change limits the applications of the Lagrangian formulation to complex 
biomolecular systems, for which it is impossible to set an initial surface $\Gamma$ that is topologically equivalent to the final optimized molecular surface. 
The Eulerian formulation is hence suggested for the investigations of the solvation energy and molecular surfaces of general biomolecular systems.

Recently,  differential geometry based implicit solvent model has been tested extensively via solvation analysis \cite{ChenZ2010a,ChenZ2011a,ChenZ2012a,Daily:2013, Thomas:2013,BaoWang:2015a}. The differential geometry based nonpolar model was found to deliver some of the best  nonpolar  solvation predictions  \cite{ChenZ2012a}. However, for general molecules with a significant polar component, our initial predictions were not up to the state of the art \cite{ChenZ2010a,ChenZ2011a}. It turns out that both the generalized Laplace-Beltrami equation and the generalized Poisson-Boltzmann equation can be  easily solved individually. However,  when these equations are coupled, there is a stability problem \cite{SZhao:2011a,SZhao:2014a}.  Essentially, when $S$ admits unphysical values beyond its physical definition $0\leq S\leq1$, the dielectric function  (\ref{eqn:epsilon_S}) will adopt unphysical (negative) values as well, which gives rise to an instability in updating   the Laplace-Beltrami equation (\ref{eqn:SE}).  This issue hinders the performance of DG based solvation models for molecules with significant polar component. To address this problem, a convex optimization algorithm \cite{BaoWang:2015a} has been developed to ensure the stability in solving coupled PDEs (\ref{eqn:GPBE}) and (\ref{eqn:SE}). As a result, the differential geometry based solvation model is found to deliver some of the most accurate prediction of experimental solvation free energies for more than 100 molecules of both polar and nonpolar types \cite{BaoWang:2015a}.

 
Most recently, Wei and coworkers have taken a different treatment of non-electrostatic interactions between the solvent  and solute in the DG based solvation models so that the resulting total energy functional and PB equations are consistent with more detailed descriptions of solvent densities at equilibrium \cite{WeiG2012a,Wei:2013}. To account for solute response to solvent polarization, a quantum mechanical (QM) treatment of solute charges was introduced to the DG-based solvation models using the Kohn-Sham density functional theory (DFT) \cite{ZhanChen:2011a}. This multiscale approach self-consistently computes the solute charge density distribution which simultaneously minimizes both the DFT energy as well as the solvation energy contributions.

Currently, efforts are  invested to improve the accuracy and robustness of DG based solvation models by combining physical models with knowledge based models, namely, machine learning approaches. Additionally, DG based solvation models and machine learning approaches are utilized for accurate  predictions of   the protein binding energies and ligand binding affinities over a wide range of conformational states. Furthermore,  it is worth noting that the method depends only on the representation of the  solvent-solute interfaces, and this representation is independent of the atomic or coarse-grained description of the biomolecules. It is  therefore possible to adopt this method to compute the potential of mean force of coarse-grained biomolecular structures along selected coordinate, and the results can be utilized for parametrization the force field for coarse-grained molecular systems as well. Finally, we would like to point out that  many    critical applications to biophysics, chemistry, and medicine mostly remain unexplored.

\section{Variational Methods for Pattern Formation in Bilayer Membranes} \label{sec:surfpattern}

As one of the most important biomolecular systems, the lipid bilayer membranes sustain the regular functions of cell and subcelluar compartments by regulating the transmembrane ion or molecular flows and by providing platforms for various essential biochemical processes \cite{BiologicalMembranes2002,Molecularbiologyofthecell}. These critical functions of bilayer membranes are determined by their lipid compositions, the specific membrane proteins, and their dynamical arrangement in the bilayers during the course of membrane morphology change as a result of various membrane-solvent, membrane-membrane, or membrane-protein interactions. Applications of the variational principle for bilayer membrane modeling have been mostly focused on four types of problems: (i) mean-curvature dependent membrane morphology \cite{DuQ2005d,DuQ2006a,McMahonH2005,CookeI2006a}, (ii) ionic or proton flows in protein channels \cite{XuS2014a,ZhengQ2011a}, (iii) lateral diffusion on membrane surfaces \cite{ZhouY2012b}, and (iv) pattern formation in bilayer membranes \cite{DuQ2011a,CamleyB2010a,WitkowskiT2012a}. Here in this section we focus  on the local pattern formation in bilayer membranes, for that there are many controversial investigations concerning the biophysical underpinning of these patterns, their spatial and temporal distributions, and their roles in modulating relevant biochemical processes \cite{TakahashiT2011a,PikeL2003a,AndersonR2002a,VieiraF2010a}. These patterns are called lipid rafts, which are small (10-200nm), heterogeneous, highly dynamic, sterol- and sphingolipid-enriched domains that compartmentalize cellular processes \cite{SimonK2000a}. Lipids move laterally within the domains mostly rather than over the entire membrane surface \cite{ApajalahtiT2010a}. Classical phase  separation models manage to minimize the total area of the domain boundaries and large domains appear at the end of the minimization; this process is usually referred to as coarsening. When these classical models are directly extended to model surface phase separation, the total arc length of the domain boundaries on the surface is minimized to generate large domains, which do not match the measured sizes of    lipid rafts \cite{DuQ2011a,CamleyB2010a,WitkowskiT2012a}. 

\subsection{Classical phase field models}
We first examine the classical phase separation model for binary systems. Consider two species of particles in $\mathbb{R}^3$ with respective mass or volume fractions 
$m_1,m_2 \in [0, 1]$. The interactions between particles of the same species are favorable while the interactions between different 
species are unfavorable. This preference can be modeled by defining a phase field function 
\begin{equation}
\bar{\phi} = \frac{m_1 - m_2}{m_1 + m_2},
\end{equation}
where { $\bar{\phi}({\bf r}) \in [-1,1]$}, ${\bf r} \in \mathbb{R}^3$ and minimizing the Ginzburg-Landau free energy functional in $\Omega\in \mathbb{R}^3$
\begin{equation}
G(\bar{\phi}) = \int_{\Omega} \left ( f(\bar{\phi}) + \frac{\sigma}{2} | \nabla \bar{\phi}|^2 \right ) d{ {\bf r}},
\end{equation}
where $f(\bar{\phi})$ is a double well potential that has two minimums at $\bar{\phi}=\pm 1$. A typical choice is
\begin{equation}
f(\bar{\phi}) = \frac{\bar{\phi}^4}{4} - \frac{\bar{\phi}^2}{2}
\end{equation}
which has two symmetric potential wells of the same depth at $\bar{\phi} = \pm 1$. It is apparent that a complete phase separation 
with $\bar{\phi}$ changing discontinuously between $1$ and $-1$ is favorable by $f(\bar{\phi})$ when $G(\bar{\phi})$ is minimized. Such an unphysical distribution of $\bar{\phi}$ is to be penalized by the term $\dd{\frac{\sigma}{2}|\nabla \bar{\phi}|^2}$ that regulates the transitional gradient of $\bar{\phi}$ between $1$ and $-1$.

\subsection{Geodesic curvature based membrane models}

\begin{figure}[!ht]
\begin{center}
\includegraphics[height=2.5cm]{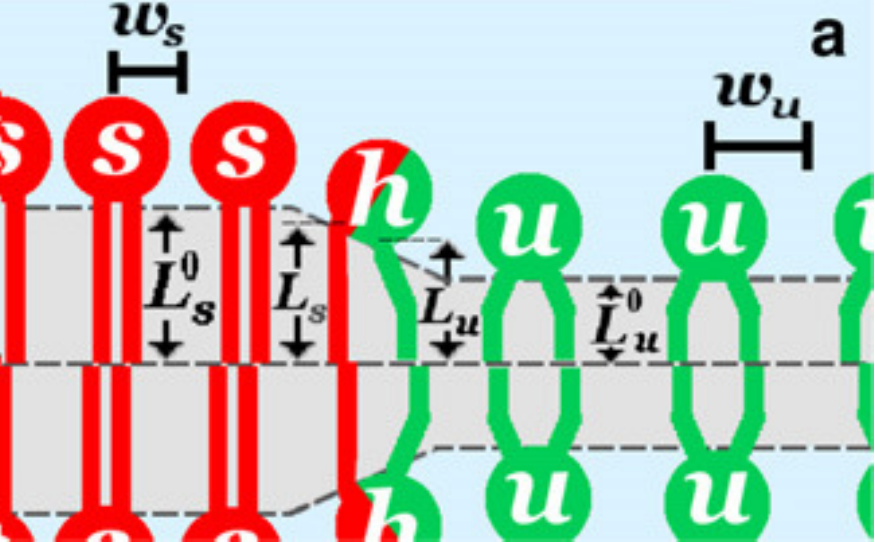} \hspace{3mm}
\includegraphics[height=2.5cm]{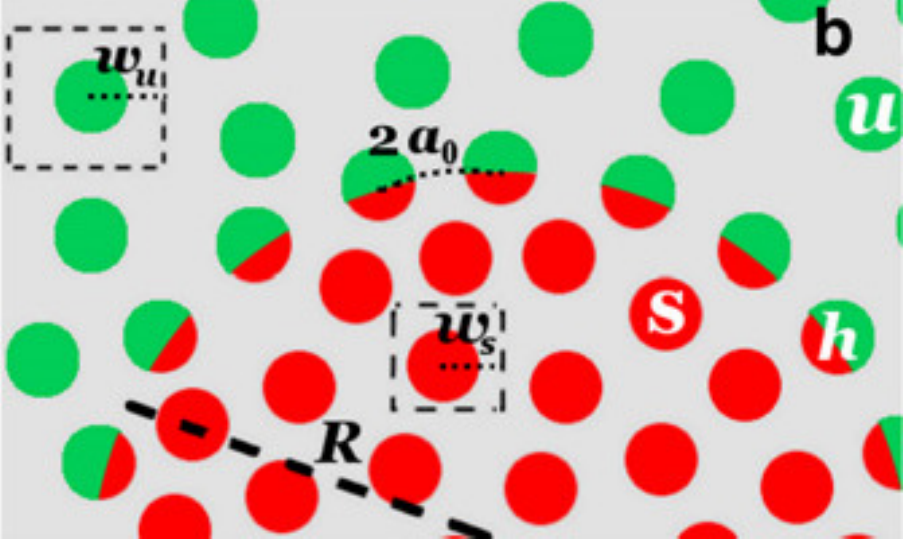} \hspace{3mm}
\includegraphics[height=2.5cm]{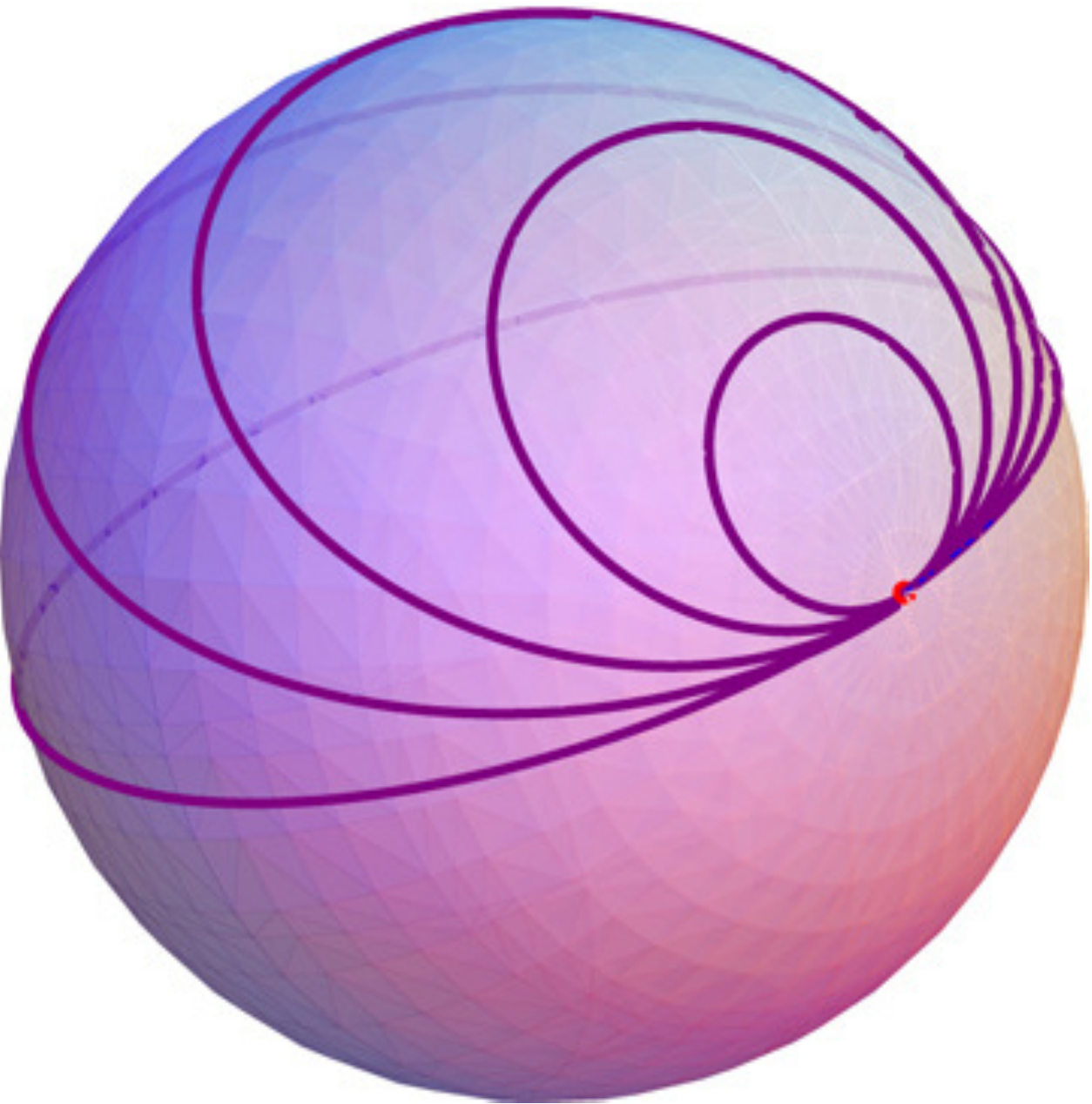}
\end{center}
\captionsetup{width=0.85\textwidth}
\caption{Left: Schematic illustration of the mismatch of the lipid structures at the interface that induces a transitional hybrid region 
between two lipid domains \cite{BrewsterR2010a}. Middle: Within the transitional hybrid layer the otherwise regular lattices of the lipids in either 
domain relax to match each other, causing a bending interface \cite{BrewsterR2010a}. Right: Circles on a sphere have constant geodesic curvatures.
The great circle, i.e., the lowest circle, has a vanishing geodesic curvature in particular.}
\label{fig:surface_pattern_model}
\end{figure}

\subsubsection{Lagrangian formulation}

Our variational model is motivated by the recent theoretical studies of the hybrid lipids saturation at the interface between saturated and unsaturated of 
lipids with geometrical and molecular mechanical mismatch \cite{BrewsterR2010a}. As illustrated in Fig.~\ref{fig:surface_pattern_model}, two 
species of lipids at their interface have different intermolecular interactions that are determined by their structures. The otherwise regular 
lattice of either species of lipids has to be relaxed in a way such that the intermolecular interactions in the transitional  region 
near the interface will fit the different lattice structure of other species. This relaxation generates curved interface between two species of 
lipids in a manner similar to the generation of surface tension. Since the domain boundary is a curve on a { two-dimensional (2D) surface 
embedded in $\mathbb{R}^3$}, it is the geodesic curvature of the interface rather than the interface curvature that determines the intermolecular interactions between two species of lipids near the 
interface. \footnote{In Sect.~\ref{sec:solvation}, we use $S$ to denote the surface function, which is a domain indicator,  and use $\Phi$ to denote the electrostatic potential following the traditional usage in the studies of biomolecular electrostatics. Here in Sect.~\ref{sec:surfpattern} and Sect.~\ref{sec:protein} the models do not involve electrostatics, and we  denote $\phi$ the phase field function, while use $S$ to denote the 2D {  surface  embedded in $\mathbb{R}^3$} when applicable.   An interface in Sect.~\ref{sec:solvation} refers to solvent-solute boundary region, whereas in Sect.~\ref{sec:surfpattern} and Sect.~\ref{sec:protein}, it refers a boundary curve on a given surface.  } 
The geodesic curvature of the interface measures how far the interface curve is from being a geodesic. We define the curvature energy of the microdomain boundary by a  {one-dimensional (1D)on-curve integration}  
\begin{equation} \label{eqn:geodesic_eng}
G = \int_{C} k(H  - H_0)^2 ds,
\end{equation} 
where $C$ is the domain boundary contour embedded in $\mathbb{R}^3$,  $H$ is the geodesic curvature, $H_0$ is the spontaneous geodesic curvature of the lipid mixture to be separated, and $k$ is the geodesic curvature energy coefficient. The spontaneous geodesic curvature $H_0$ is an intrinsic property of the combination of any two species of lipids in the bilayer membrane that will be separated to form local microdomains as a result of geometric and molecular mechanical mismatch. In the transitional region near the interface two species of lipids are arranged in a hybrid state rather than the regular lattice structure. Indeed a recent theoretical study adopted a free energy for the hybrid packing of two species of lipids (denoted by the subscript $1$ and $2$ below) at the interface \cite{BrewsterR2009a,BrewsterR2010a}:
\begin{equation} \label{eqn:bending_energy}
\mathcal{F} = k_s(L_1 - L_1^0)^2 + k_u(L_2 - L_2^0) ^2 + \gamma (L_1 - L_2)^2, 
\end{equation}
where $L_i$ is the length of the lipid chains in the transitional region and $L_i^0$ is the length of the equilibrium chain in the bulk.  Parameters $k_s $ and $  k_u$  are the free energetic costs of mismatch between two species and their hybrids at the interface, respectively and similarly, $\gamma$ is the energy cost of  mismatch between two chains of the hybrid.  Furthermore, the following relations are identified to related the domain curvature and lipid geometrical properties:
\begin{equation} \label{eqn:V_i}
V_i = L_i a_0 w_i \left(1 \pm \frac{w_i H}{2} \right), \quad i = 1, 2,
\end{equation}
where $V_i$ is the molecular volume of the lipid chains, $w_i$ is  the length that characterizes the molecular spacing of the lipid head groups, and $a_0 = (w_1 + w_2)/2$ is the headgroup spacing of the hybrids along the interface. Here the subtraction sign is chosen if the species is included in the microdomain, otherwise the addition sign is used. The chain length in the equilibrium bulk state, $L_i^0$, can be computed from the molecular volume divided by the head group area in the bulk state
\begin{equation}
L_i^0 = \frac{V_i}{w_i^2}.
\end{equation}
Eqs. (\ref{eqn:bending_energy}-\ref{eqn:V_i}) represent the interface bending energy $\mathcal{F}$ as a function of it geodesic curvature $H$. The minimizer $H_0$ can be analytically calculated to the linear order:  
\begin{equation}
H_0 = \frac{1}{w_T} \left [ \frac{(1 - 2 B)w_d }{(1+2 B) w_T} + \frac{2 B V_d }{(1+2 B) V_T}\right],
\end{equation}
where $B$ is a constant characterizing the free energetic cost of lipid mismatch at the interface,  
$w_T = (w_1 + w_2)/2, w_d = w_1 - w_2, V_T = (V_1 + V_2)/2,$ and $  V_d = V_1 - V_2$. By truncating the Taylor series approximation of
$\mathcal{F}(H)$ with respect to $H_0$ to the second order we   get an energy functional in the form of Eq. (\ref{eqn:geodesic_eng}).

\subsubsection{Eulerian formulation}

It has been seen in Sect.~\ref{sec:solvation} that the parametrization of solvation energy using the surface function allows  one to implicitly track the molecular surface by following the iso-surface extraction during the evolution of the surface function, which  is referred as to the Eulerian formulation.  We could also evolve a phase field function to minimize the energy in Eq. (\ref{eqn:geodesic_eng}) and to obtain the configuration of microdomains. This is achieved by using the following 2D Eulerian formulation of the microdomain geodesic curvature energy defined on the entire membrane surface $S$:
\begin{equation} \label{eqn:geo_curv_energy}
G(\phi) = \int_{S} \frac{k \varepsilon}{2} \left( \Delta_{\vecx} \phi + \frac{1}{\varepsilon^2} (\phi + H_c \varepsilon) (1 - \phi^2)\right)^2 d { \vecx}
\end{equation} 
where $H_c = \sqrt{2} H_0$ and   $\varepsilon$ is a small positive parameter that characterizes the width of the transitional layer from $\phi(\vecx)=-1$ to $\phi(\vecx)=1$. Here $S$ is a surface embedded in ${\mathbb R}^3$, $\vecx=(x_1,x_2)$  and $d{\bf x}$  is an infinitesimal surface element.  The equivalence of this Eulerian formulation (\ref{eqn:geo_curv_energy}) to the Lagrangian formulation (\ref{eqn:geodesic_eng}) is 
analogous to the equivalence between the Canham-Helfrich-Evans curvature energy and the membrane elastic energy \cite{DuQ2005d,Melissa_PhDThesis}. In particular, if the phase field function is defined by
\begin{equation}
\phi(\vecx) = \tanh \left ( \frac{d(\vecx)}{\sqrt{2} \varepsilon} \right )
\end{equation}
with $d(\vecx)$ being the signed geodesic distance at the surface point $\vecx$ to the interface contour $C$ where $\phi=0$, then 
$$\nabla_{\vecx} \phi = \frac{1}{\varepsilon} q'(d(\vecx)) \nabla_{\vecx}d, \quad  \Delta_{\vecx} \phi = \frac{1}{\varepsilon} q''(d(\vecx)) |\nabla_{\vecx} d |^2 + 
\frac{1}{\varepsilon} q'(d(\vecx)) \Delta_{\vecx} d,$$
where
$$q(x) = \tanh \left( \frac{{ x}}{\sqrt{2}{ \varepsilon}} \right ),  
q'(\vecx) = \frac{1}{\sqrt{2}} \left [ 1 - \tanh ^2 \left ( \frac{{ x}}{\sqrt{2}{ \varepsilon}} \right ) \right ], $$  
$$
q''(\vecx) = -{ \f{1} {\varepsilon}} \tanh \left ( \frac{{ x}}{\sqrt{2}{ \varepsilon}} \right ) \sech^2 \left ( \frac{{ x}}{\sqrt{2}{ \varepsilon}} \right ),$$
and $\nabla_{\vecx}, \nabla_{\vecx} \cdot$ are surface gradient and surface divergence operators, respectively.  The geodesic curvature of a contour is given by
\begin{equation}
H = \nabla_{\vecx} \cdot \vecn,
\end{equation}
where $\vecn$ is the normal vector to the contour $C$. Since $\vecn = \nabla_{\vecx} d$ we 
have $H = \nabla_{\vecx} \cdot \nabla_{\vecx} d = \Delta_{\vecx} d$ and
$$ \Delta_{\vecx} d = \frac{\varepsilon}{q'} \Delta_{\vecx} \phi - \frac{q''}{q'} | \nabla_{\vecx} d|^2, 
\quad \nabla_{\vecx} d = \frac{\varepsilon}{q'} \nabla_{\vecx} \phi.$$
Therefore, one has 
$$ \Delta_{\vecx} d = \frac{\varepsilon}{q'} \Delta_{\vecx} \phi - \frac{q''}{q'} \left |\frac{\varepsilon}{q'} \nabla_{\vecx} \phi \right |^2.$$
Writing $q'(\vecx)$ and $q''(\vecx)$ in terms of $q(\vecx)$ we can convert the above representation to  
$$ \Delta_{\vecx} d = \frac{\sqrt{2} \varepsilon}{1-q^2} \left( \Delta_{\vecx} \phi + \frac{2q}{1-q^2} |\nabla_{\vecx} \phi|^2 \right ),$$
which is the geodesic curvature $H = \Delta_{\vecx} d$. Replacing $q(\vecx)$ with $\phi$ one obtains the final form of $H$ as 
\begin{align}
H & = \frac{\sqrt{2} \varepsilon}{1-\phi^2} \left( \Delta_{\vecx} \phi + \frac{2\phi}{1-\phi^2} |\nabla_{\vecx} \phi|^2 \right ) \nonumber \\
  & = \frac{\sqrt{2} \varepsilon}{1-\phi^2} \left( \Delta_{\vecx} \phi + \frac{1}{\varepsilon^2} (1-\phi^2) \phi \right ), 
\end{align}
where we assume $\|\vecn\|  = 1$ in the last step of derivation. When minimizing the curvature energy in Eq. (\ref{eqn:geo_curv_energy}) the following constraint 
\begin{equation}
A(\phi) = \int_{S} \phi(\vecx) d{ \vecx} = \mbox{constant}
\end{equation}
must be enforced such that quantities of both species of lipids are conserved.

To derive the equation of the geometric flow for the energy $G(\phi)$ we compute its first variation with respect to $\phi$:
\begin{equation}
\frac{\delta G}{\delta \phi} = k \left[ \Delta_{\vecx} W - \frac{1}{\varepsilon^2} (3 \phi^2 + 2H_c \varepsilon \phi - 1) W \right]
\end{equation} 
where
$$ W = \varepsilon \Delta_{\vecx} \phi - \frac{1}{\varepsilon} ( \phi + H_c \varepsilon)(\phi^2 - 1).$$
We then split the linear and nonlinear components ($W_L $ and $ W_N$) of $W$ to facilitate the numerical treatments. They are given respectively by
$$ W_L = \varepsilon \Delta_{\vecx} \phi + \frac{1}{\varepsilon} \phi + H_c, \quad W_N = -\frac{1}{\varepsilon} \phi^3 - 
H_c \phi^2.$$
We then have the full expansion of the variation
\begin{eqnarray}
\frac{\delta G}{\delta \phi} & = & k \Delta_{\vecx} W_L + \frac{k}{\varepsilon^2} W_L + k \Delta_{\vecx} W_N  - 
\frac{k}{\varepsilon^2} (3 \phi^2 + 2H_c \varepsilon \phi) (W_N + W_L) + \frac{k}{\varepsilon^2} W_N \nonumber \\
    & = & 
	k \varepsilon \Delta_{\vecx}^2 \phi 
	+ \frac{k}{\varepsilon}\left(2- 6 \phi^2 - 4 k H_c \epsilon \right) \Delta_{\vecx} \phi 
	-   \left( \frac{6k}{\varepsilon} \phi+ 2k H_c \right) | \nabla_{\vecx} \phi|^2
	\nonumber \\ 
       & + & k \left( - \frac{2H_c^2 }{\varepsilon}+ \frac{1}{\varepsilon^3} \right) \phi 
	- \frac{ 3 k H_c}{\varepsilon^2} \phi^2 
	- k\left( \frac{4}{\varepsilon^3} - \frac{2  H_c^2}{\varepsilon}\right) \phi^3 
	+ \frac{5kH_c}{\varepsilon^2} \phi^4  
	+   \frac{3k}{\varepsilon^3} \phi^5 
	\nonumber \\
	& +& \frac{k H_c}{\varepsilon^2}.
\end{eqnarray}
 Also note that the variation of the mass conservation constraint is
\begin{equation}
\frac{\delta A}{\delta \phi} = 1.
\end{equation}
The appearance of fourth order derivative in the variation $\delta G/\delta \phi$ motivates us to adopt the following equation of the geometric flow with an artificial time for $\phi$:
\begin{equation} \label{eqn:gradient_flow_surface_pattern}
\frac{\p \phi}{\p t} =- \frac{\delta G}{\delta \phi} + \lambda \frac{\delta A}{\delta \phi},
\end{equation}
where $\lambda$ is a Lagrangian multiplier used to ensure the conservation of $\phi$. We can derive a representation of $\lambda$ by integrating Eq. (\ref{eqn:gradient_flow_surface_pattern}) and noting that $\dd{\int_S \frac{\p \phi}{\p t} d{  \vecx} = 0}$, hence
$$ 0 = -\int_S \frac{\delta G}{\delta \phi} d{  \vecx} + \int_S \lambda d{  \vecx},$$
and consequently
$$ \lambda =  \frac{1}{|S|} \int_S \frac{\delta G}{\delta \phi} d{  \vecx},$$ 
which yields
\begin{equation} \label{eqn:gradient_flow_surface_pattern_2}
\frac{\p \phi}{\p t} =- \frac{\delta G}{\delta \phi} + \frac{1}{|S|} \int_S \frac{\delta G}{\delta \phi} d{  \vecx}.
\end{equation}
Eq. (\ref{eqn:gradient_flow_surface_pattern_2}) is a fourth-order nonlinear surface diffusion equation. Alternatively, one could   derive a Cahn-Hilliard equation for the surface phase field function $\phi$ as
\begin{equation} \label{eqn:surface_CHE}
\frac{\p \phi}{\p t} = \Delta_{\vecx} \left( \frac{\delta G}{\delta \phi} \right),
\end{equation}
which guarantees the conservation of $\phi$ and thus does not need a Lagrangian multiplier. However, it involves a sixth order surface derivative and thus 
is more complicated when the equation is to be solved numerically on a discretized surface $S$. 

To simplify the exposition of numerical treatments we adopt  $\lambda = \frac{1}{|S|}\int_{S} \frac{\delta A}{\delta \phi} d{  \vecx} $ and define $g = \frac{\delta G}{\delta \phi}$. Then we write Eq. (\ref{eqn:gradient_flow_surface_pattern_2}) as
\begin{equation}
\phi_t = -g + \lambda.
\end{equation}
To implement the time discretization we average the nonlinear function $g(\phi)$ over the current and next time steps $\phi_n,\phi_{n+1}$ to implement a Crank-Nicolson approximation
\begin{equation} \label{eqn:surface_pattern_CN}
\frac{\phi_{n+1} - \phi_{n}}{\Delta t} + g(\phi_{n+1}, \phi_{n}) - \lambda(\phi_n) = 0,
\end{equation}
where the averaged function is defined by
\begin{eqnarray*}
g(\phi_{n+1}, \phi_n) & = & \frac{k}{2} \Delta_{\vecx} (f_c(\phi_{n+1}) + f_c(\phi_n)) - \\
    &  & \frac{k}{2\varepsilon^2} (\phi_{n+1}^2 + \phi_{n+1} \phi_n + \phi_n^2 + \varepsilon H_c (\phi_{n+1} + \phi_n) - 1) (f_c(\phi_{n+1}) + f_c(\phi_n)),
\end{eqnarray*} 
and 
\begin{equation*}
f_c(\phi) = k\left ( \varepsilon \Delta_{\vecx} \phi 
- (\frac{1}{\varepsilon} + \varepsilon H_c) (\phi^2 -1) \right).
\end{equation*}
 To numerically solve Eq. (\ref{eqn:surface_pattern_CN}) which is an implicit scheme for $\phi_{n+1}$, we define an interior iteration for computing $\psi_m$ such that $\psi_m \rightarrow \phi_{n+1}$ as $m \rightarrow \8$. The equation for $\psi_m$ reads as
\begin{equation}  \label{eqn:surface_pattern_CN_2}
\frac{\psi_{m+1} - \phi_{n}}{\Delta t} + g(\psi_{m+1}, \psi_m, \phi_{n}) - \lambda(\psi_m) = 0,
\end{equation}
where new averaged functions are defined by
\begin{eqnarray*}
g(\psi_{m+1}, \psi_m, \phi_n) & = & \frac{k}{2} \Delta_{\vecx} \tilde{f}_c(\psi_{m+1}, \psi_{m}, \phi_n) - \\
    &  & \frac{k}{2\varepsilon^2} (\psi_{m}^2 + \psi_{m} \phi_n + \phi_n^2 + \varepsilon H_c (\psi_{m} + \phi_n) - 1) 
(f_c(\psi_{m}) + f_c(\phi_n)), \\
\tilde{f}_c(\psi_{m+1}, \psi_m, \phi_n) & = &  \frac{\varepsilon}{2} \Delta_{\vecx}(\psi_{m+1} + \phi_n) - \frac{1}{4 \varepsilon} (\psi_m^2 + \phi_n^2 -2)
(\psi_m + \phi_n + 2\varepsilon H_c).
\end{eqnarray*}
Convergent $\psi_m$ is obtained by iterating over the interior index $m$, usually up to a tolerance $\| \psi_{m+1} - \psi_{m} \| \le \varepsilon_{\psi}$ for some small $\varepsilon_{\psi} >0$. This convergent $\psi_{m}$ is assigned to $\phi_{n+1}$, and computation is advanced to the next time step. The linear and nonlinear components of $\psi_{m+1}$ in Eq.(\ref{eqn:surface_pattern_CN_2}) are further split. The nonlinear components are updated slower than the linear components, allowing an efficient numerical solution. The spatial approximation of the equation is obtained by a newly developed a $C^0$ interior penalty surface finite element method \cite{BrennerS2012a,Melissa_PhDThesis}.

\subsection{Computational simulations and summary}

We apply the geodesic curvature driven phase separation model to simulate the microdomain formation on surfaces. We present four simulations on different surfaces or with different spontaneous geodesic curvatures. The energetic histogram and the dynamics of the domain formation in each simulation are compared to those generated by the Allen-Cahn equation obtained by the direct extension of the Ginzburg-Landau energy based a classical phase separation model on surfaces \cite{DuQ2011a}. We also compute the radii of the microdomains which are expected to approximate the reciprocal of the given spontaneous geodesic curvature.

\begin{figure}[!htb]
\begin{center}
\includegraphics[angle=270,width=3.8cm]{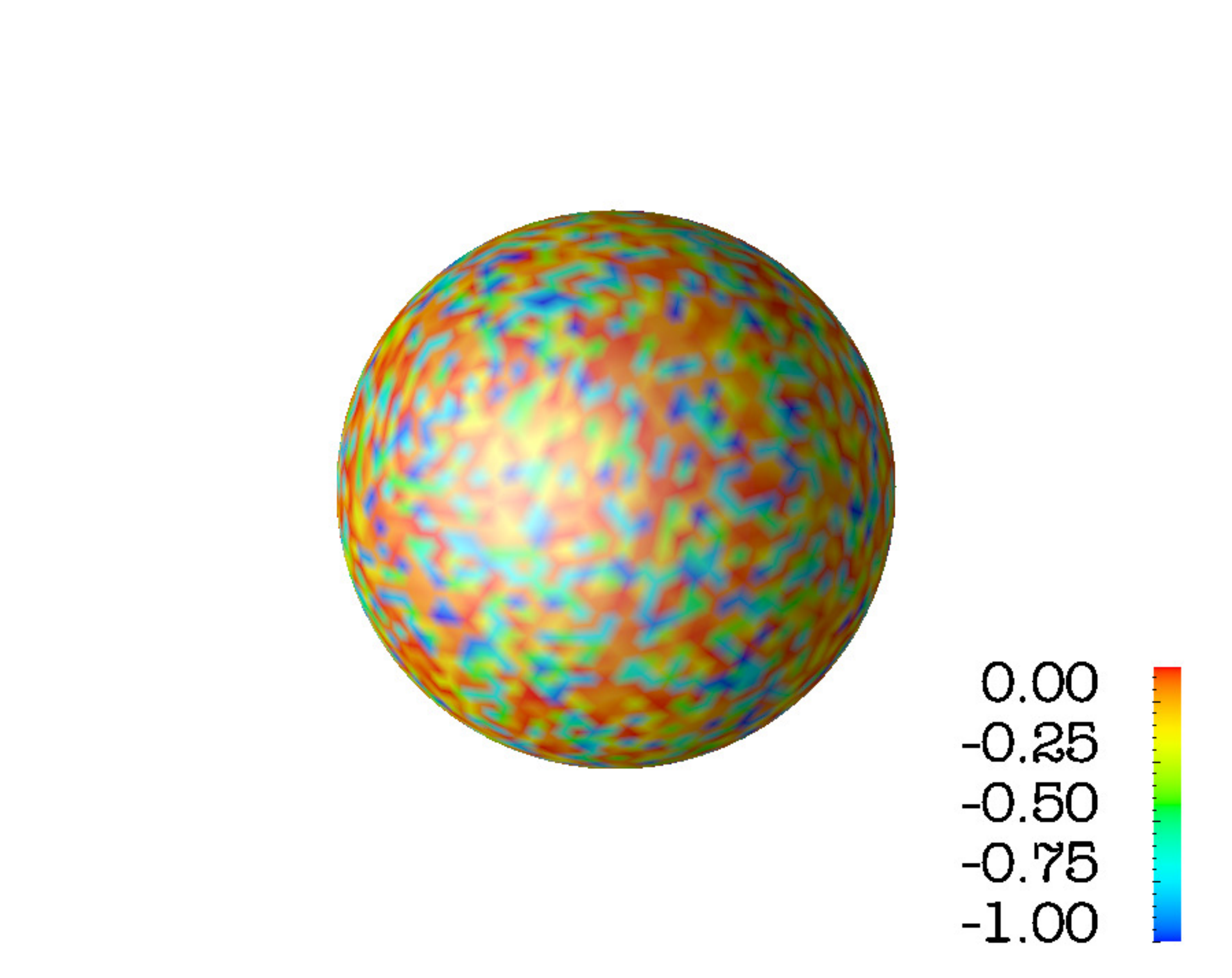} 
\includegraphics[angle=270,width=3.8cm]{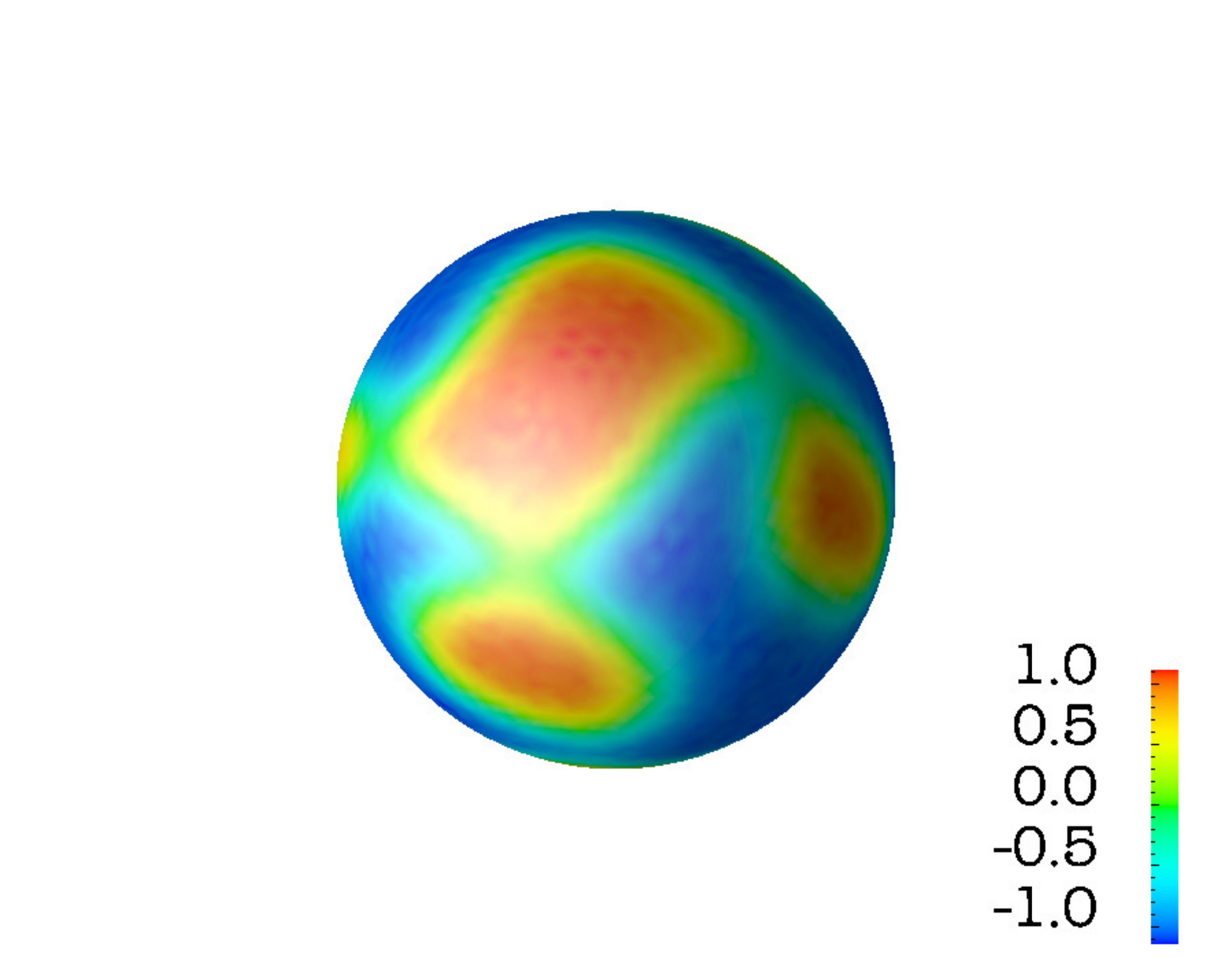}
\includegraphics[angle=270,width=3.8cm]{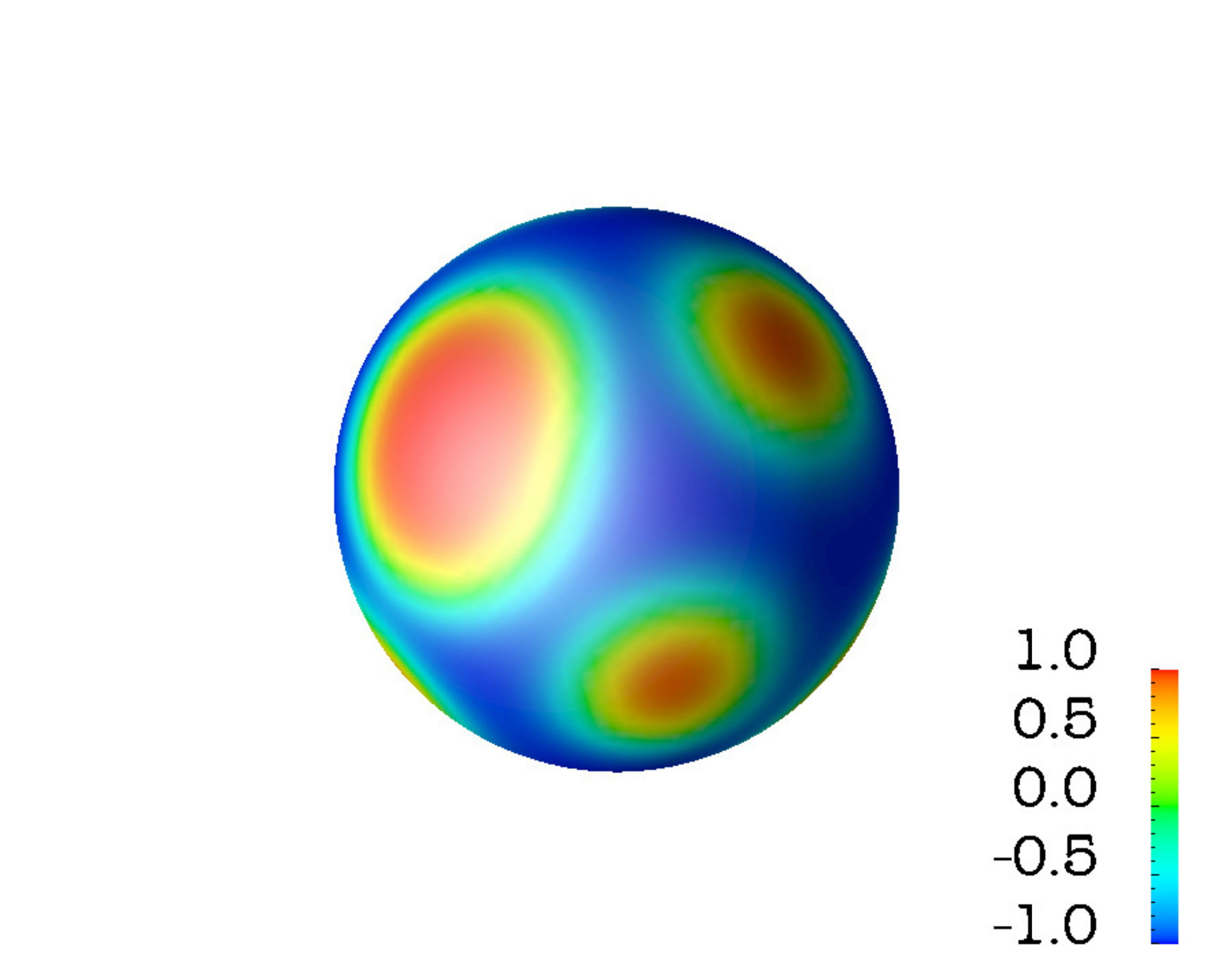} \\
\includegraphics[angle=270,width=3.8cm]{Sim1CurvIni.pdf} 
\includegraphics[angle=270,width=3.8cm]{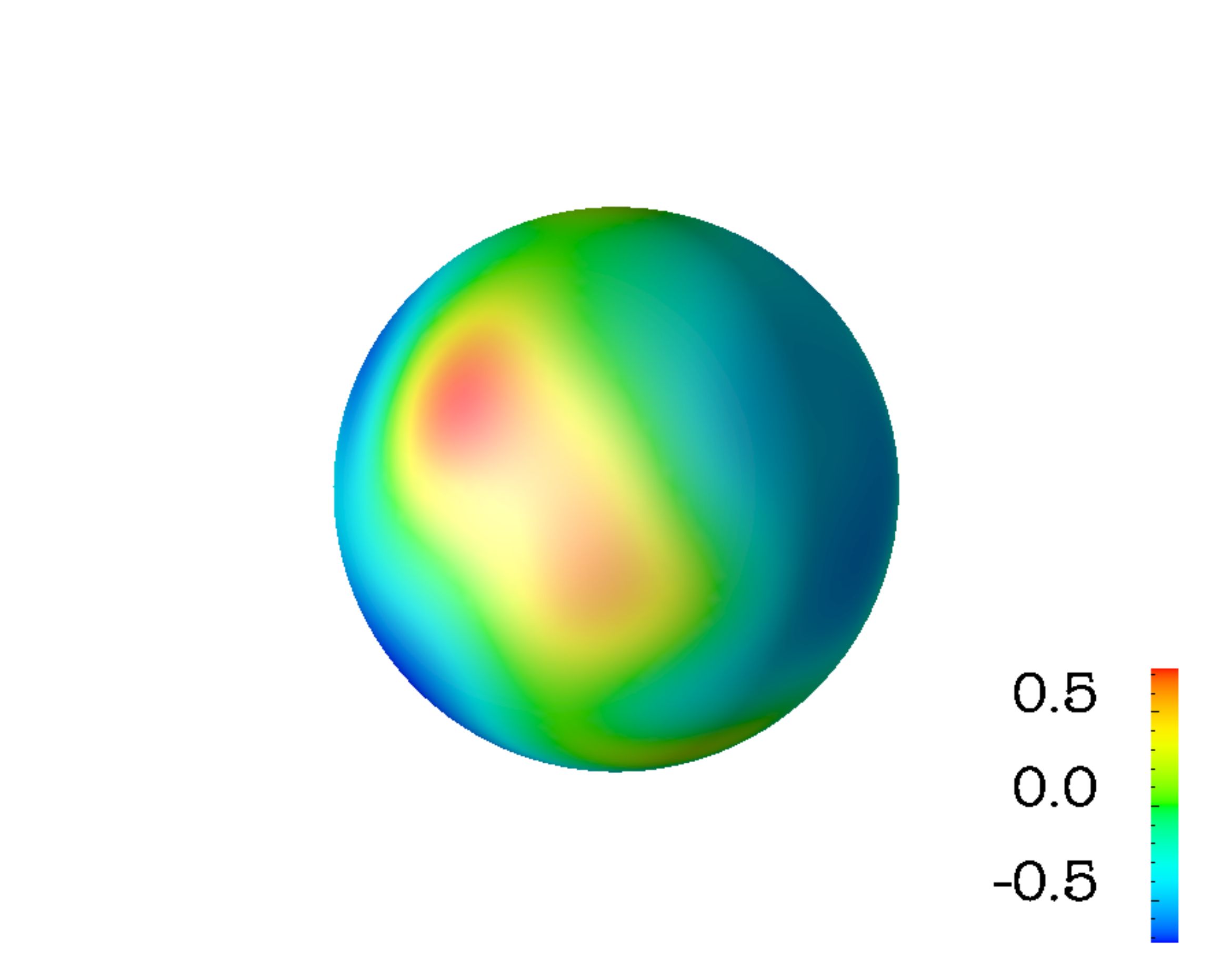} 
\includegraphics[angle=270,width=3.8cm]{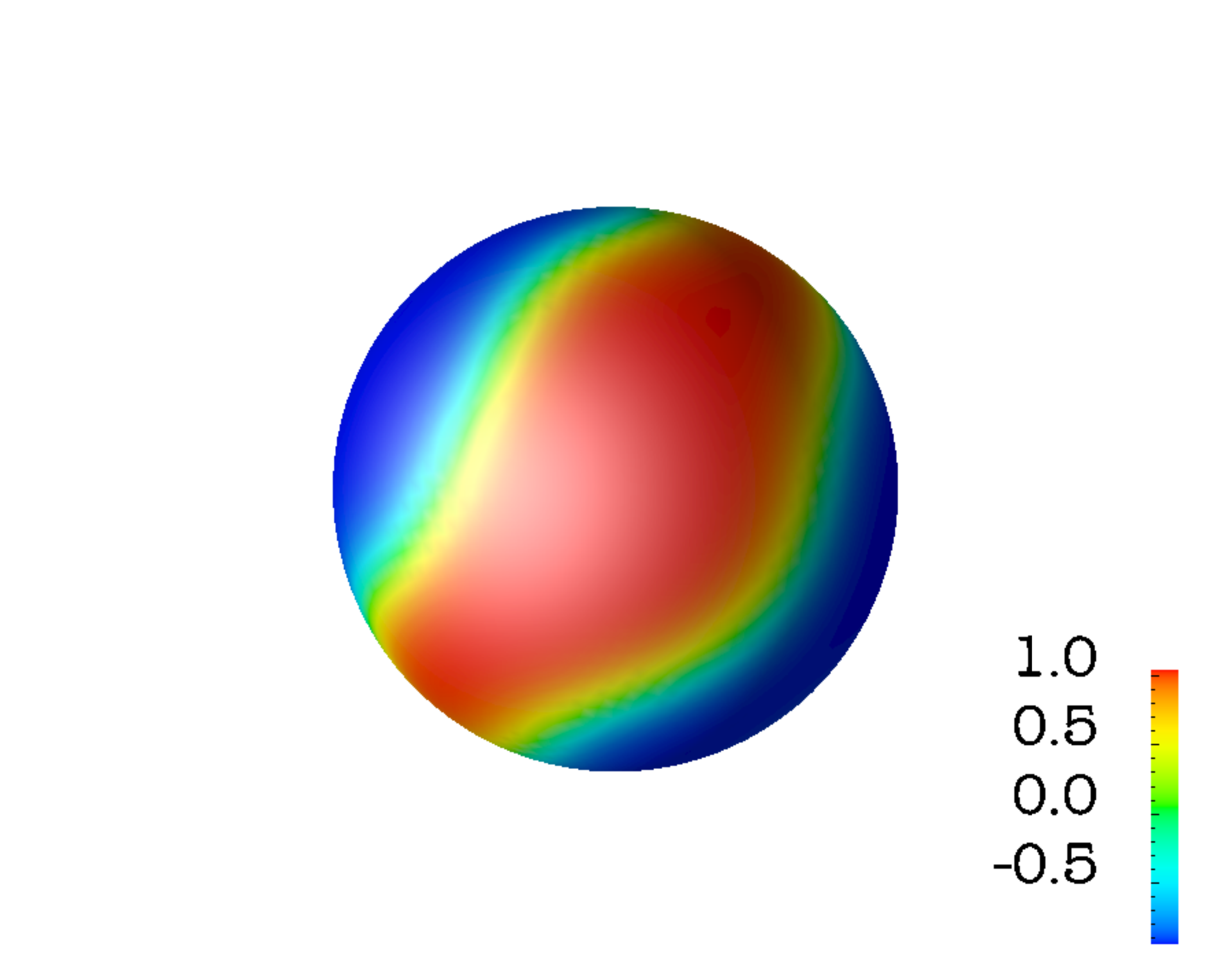}  
\end{center}
\captionsetup{width=0.85\textwidth}
\caption{Simulation \#1. Formation of local microdomains simulated by the geodesic curvature energy (top row) and domain separation simulated by the classical Ginzburg-Landau energy (bottom row) from the same initial random field (left column) on the unite sphere with $3963$ nodes. Sampling time  from left to right is: $t=0, 3,$ and 7.}
\label{fig:Sim1}
\end{figure}

In the first simulation (\#1) on unit sphere with $3963$ approximately uniformly distributed nodes,  we choose $\varepsilon = 0.1, H_c = \frac{1}{0.3}, k = 0.01$  and $\Delta t = 0.001$. A random field is initialized on the surface such that $\int_{S} \phi ds = 0$. The results are compared side by side with those of the classical Allen-Cahn equation in Fig.~\ref{fig:Sim1}. Using a K-means clustering method we are able to identify a number of
microdomains whose radii are then calculated. The radius associated with each microdomain is approximately $0.23$. This means the curvature is  approximately $\frac{1}{0.23}$, close to the specified spontaneous geodesic curvature. 

\begin{figure}[!htb]
\begin{center}
\includegraphics[width=5.5cm]{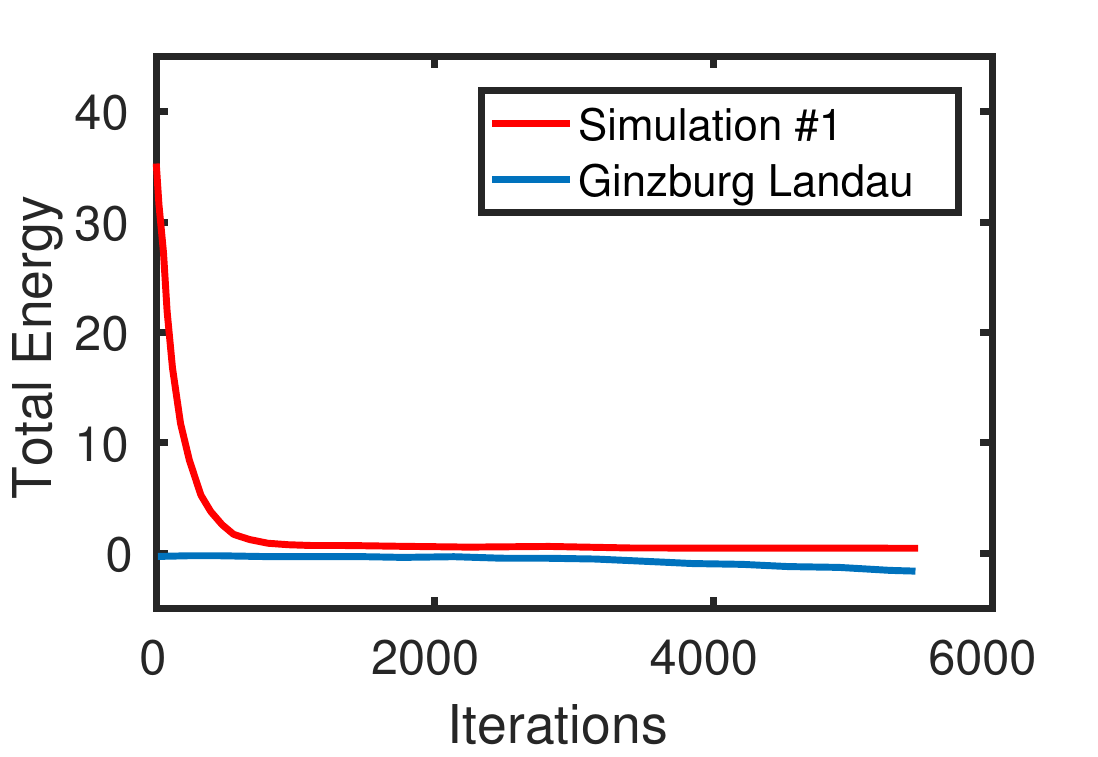}
\includegraphics[width=5.5cm]{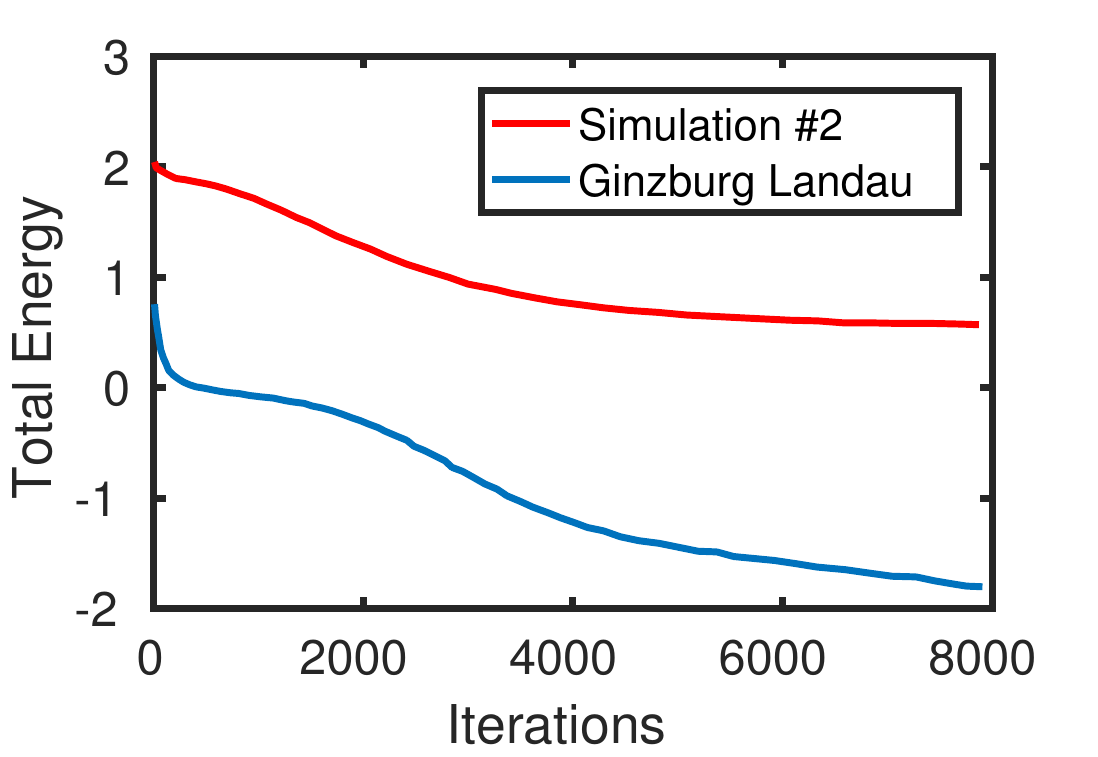}
\end{center}
\captionsetup{width=0.85\textwidth}
\caption{Minimization of the geodesic curvature total energy and the Ginzburg-Landau Energy. Left: Simulation \#1 on unit sphere with $3963$ nodes and $H_c = \frac{1}{0.3}$.
Right: Simulation \#2 on unite sphere with $984$ nodes and $H_c=1/0.4$.} 
\label{fig:energy_sim_12}
\end{figure}

The total energies for the geodesic curvature model and the classical Allen-Cahn model are plotted in Fig.~\ref{fig:energy_sim_12}. Both converge as time evolves. The number of iterations is large because of the small $\Delta t$, which is constrained by the stability of our numerical method  for the fourth-order nonlinear partial differential equation. 

\begin{figure}[!htb]
\begin{center}
\includegraphics[angle=270,width=3.8cm]{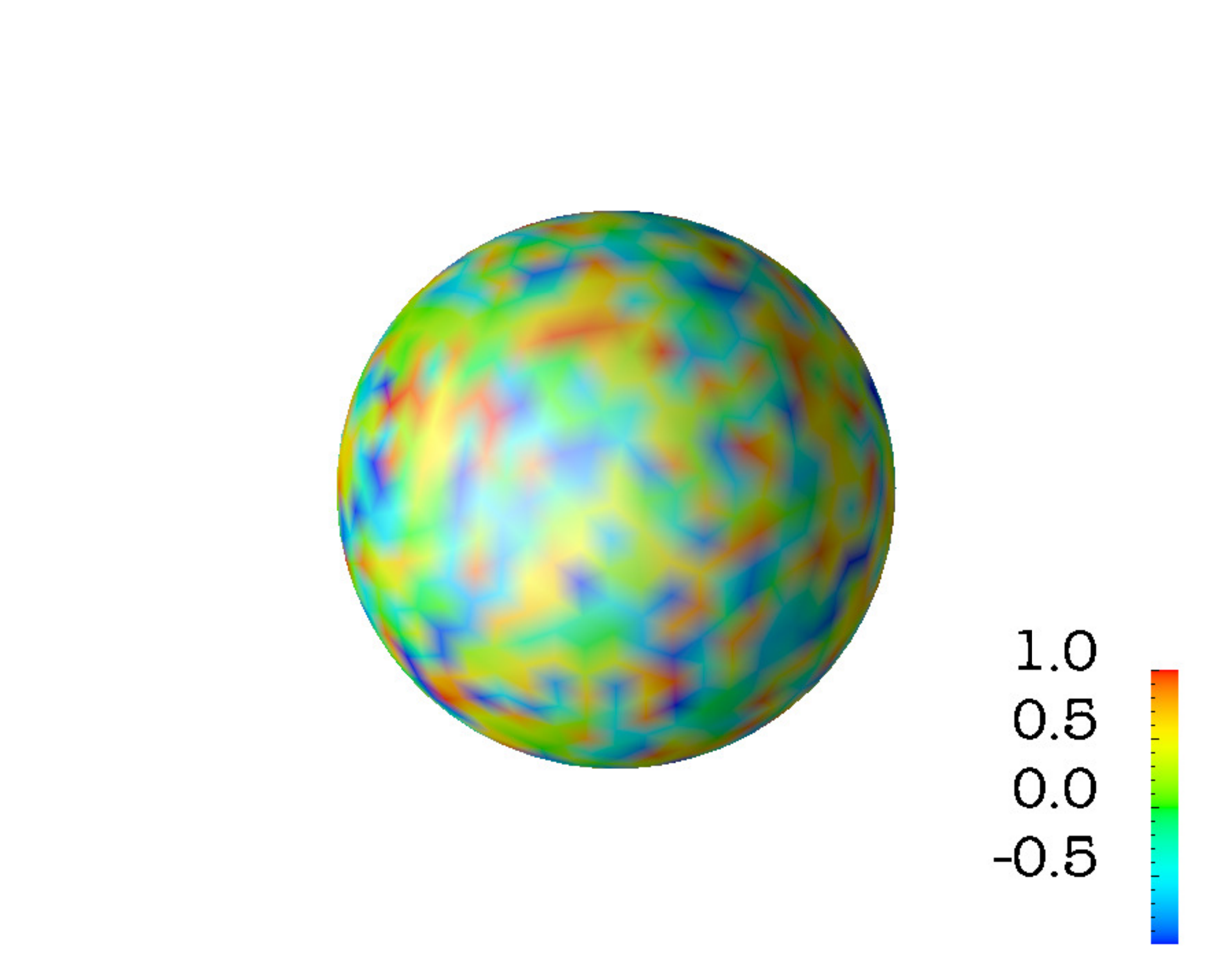} 
\includegraphics[angle=270,width=3.8cm]{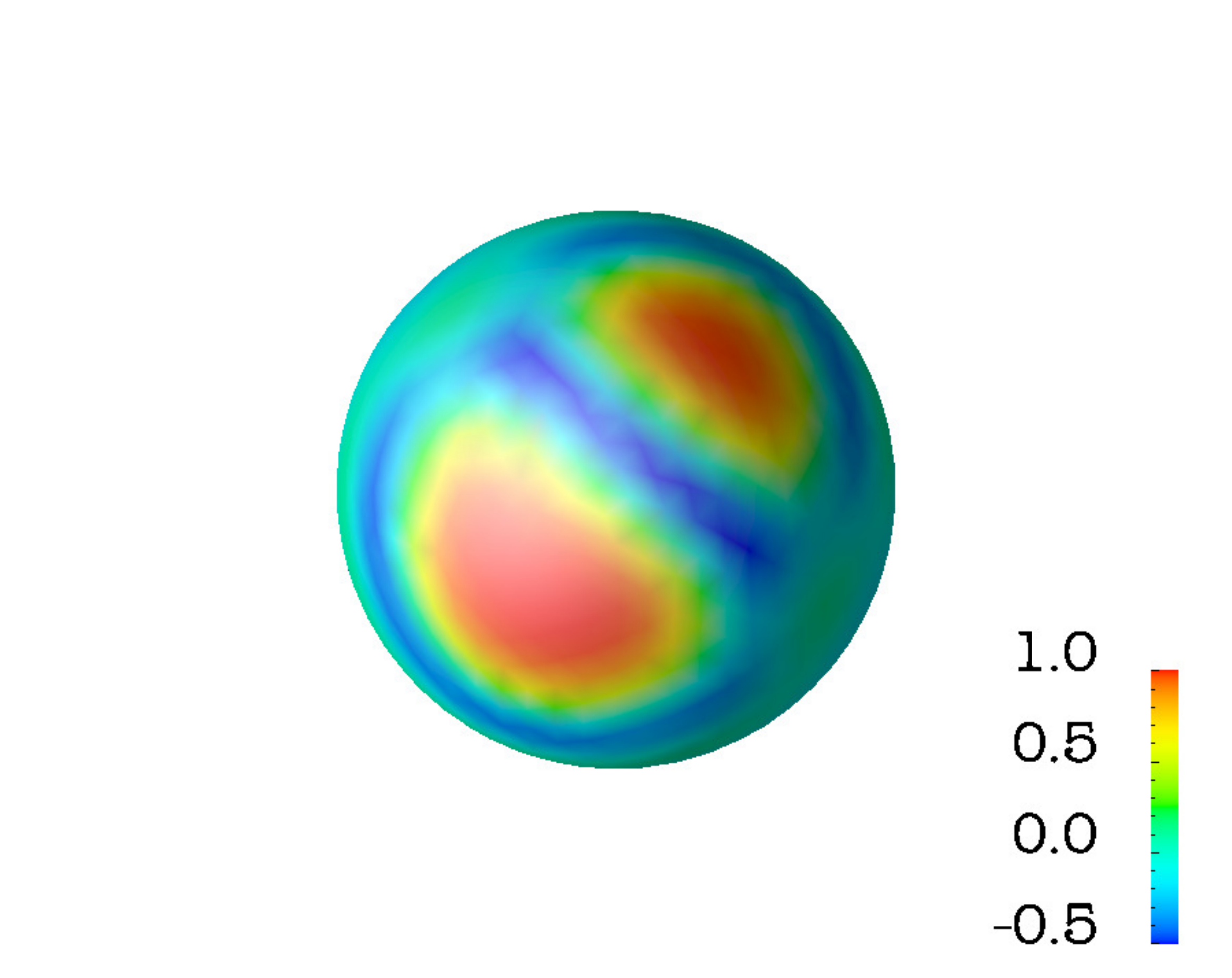} 
\includegraphics[angle=270,width=3.8cm]{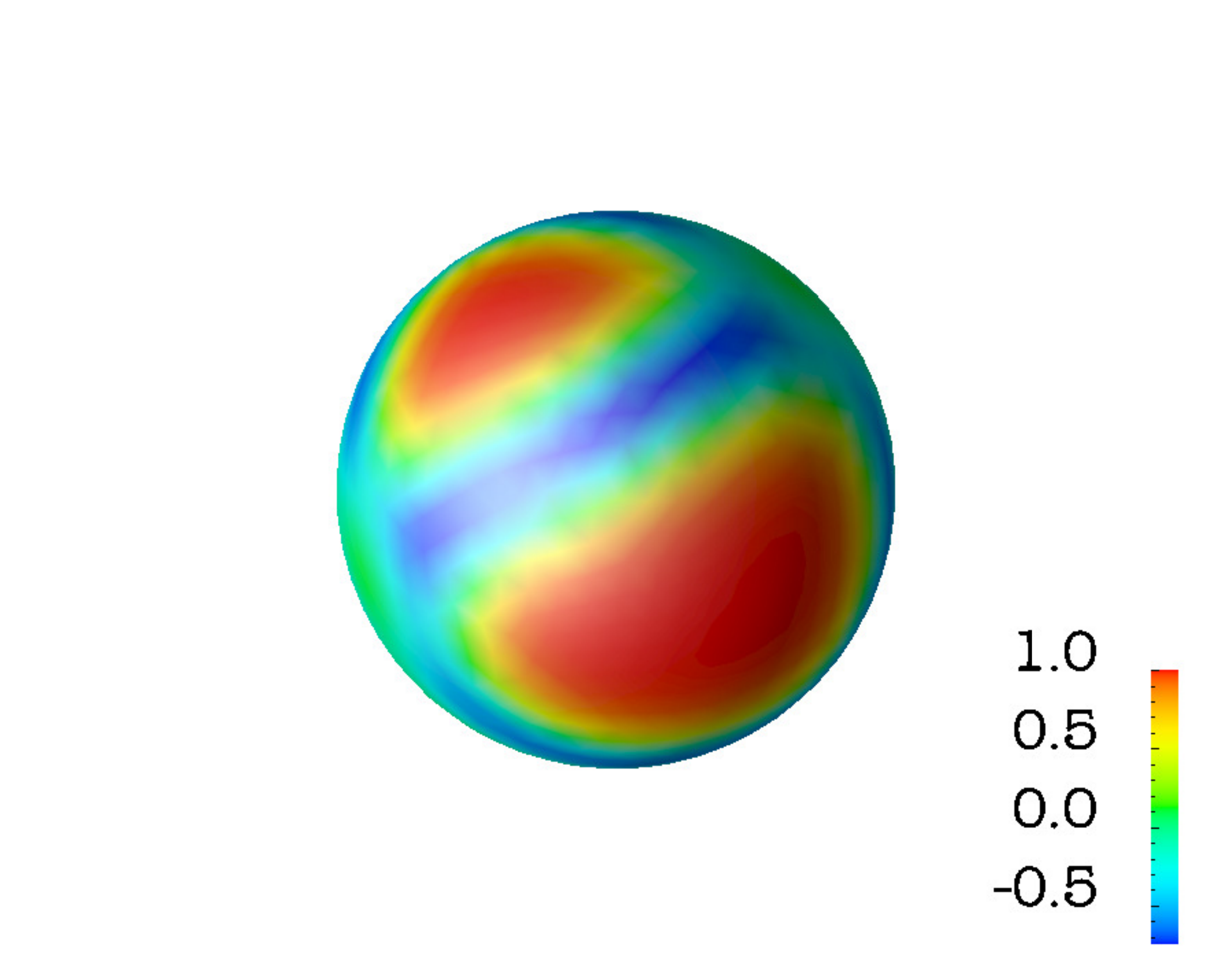} \\
\includegraphics[angle=270,width=3.8cm]{Sim3CurvIni.pdf} 
\includegraphics[angle=270,width=3.8cm]{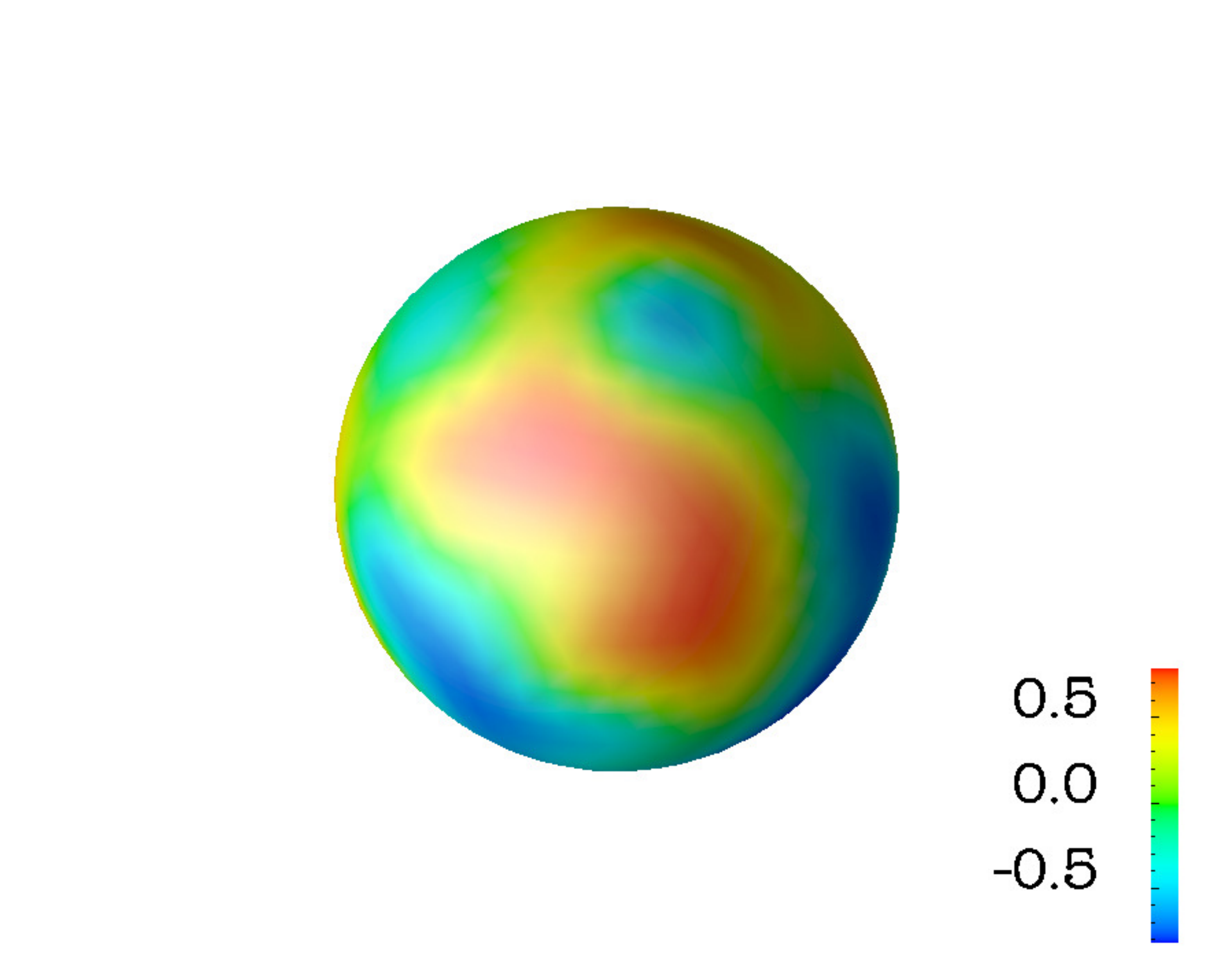} 
\includegraphics[angle=270,width=3.8cm]{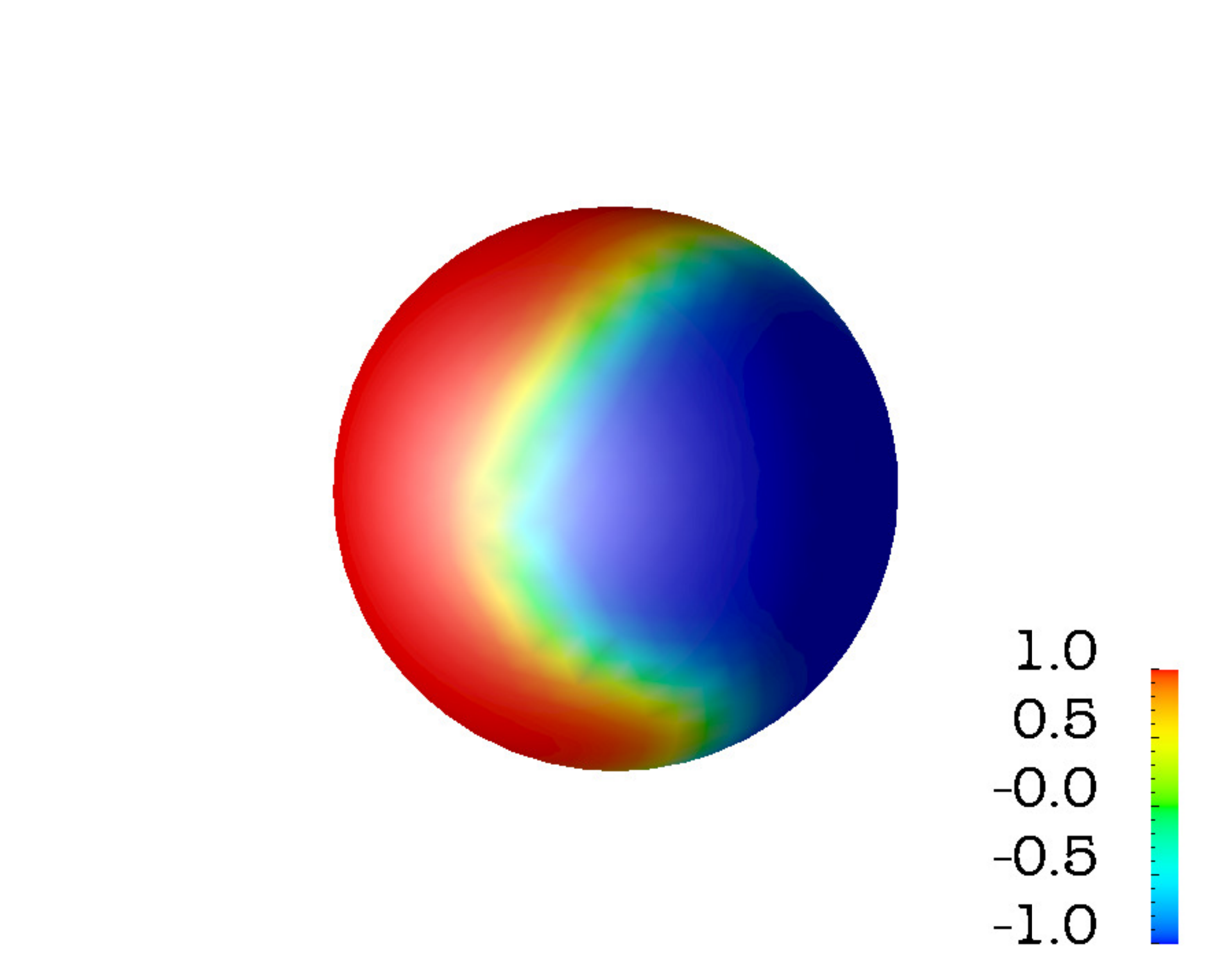}
\end{center}
\captionsetup{width=0.85\textwidth}
\caption{Simulation \#2. Formation of local microdomains simulated by the geodesic curvature energy (top row) and domain separation simulated by the classical Ginzburg-Landau energy (bottom row) from the same initial random field (left column) on unit sphere with $984$ nodes. 
Sampling time from left to right is: $t=0, 3,$ and 7.}
\label{fig:Sim2}
\end{figure}

In the second simulation (\#2) on the unit sphere as shown in  Fig. \ref{fig:Sim2}, we choose $\varepsilon = 0.1, H_c = \frac{1}{0.40}, k = 0.01$ and $\Delta t = 0.002$. This spontaneous curvature matches the reported spontaneous curvature for DOPE/DOPS mixture \cite{FullerN2003a}. A coarser while quasi-uniform  mesh with $984$ nodes is deployed on the unit sphere. The radius associated with the each microdomain is approximately $0.37$, indicating  a curvature approximately $\frac{1}{0.37}$. The convergence of the energies of the geodesic curvature model and the classical Allen-Cahn mode are plotted in Fig.~\ref{fig:energy_sim_12} as well. The lower resolution resulting from the coarser mesh in the second simulation can be seen in the larger spots  in the initial field and the wider transitional layers between different domains.

\begin{figure}
\begin{center}
\includegraphics[angle=270,width=3.8cm]{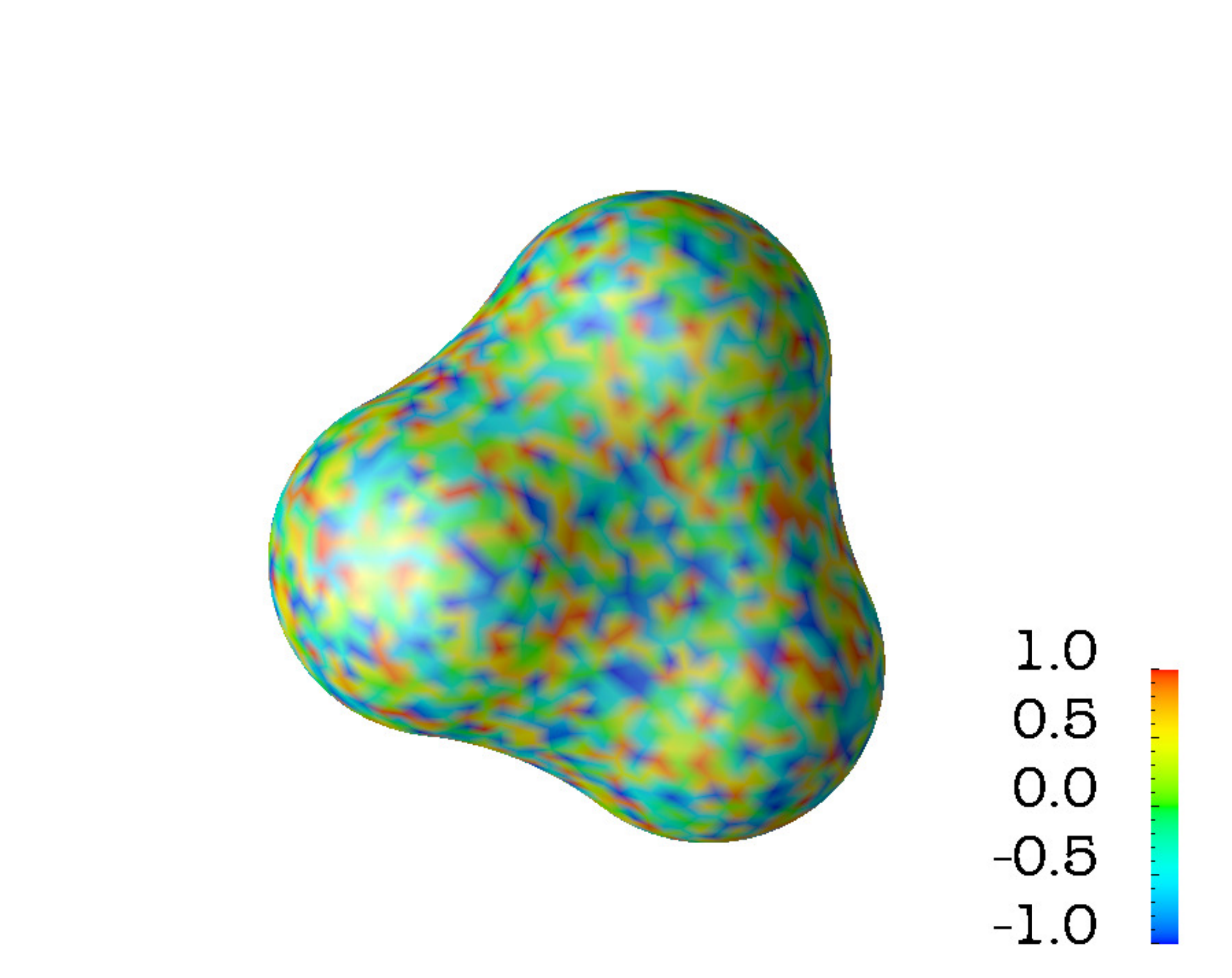}
\includegraphics[angle=270,width=3.8cm]{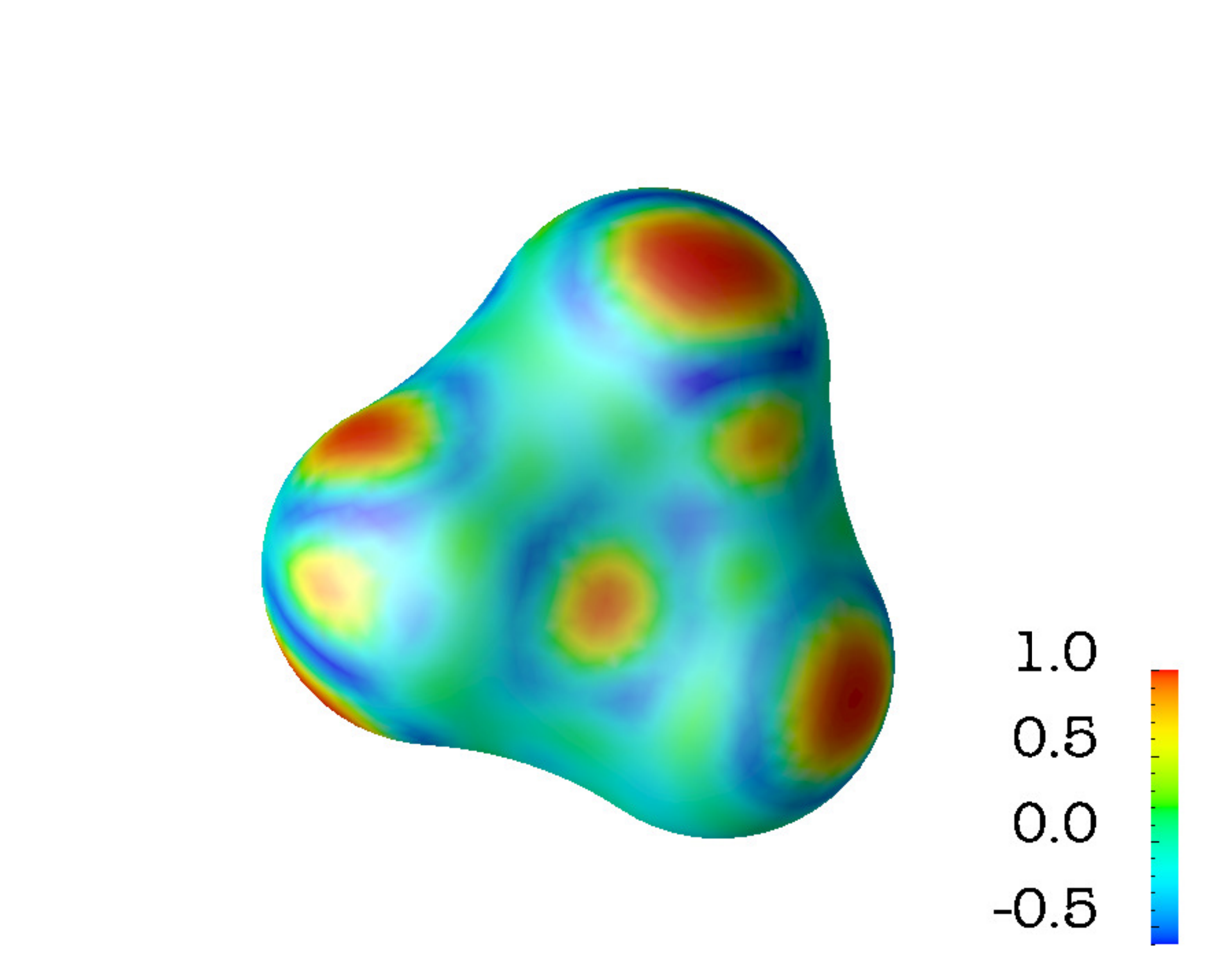} 
\includegraphics[angle=270,width=3.8cm]{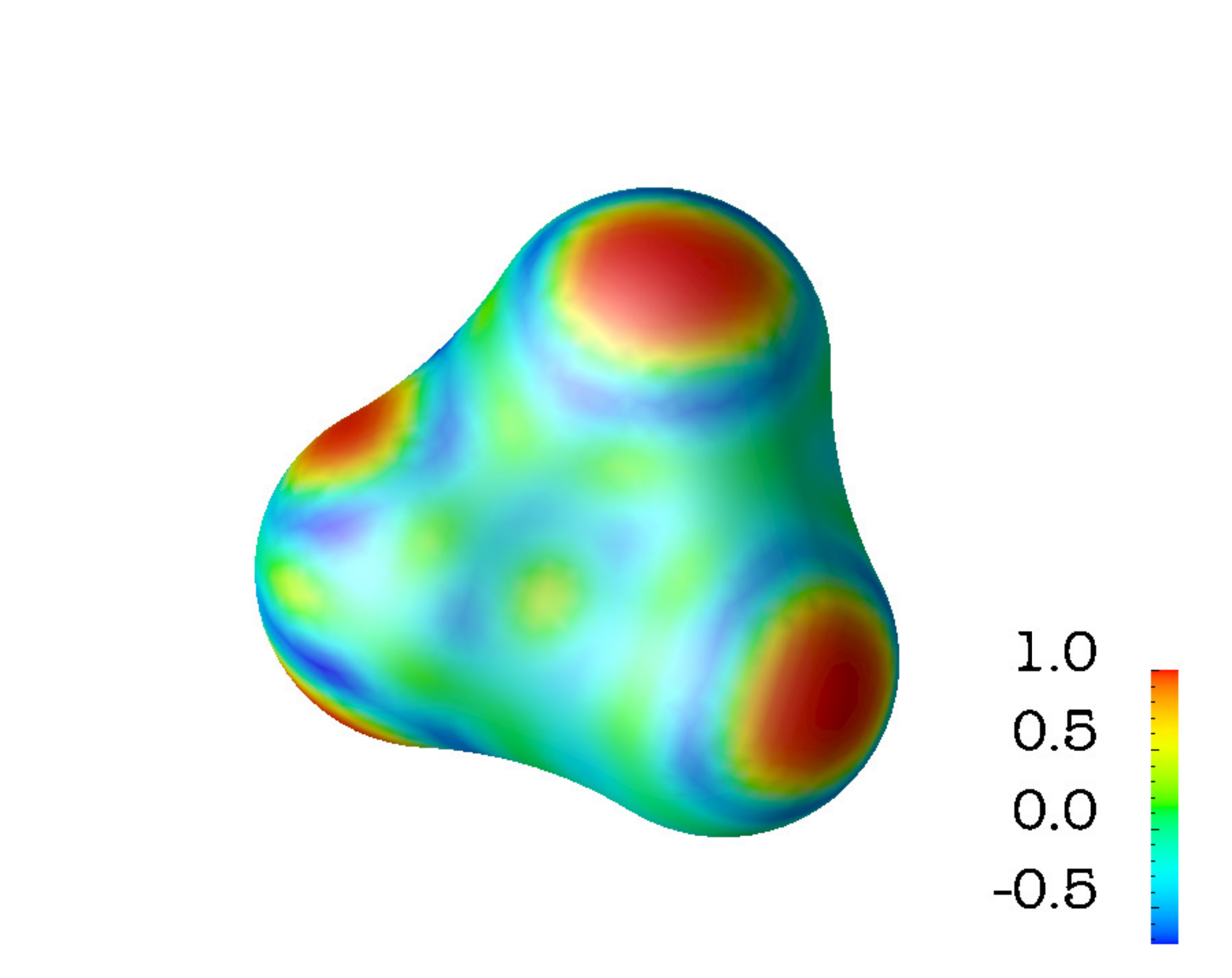}  \\
\includegraphics[angle=270,width=3.8cm]{Sim4CurvIni.pdf}
\includegraphics[angle=270,width=3.8cm]{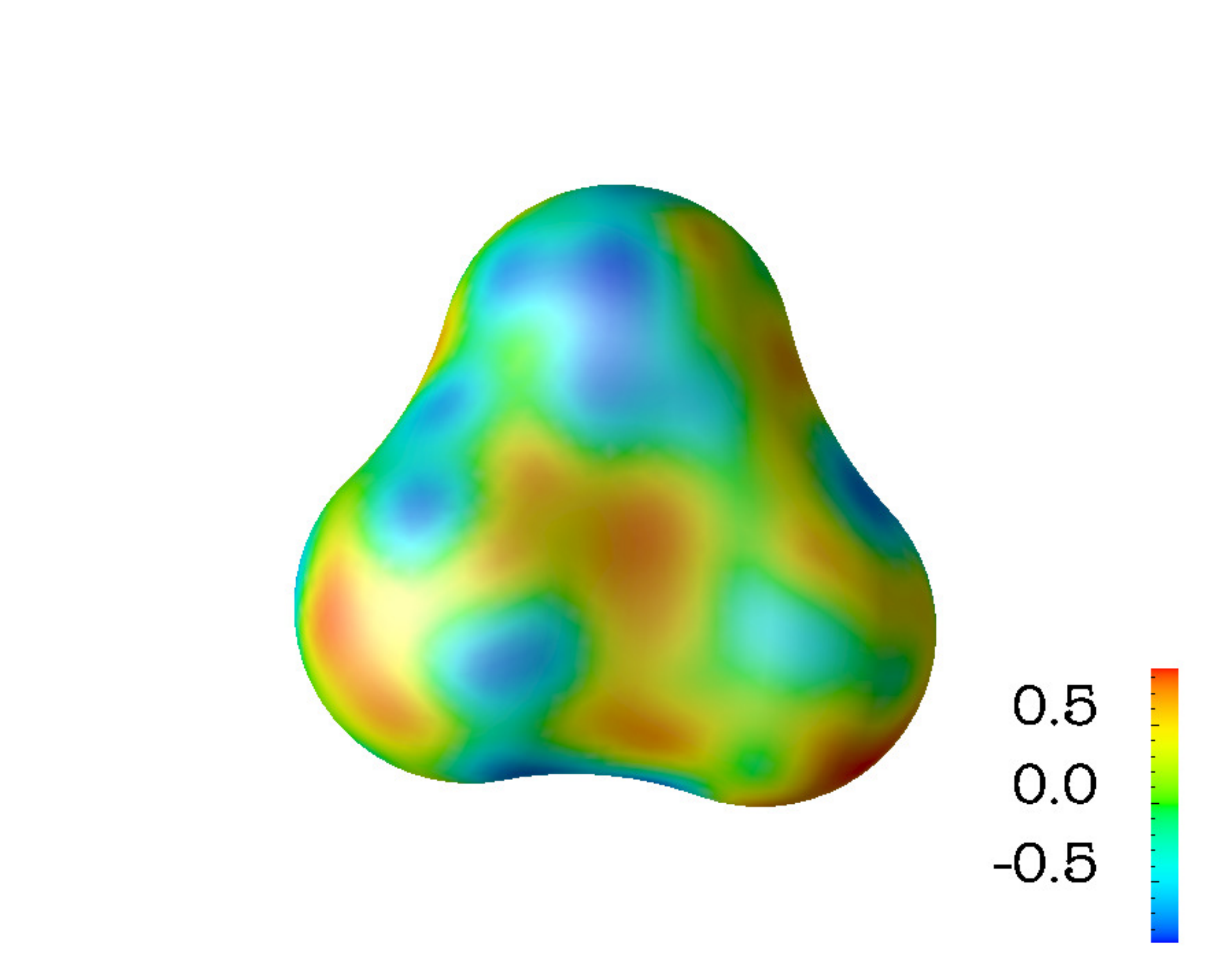}
\includegraphics[angle=270,width=3.8cm]{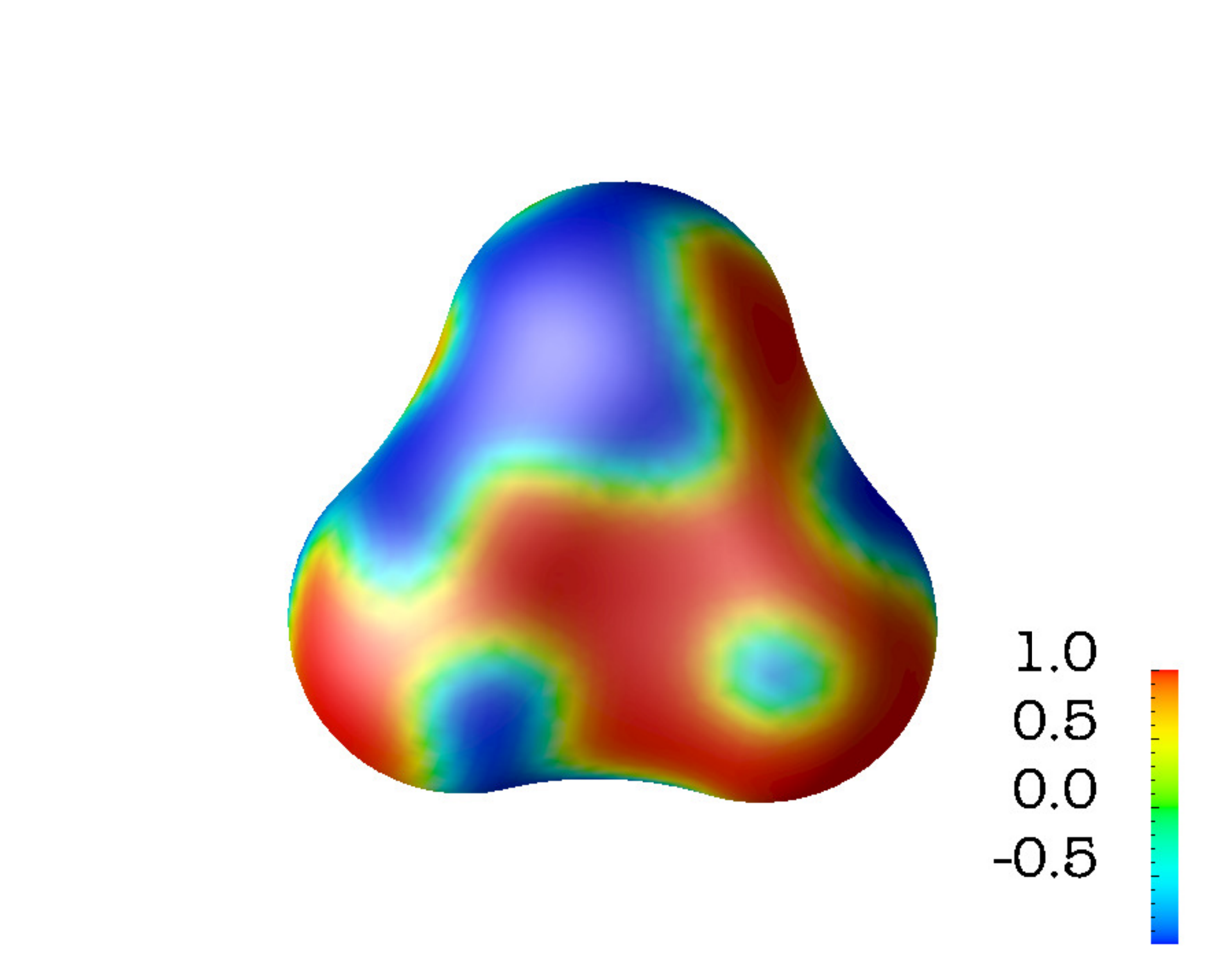} 
\end{center}
\captionsetup{width=0.85\textwidth}
\caption{Simulation \#3. Formation of local microdomains simulated by the geodesic curvature energy (top row) and domain separation simulated by the classical Ginzburg-Landau energy (bottom row) from the same initial random field (left column) on the molecular surface of three-atom with $2974$ nodes. Sampling time from left to right is: $t=0, 3, $ and 7.}
\label{fig:Sim3}
\end{figure}

\begin{figure}[!htb]
\begin{center}
\includegraphics[width=9cm]{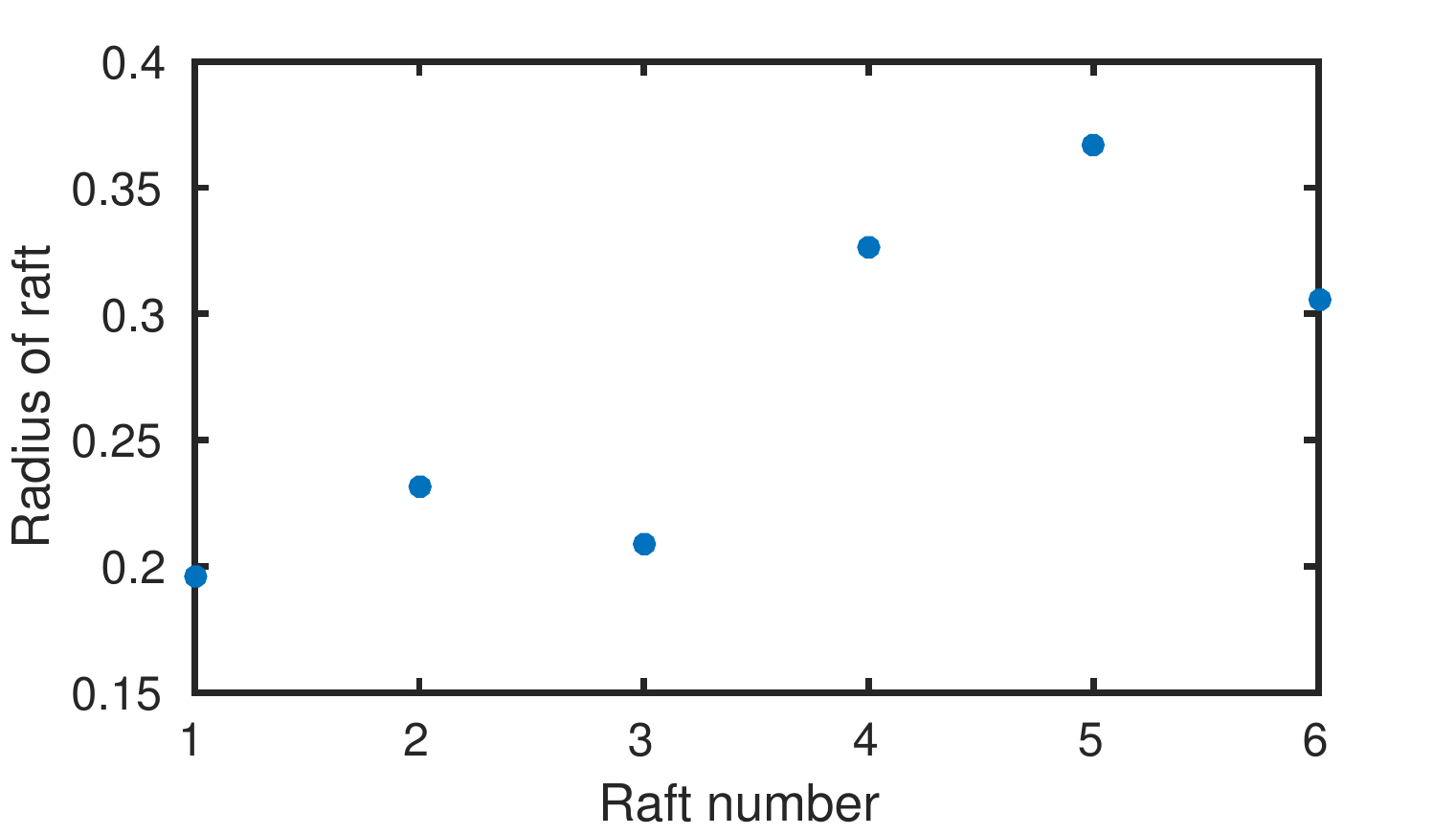}
\end{center}
\captionsetup{width=0.85\textwidth}
\caption{The radii of the prominent 6 microdomains produced in Simulation \#3}
\label{Sim4Rad}
\end{figure}

The third simulation (\#3) is conducted on a more complicated surface as shwon in Fig. \ref{fig:Sim3}. We choose the molecular surface of three particles of unit radius respectively centered at $(0, 1, 0), (-0.864, -0.5, 0)$ and $(0.864, -0.5, 0)$. The surface is quasi-uniformly meshed with $2974$ nodes and  we set $\varepsilon = 0.1, H_c = \frac{1}{0.4}, k = 0.01$ and $\Delta t = 0.001$. Starting with a random initial field we finally identified six  microdomains using the K-means clustering method at the equilibrium state, whose radii are estimated. As seen in Fig.~\ref{Sim4Rad}, the radii of the microdomains approximate 
the given spontaneous geodesic curvatures. 

\begin{figure}[!htb]
\begin{center}
\includegraphics[width=9cm]{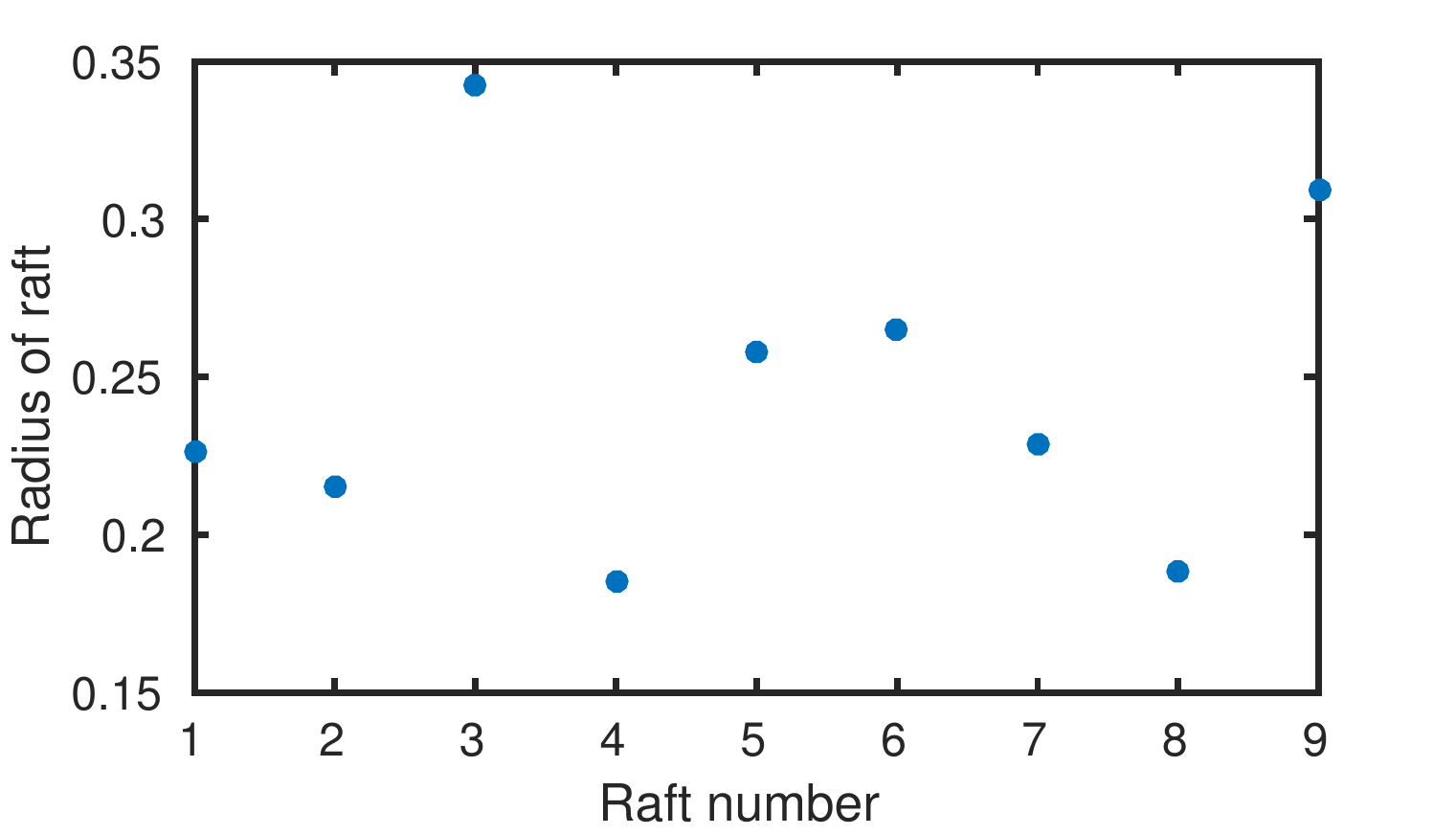}
\end{center}
\captionsetup{width=0.85\textwidth}
\caption{The radii of the prominent 9 rafts produced by Simulation \#4}
\label{Sim5Rad}
\end{figure}

In the last simulation (\#4) we choose the molecular surface of six particles of unit radius respectively centered at $(1,0,0),(-1,0,0),(0,1,0),(0,-1,0),(0,0,1)$ and $(0,0,-1)$. The quai-uniform surface mesh has $3903$ nodes and we  set $\varepsilon = 0.1, H_c = \frac{1}{0.4}, k = 0.01$ and $\Delta t = 0.001$ for the simulation. One can see from Fig. \ref{Sim5Rad} that the largest raft radius obtained by the simulation is about $0.35$ which means the curvature of that raft is about $\frac{1}{0.35}$, a value close to given spontaneous geodesic curvature. 

\begin{figure}[!htb]
\begin{center}
\includegraphics[angle=270,width=3.8cm]{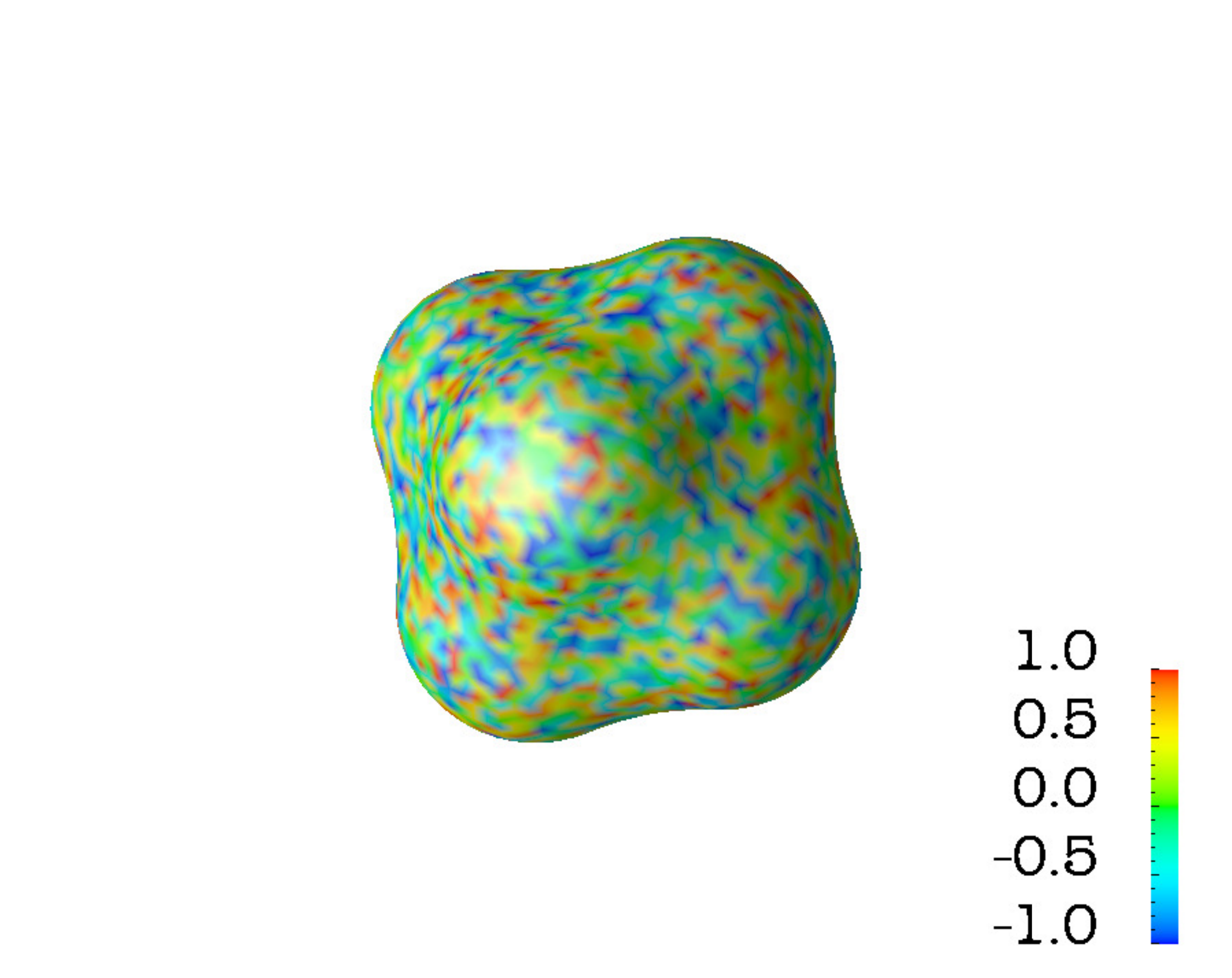} 
\includegraphics[angle=270,width=3.8cm]{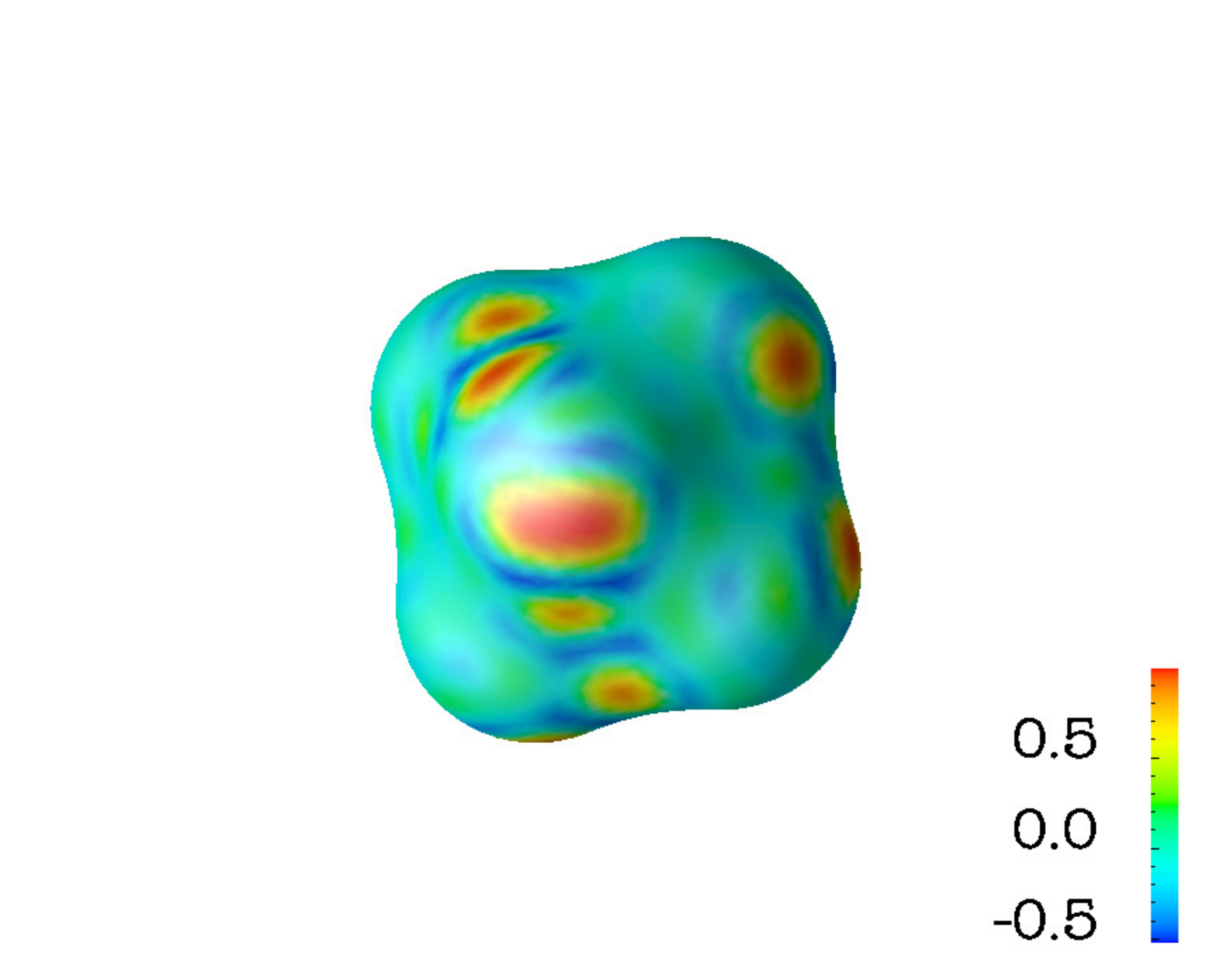} 
\includegraphics[angle=270,width=3.8cm]{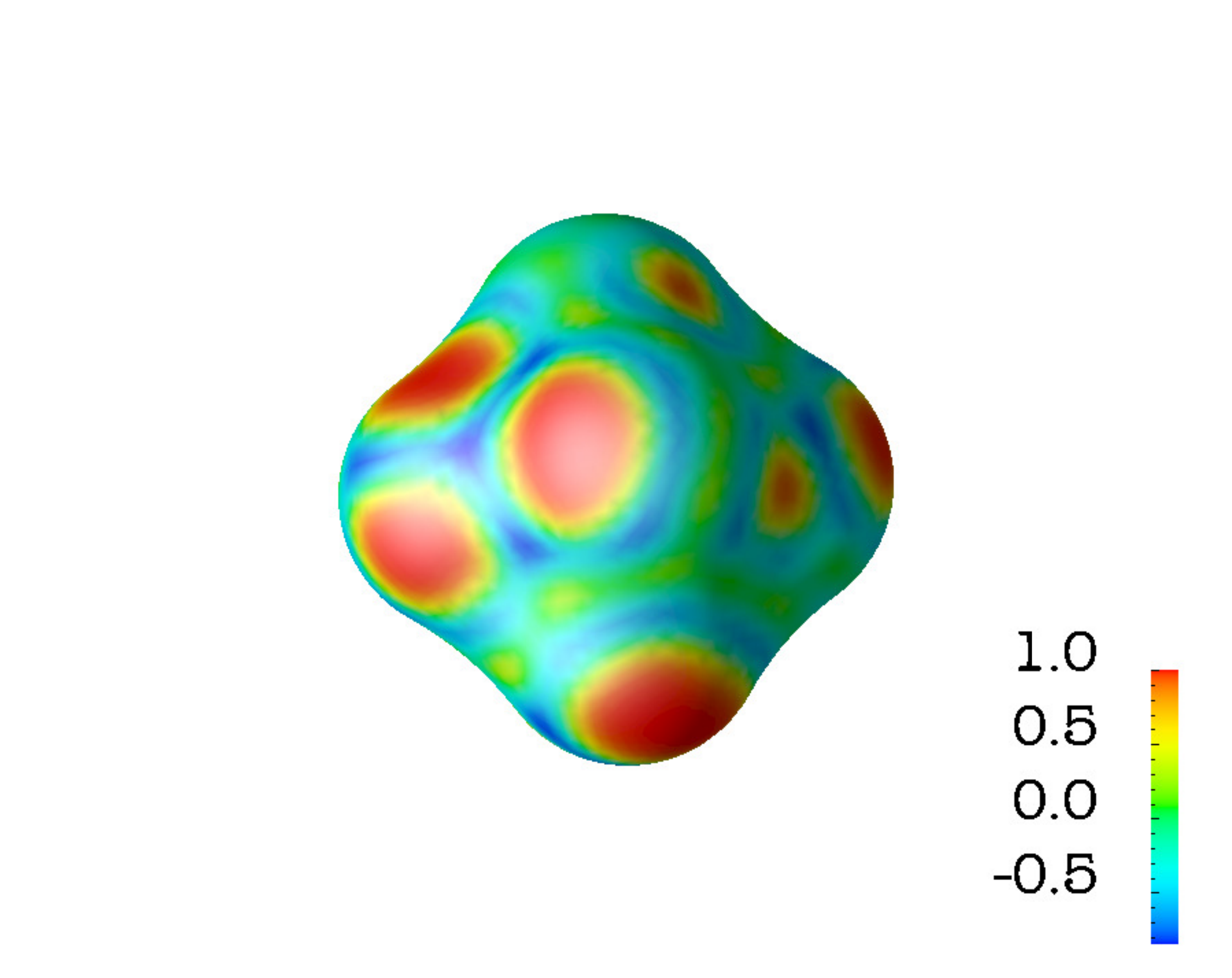} \\  
\includegraphics[angle=270,width=3.8cm]{Sim5CurvIni.pdf} 
\includegraphics[angle=270,width=3.8cm]{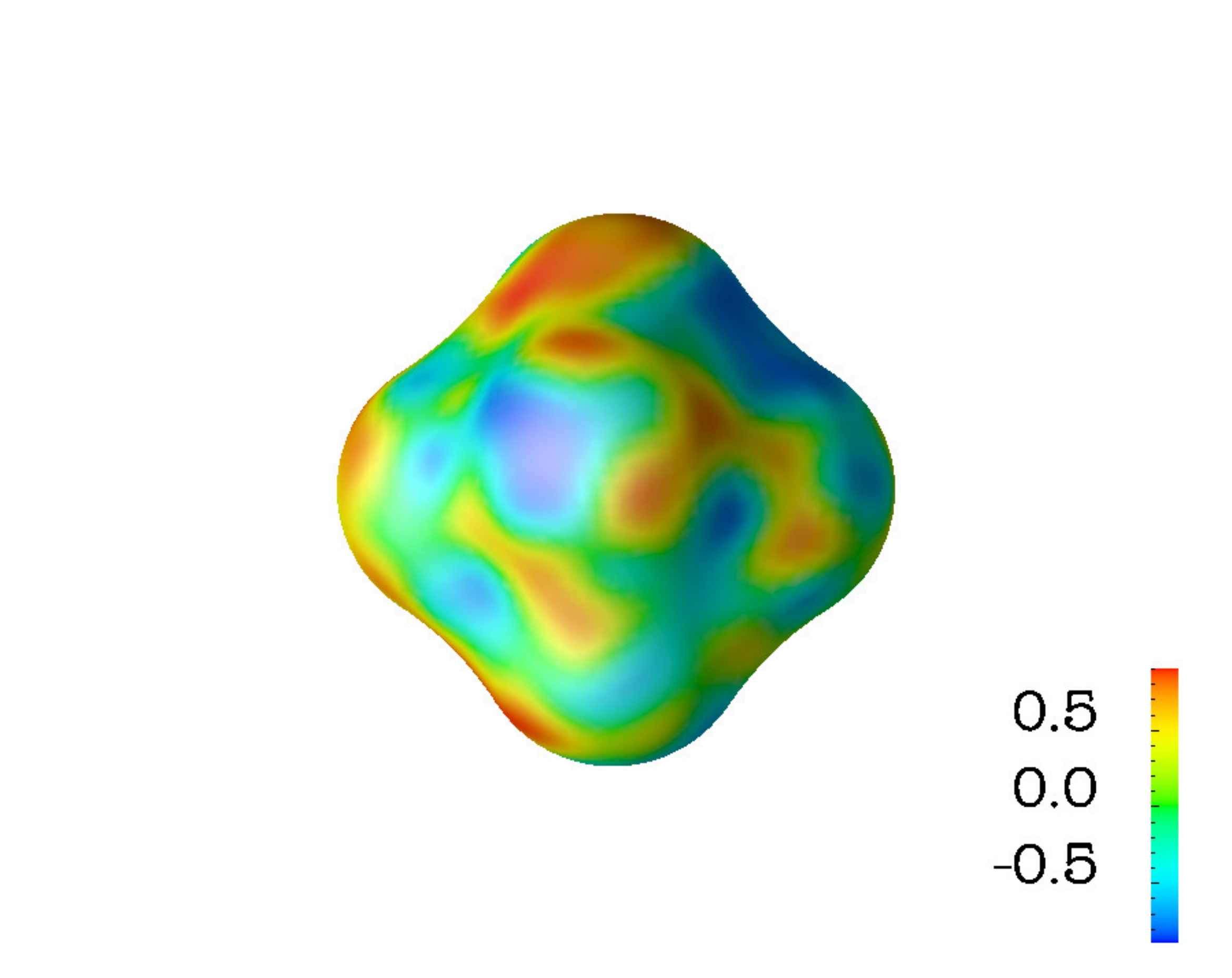} 
\includegraphics[angle=270,width=3.8cm]{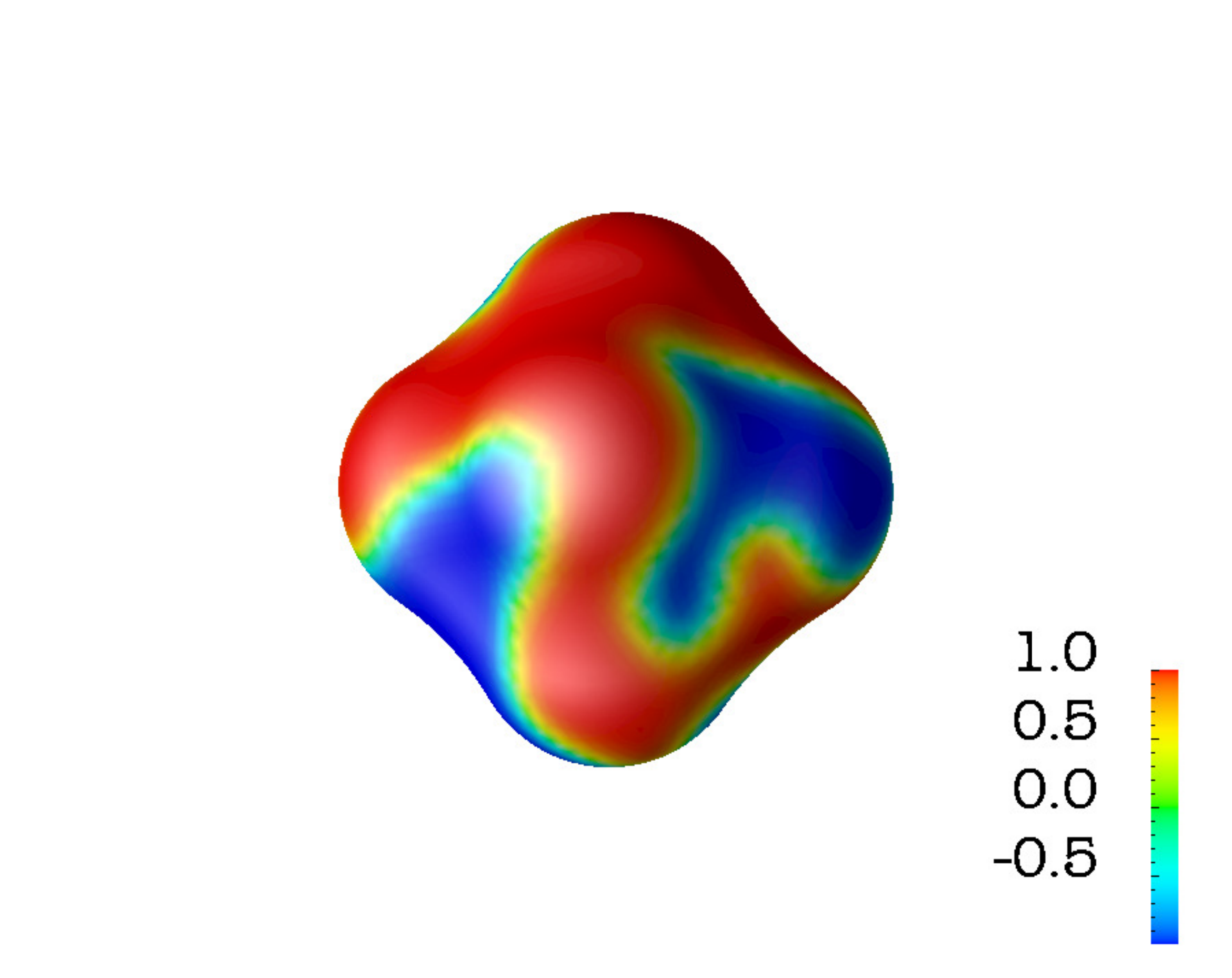} 
\end{center}
\captionsetup{width=0.85\textwidth}
\caption{Simulation \#4. Formation of local microdomains simulated by the geodesic curvature energy (top row) and domain separation simulated by the classical Ginzburg-Landau energy (bottom row) from the same initial random field (left column) on the molecular surface of six-atom with $3903$ nodes. Sampling time from left to right is: $t=0, 3, $ and 7.}
\label{fig:Sim4}
\end{figure}

The radii of the microdomains generated in our simulations are not exactly the given spontaneous geodesic curvature. Rather they are distributed around the given curvature. Apart from the numerical error in simulation and in K-mean clustering and radii estimate, this non-uniform distribution of domain radii is mostly related to the total quantity of the lipid phases in the initial random field. The initial quantity may not exactly cover an integer number of microdomains with the given radius. However, the overall distribution of radii around the given radius of curvature demonstrated that our geodesic curvature model is capable of predicting the formation of microdomains that are caused by the geometrical and molecular mechanical mismatch of lipid mixtures. The predicted microdomains can be compared to the observed lipid rafts, and the boundaries of these microdomains can be identified to provide locations where specific proteins can aggregate. Coupling of our model of geodesic curvature driven microdomains formation to the localization of proteins will provide a very useful quantitative technique for studying the crucial roles of these proteins in high-fidelity signal transduction in cells \cite{LiH2013a,HockerH2013a}.

\section{Variational Methods for Curvature Induced Protein Localization in Bilayer Membranes} \label{sec:protein}

Rather than forming distinct domains in a way similar to lipids as modeled in Sect. \ref{sec:surfpattern}, many membrane proteins do not form distinct domains in membranes \footnote{A protein unit consisting of several segments such as most ion channel proteins or  G-protein-coupled receptors (GPCRs) 
 is not taken as a distinct domain in this study. The whole unit is considered as a single protein instead.}. Given the fact that their distribution on bilayer membranes is not uniform,  molecular mechanisms need to be identified to quantitatively investigate this distribution and its biological consequences. On the one hand, approximately 30-90\% of all membrane proteins can freely diffuse along the membrane \cite{FaraidoJ2002a,KimK1998a,MartinireA2012a,RamadS2009a}. On the other hand, insertion or tethering of the membrane proteins to bilayer membrane will cause membrane curvature \cite{ZimmerbergJ2006a,ReynwarB2006a,GumbartJ2012a}. For instances, the rigid proteins such as those in the BAR (Bin/Amphiphysin/Rvs) domain family can act as a scaffold to the membrane. These proteins have an intrinsic curvature and, upon attaching, the membrane bends to match the protein curvature \cite{MimC2012b}. In a similar fashion, several proteins can oligomerize to create a rigid shape and bend the membrane. Protein coats such as clathrin, COPI (COat Protein I) and COPII (COat Protein II) are examples of this type \cite{FarsadK2003a,KirchhausenT2012a}. Other proteins may insert themselves into the membrane. Membrane curvature is also induced when there is a difference between the length of the hydrophobic region of a membrane protein and the thickness of the hydrophobic core of the lipid bilayer in which it is embedded \cite{PartonD2011a}. Epsin proteins do this by forming an alpha-helix known as H0 upon binding to the membrane, and this helix inserts itself into the membrane \cite{BaumgartT2011a}. Moreover, local crowding of peripheral proteins can cause membrane bending by creating an asymmetry of the monolayer areas and thereby curling the membrane away from the side on which the crowding occurred. This effect is experimentally demonstrated in \cite{StachowiakJ2010a}. Further illustrating the importance of proteins in membranes, Schmidt {\it et. al.} showed that the M2 protein plays an essential role in generating regions of high curvature in the influenza A virus membrane \cite{SchmidtN2013a}. This specific protein accumulates in regions of negative Gaussian curvature and can generate curvature in the membrane itself, allowing the replicated virus to be wrapped and released from the infected cells. While these examples should provide sufficient motivation to include proteins to the model, we note that all endocytosis and exocytosis  processes are promoted in one way or another by proteins. Therefore, any viral replication process requires proteins. Antagonizing the  curvature effects of proteins is a viable antiviral strategy \cite{SchmidtN2013a}. This motivates the necessity for a model coupling membrane curvature and lateral diffusion of proteins. We shall 
observe 
below that the final governing equation for this curvature-driven lateral transportation appears a drift-diffusion equation in its essential form. This mechanism is different from the transporation of surfactants on interfaces moving with the fluid flow as investigated in the literature \cite{TeigenK2009a,TeigenK2011a,XuJ2012a}.

\subsection{Lagrangian formulation}

Modeling generation of membrane curvature using energetic variational principle has been well established in the past few decades \cite{CanhamP1970,HelfrichW1973,EvansE1974,DuQ2006a}. These research have been inspiritional to our work. However, the focus of our discussion in this section is on the curvature driven protein localization.
We sketch the framework of the integration of these two components. The numerical implementation is computationally intensive because of the coupling of dynamical membrane morphology and the varying surface concentration of proteins. Consider a membrane with $(m+1)$ distinct lipid species with concentrations $\rho_l^{\rm lip}, l=0, \cdots, m$ and a single type of diffusive membrane proteins with a concentration $\rho^{\rm pro}$. A closed membrane is modeled as a structureless surface $S$ contained in a 3D domain $\Omega\in {\mathbb R}^3$ and separated $\Omega$ into two subdomains, one inside the membrane and the other outside. The total energy of the system is composed of the membrane curvature energy and the entropic energy from the lipids and proteins 
\begin{equation} \label{eqn:total_eng_curvprotein}
G_{\rm tot} = G_{\rm mem} + G_{\rm ent},
\end{equation} 
where the membrane curvature energy is given in the classical Canham-Helfrich-Evans form 
\begin{equation} \label{eqn:G_mem}
G_{\rm mem} = \int_{C} k (H - H_0(\rho_l^{\rm lip}, \rho^{\rm pro}))^2 ds,
\end{equation}
and the entropic energy for the membrane with membrane protein attachments is 
\begin{equation} \label{eqn:G_ent}
G_{\rm ent} = \frac{1}{\beta} \int_{C} \left( \sum_{l=0}^m \rho_l^{\rm lip} \left[ \ln(\rho_l^{\rm lip} (a_l^{\rm lip})^2) -1 \right ] + 
    \rho^{\rm pro} \left[ \ln(\rho^{\rm pro}(a^{\rm pro})^2 ) - 1 \right ] \right )ds,
\end{equation}
Here $H$ is the membrane mean curvature and $H_0$ is the spontaneous membrane curvature, $k$ is a curvature energy coefficient, and $\beta = 1/(k_BT)$ is the inverse of thermal energy. The effective sizes of lipids and proteins are respectively given by $a_l^{\rm lip}$ and $a^{\rm pro}$. By   modeling lipids and proteins as hard disks, the occupied surface areas in the membranes are taken as $(a_l^{\rm lip})^2$ and $(a^{\rm pro})^2$, respectively. The essential feature of our model is seen in the dependence of the membrane spontaneous curvature on the local lipid composition $\rho_{l}^{\rm lip}$ and the protein concentration $\rho^{\rm pro}$. This dependence is justifiable considering that (i)  each lipid species $l$ has its own spontaneous curvature \cite{MarshD2007} therefore the membrane spontaneous curvature must be a function of the local lipid composition, and (ii) membrane proteins will induce membrane curvature so that  the observed spontaneous curvature must be a function of the local protein concentration \cite{SvetinaS2015a,MijoS2015a,IvankinA2012a,PartonD2011a,SchmidtN2013a}.  We define the membrane curvature induced by a single membrane protein as the spontaneous (membrane) curvature of the protein. Here we define $H_0$ as the average spontaneous curvature of lipids and proteins weight by their respective surface coverage fraction:
\begin{equation} \label{eqn:H_1}
H_0 = \sqrt{2} \frac{\dd{\sum_{l=0}^m C_0^l (a_{l}^{\rm lip})^2 \rho_l^{\rm lip} + C_0^{\rm pro} (a^{\rm pro})^2 \rho^{\rm pro}}}
{\dd{\sum_{l=0}^m (a_{l}^{\rm lip})^2 \rho_l^{\rm lip} + (a^{\rm pro})^2 \rho^{\rm pro}}},
\end{equation}
where $C_0^l $ and $ C_0^{\rm pro}$ are the spontaneous curvature of the $l^{th}$ species of lipids and proteins, respectively. Considering that the membrane surface is completely covered by the lipids and proteins, the following saturation constraint holds true:
\begin{equation} \label{eqn:saturation}
\sum_{l=0}^{m} (a_l^{\rm lip})^2 \rho_l^{\rm lip} + (a^{\rm pro})^2 \rho^{\rm pro} = 1.
\end{equation} 
With this constraint we can write the spontaneous curvature in Eq. (\ref{eqn:H_1}) as
\begin{equation} \label{eqn:H_2}
H_0 = \sqrt{2} \left( \sum_{l=0}^m C_0^l (a_{l}^{\rm lip})^2 \rho_l^{\rm lip} + C_0^{\rm pro} (a^{\rm pro})^2 \rho^{\rm pro} \right)
\end{equation}
and the membrane entropic energy as
\begin{eqnarray} \label{eqn:Gent_2}
G_{\rm ent} & = & \frac{1}{\beta} \int_{C} \left \{ \frac{1}{(a_0^{\rm lip})^2} \left( 1 - \rho^{\rm pro} (a^{\rm pro})^2 - \sum_{l=1}^m \rho_l^{\rm lip} (a_l^{\rm lip})^2 \right ) 
\right . \nonumber \\
 & & \times \left [  \ln \left(  1 - \rho^{\rm pro} (a^{\rm pro})^2 - \sum_{l=1}^m \rho_l^{\rm lip} (a_l^{\rm lip})^2 \right) - 1 \right ] + \nonumber \\
 & & \left . \sum_{l=1}^m \rho_l^{\rm lip} \left[ \ln(\rho_l^{\rm lip} (a_l^{\rm lip})^2) -1 \right ]  + \rho^{\rm pro} \left( \ln(\rho^{\rm pro}(a^{\rm pro})^2 ) - 1 \right ) \right \}ds.
\end{eqnarray}
To obtain the dynamics of the membrane morphology, one can calculate the variation of the total energy $G_{\rm tot}$ in Eq. (\ref{eqn:total_eng_curvprotein}) and solve the resulting equation for the gradient flow of $\phi$. This process is routine and can be found in the studies of spontaneous curvature effects of pure or multi-component membranes without proteins \cite{DuQ2006a,DuQ2005d}. Since our interest here is to investigate the protein localization on membrane surfaces, we choose to fix the membrane morphology, i.e., $H_0$ is a time-independent function. We then only need to calculate the variation of the total energy with respect to the membrane protein concentration, which turns out to be
\begin{eqnarray} \label{eqn:var_Gtot}
\frac{\delta G_{\rm tot}}{\delta \rho^{\rm pro}} & = & \frac{\delta G_{\rm mem}}{\delta \rho^{\rm pro}} + \frac{\delta G_{\rm ent}}{\delta \rho^{\rm pro}} \nonumber \\
 & = & k_B T \left [-\left( \frac{a^{\rm pro}}{a_0^{\rm lip}} \right)^2  \ln \left( 1 - \rho^{\rm pro} (a^{\rm pro})^2 - \sum_{l=1}^m \rho_l^{\rm lip} (a_l^{\rm lip})^2 \right) 
+ \ln (\rho^{\rm pro} (a^{\rm pro})^2) \right ]  \nonumber \\
& & + 2 C_0^{\rm pro} (a^{\rm pro})^2 ( H - H_0).
\end{eqnarray}

\subsection{Eulerian formulation}

While we are working on the membrane with fixed morphology, the formulation of the curvature driven protein localization is expected to interface with dynamical morphology where the membrane surface is not {\it a prior} known. For that purpose one could trace the position of membrane implicitly by evolving a phase field function $\phi(\vecx)$ on surface $S$ embedded in $\Omega\in{\mathbb R}^3$, where $\phi$ takes the value of $-1$ in the exterior of the membrane enclosure and $1$ inside \cite{DuQ2006a,DuQ2005d}. The membrane mean curvature at $\phi=0$ can be computed as a function of $\phi$ following
\begin{equation}
H = \frac{\sqrt{2} \varepsilon}{2 { (1-\phi^2)}} \left( \Delta_{\vecx} \phi + \frac{1}{\varepsilon^2} (1-\phi^2)\phi \right ),
\end{equation}
where $\varepsilon>0$ is a small parameter that adjust the transition of $\phi$ from $-1$ to $1$ near the membrane as in Eq. (\ref{eqn:geo_curv_energy}). We then identify three components of the chemical potential defined by the variation in Eq. (\ref{eqn:var_Gtot}) 
\begin{eqnarray}
L^{\rm pro} & = & \ln (\rho^{\rm pro} (a^{\rm pro})^2),  \label{eqn:L} \\
R^{\rm pro} & = & -\left( \frac{a^{\rm pro}}{a_0^{\rm lip}} \right)^2 \ln \left( 1 - \rho^{\rm pro} (a^{\rm pro})^2 - \sum_{j=1}^m  \rho_j^{\rm lip} (a_j^{\rm lip})^2\right ), \label{eqn:R} \\
P^{\rm pro} & = &  \frac{\varepsilon}{\sqrt{2}{ (1-\phi^2)}} \left( \Delta_{\vecx} \phi+ \frac{1}{\varepsilon^2} \phi(1-\phi^2  ) \right ) - H_0 \label{eqn:P} 
\end{eqnarray}
to write this chemical potential as
\begin{eqnarray}
\mu^{\rm pro} & = & \frac{\delta G_{\rm tot}}{\delta \rho^{\rm pro}} = k_BT ( L^{\rm pro} + R^{\rm pro}) + 2 C_0^{\rm pro} (a^{\rm pro})^2 \nabla_{\vecx} P^{\rm pro}.
\end{eqnarray}
This chemical potential allows us to define the diffusion flux vector and the transportation equation. Two options are available for the definition of the transportation equation. One could extract the membrane surface $S$ from the phase field function $\phi$ and solve a surface transportation on $S$. This involves the dynamic meshing or mesh deformation if $\phi$ is evolving in time, and singularity will arise if there is topological change in $S$ as $\phi$ evolves. 

Alternatively, one could formally define a 3D transportation equation in the entire domain $\Omega$ but practically restrict the transportation of membrane proteins to a very small neighborhood near the membrane surface $S$. This is accomplished by introducing to the flux vector
\begin{equation}
\vec{J}^{\rm pro}(\vecr) =-D^{\rm pro} \delta_{S} \beta \rho^{\rm pro}(\vecr) \nabla \mu^{\rm pro}
\end{equation}  
a function $\delta_{S}$ which is concentrated at the membrane $S$ where $\phi=0$. Various choices of such functions are available and their numerical properties differ subtly \cite{LeeH2012a}. We choose 
\begin{equation}
\delta_{S} = \begin{cases}
\tanh (10(\phi+1)), & -1 \le \phi \le 0, \\
-\tanh (10(\phi-1)), & 0 \le \phi \le 1, \\
\end{cases}
\end{equation}
so that effective domain near $\phi=0$ can be automatically identified as $\phi$ evolves. The general transportation equation for membrane proteins reads
\begin{equation} \label{eqn:membrane_protein_eq}
\frac{\p \rho^{\rm pro}(\vecr)}{\p t} + \nabla \cdot (\vec{v} \nabla \rho^{\rm pro}(\vecr)) ={ -} \nabla \cdot \vec{J}^{\rm pro}(\vecr),
\end{equation}
where $\vec{v}$ is the velocity of the membrane in which the membrane proteins move. Although this velocity is taken to be zero in our computations simulations to be presented here, it can be computed if the membrane moves with the evolving phase field function. The nature of the equation can be seen if the size effects of  lipids and membrane proteins are not considered, i.e., $a_l^{\rm lip} = a^{\rm pro} = 0$. In this case $R^{\rm pro}=0$ and
\begin{equation} \label{eqn:membrane_protein_eq_2}
\frac{\p \rho^{\rm pro}}{\p t} = \nabla \cdot ( D^{\rm pro} \delta_{S} \nabla \rho^{\rm pro} + 2 k_B T D^{\rm pro} C_0^{\rm pro} (a^{\rm pro})^2 \delta_{S} \rho^{\rm pro} \nabla P^{\rm pro}),
\end{equation}
which is a drift-diffusion equation with a potential $P^{\rm pro}$. The mean curvature of the membrane therefore appears a potential that drives the 
transportation of membrane proteins to membrane surfaces where its mean curvature well fits the spontaneous membrane curvature of proteins. To 
numerically solve the equation, we separate the linear and nonlinear components of the equation, which are then treated using an implicit-explicit splitting 
interaction methods similar to the treatment of Eq. (\ref{eqn:surface_CHE}) presented in Sect.~\ref{sec:surfpattern}. The spatial approximation of 
the equation is obtained by using the Fourier spectral method, and a change of variable is necessary to convert the equation with variable 
diffusion coefficient $D \delta_{S}$ to a constant diffusion coefficient so that the Fourier spectral
method is applicable \cite{ConcusP1973a,StrainJ1994a}. 

\subsection{Computational simulations and summary}
\begin{figure}[!htb]
 \centering
\begin{tabular}{lll}
 \includegraphics[width=3.125cm]{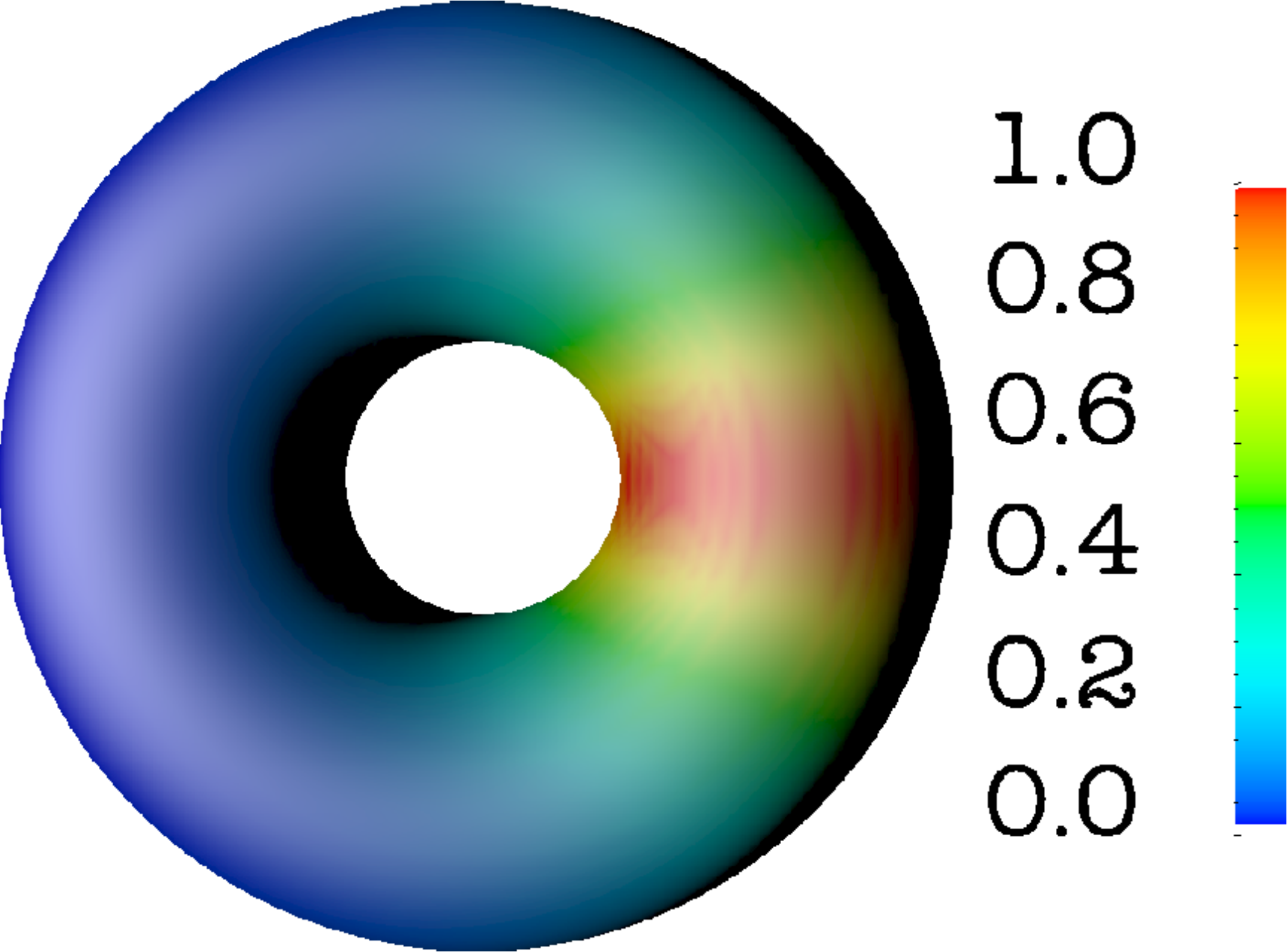}   &  
 \includegraphics[width=3.125cm]{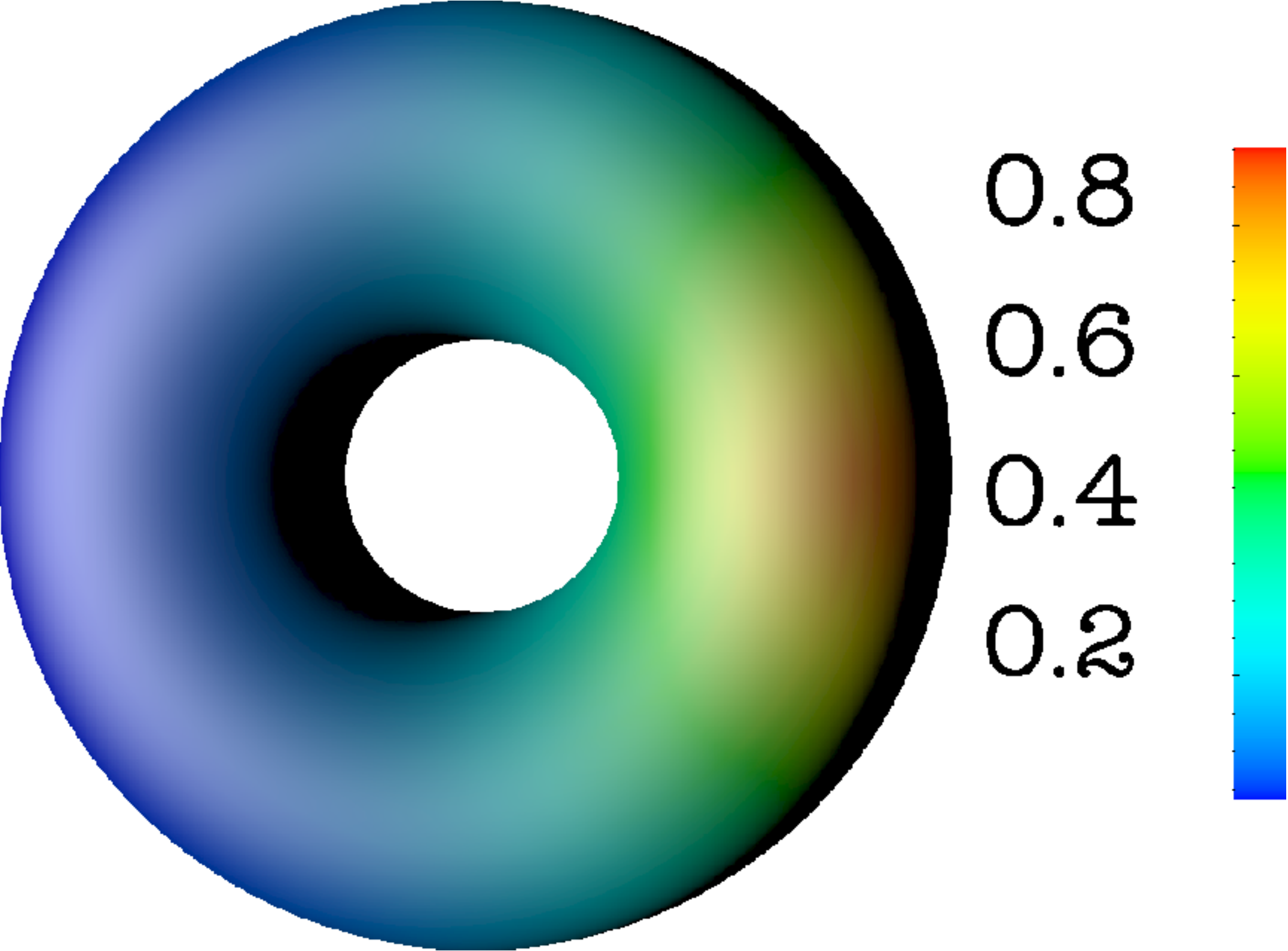} &
 \includegraphics[width=3.125cm]{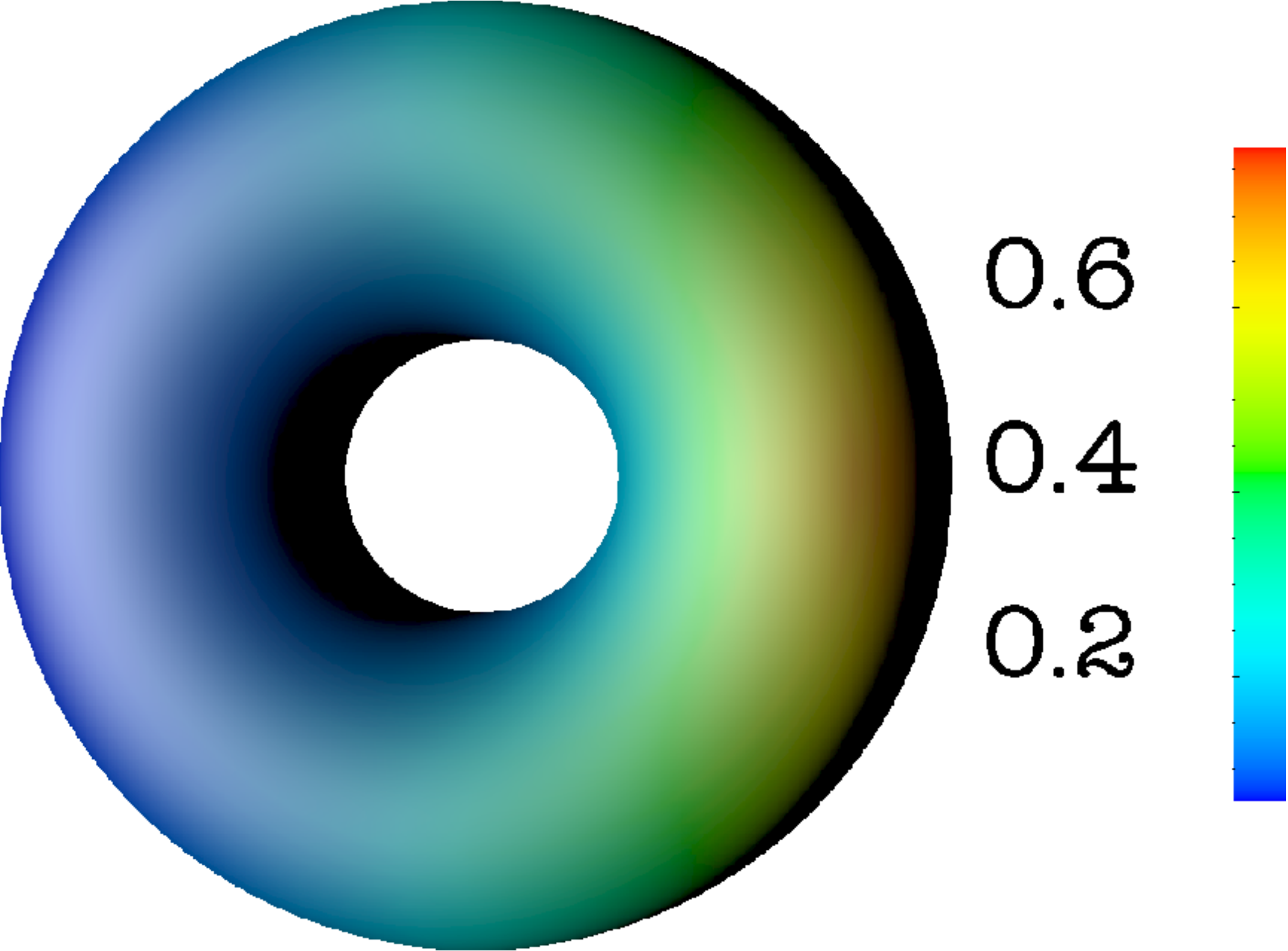} \\ 
 \includegraphics[width=3.125cm]{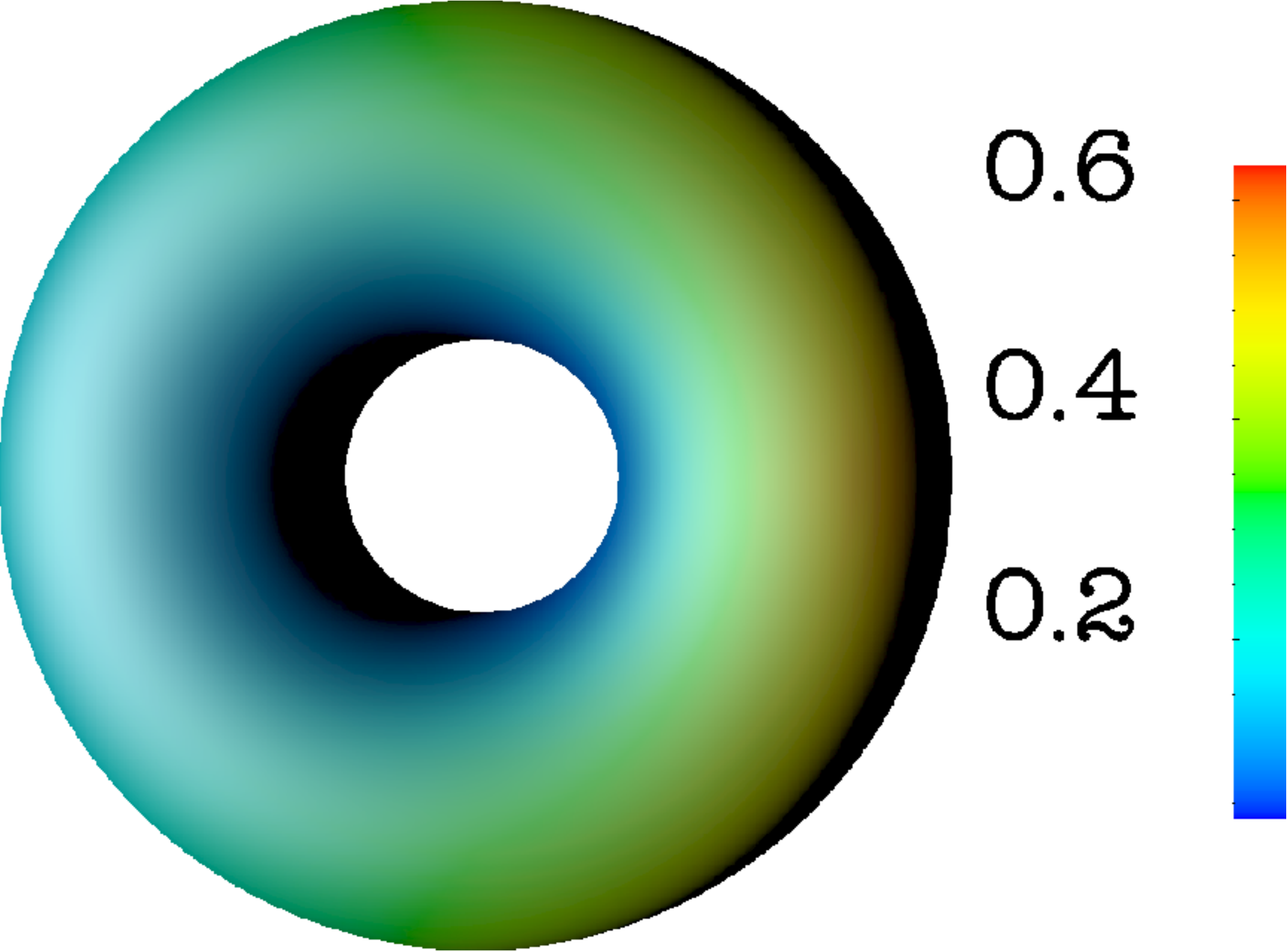} &
 \includegraphics[width=3.125cm]{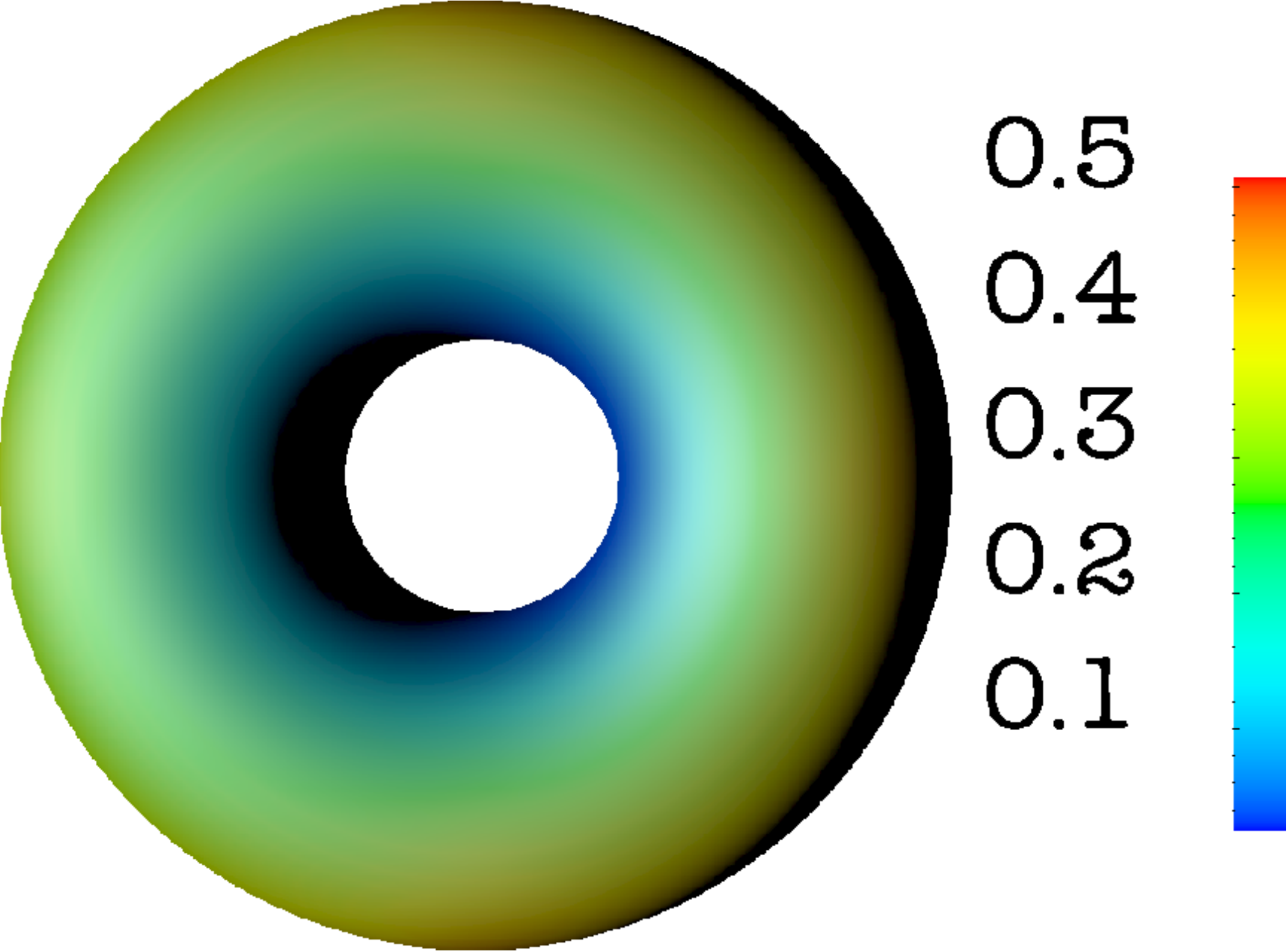} &
 \includegraphics[width=3.125cm]{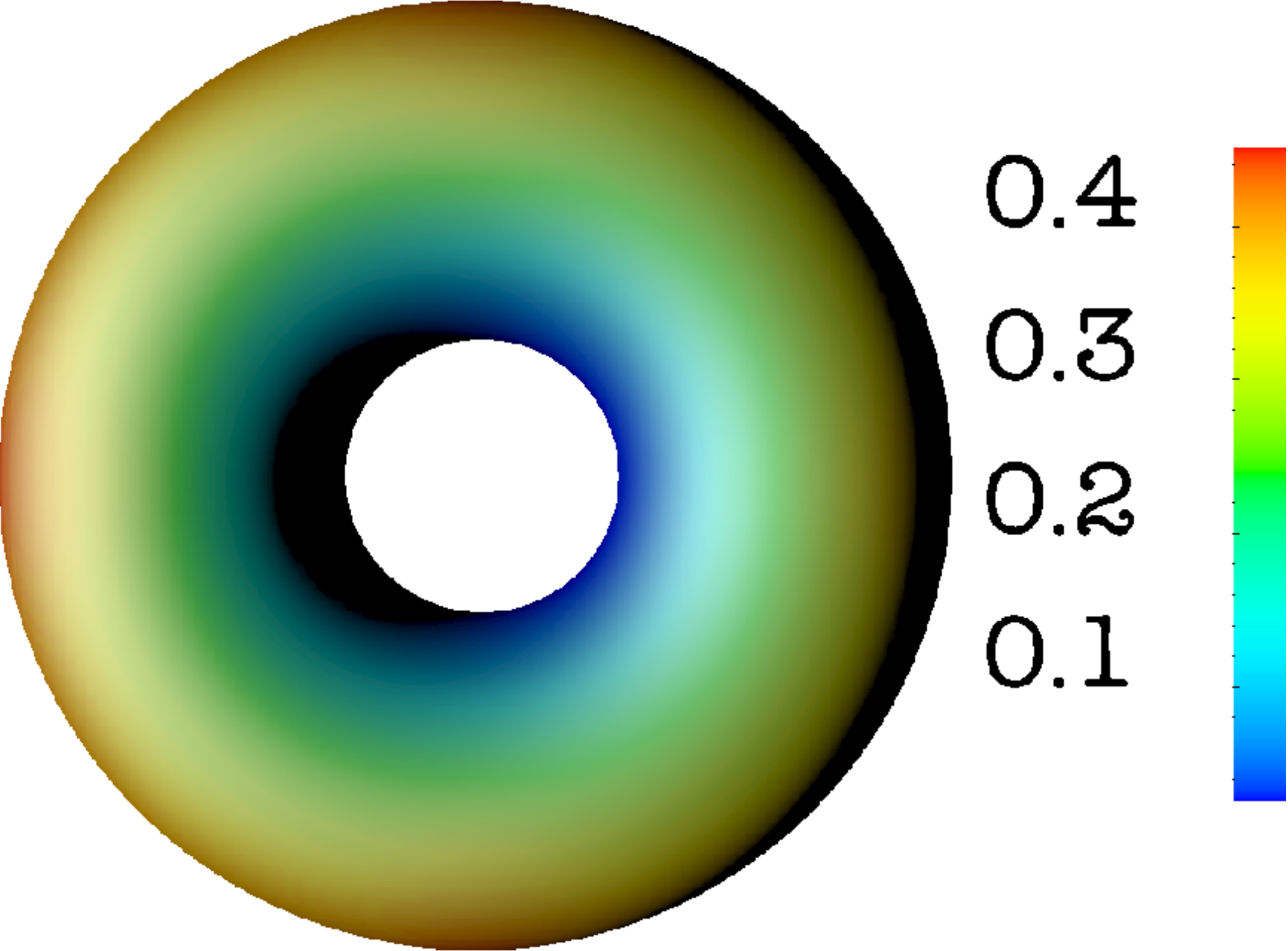} \\ 
  & & \\
 \includegraphics[width=2.75cm]{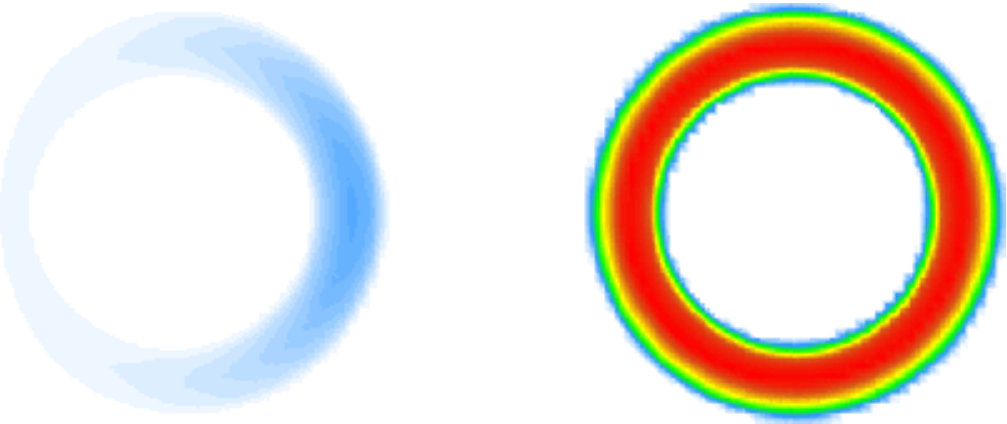}   &  
 \includegraphics[width=2.75cm]{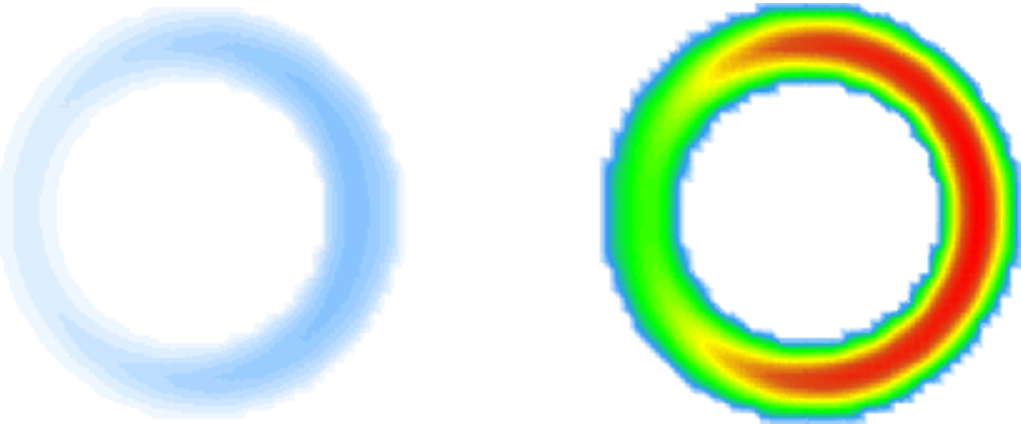} &
 \includegraphics[width=2.75cm]{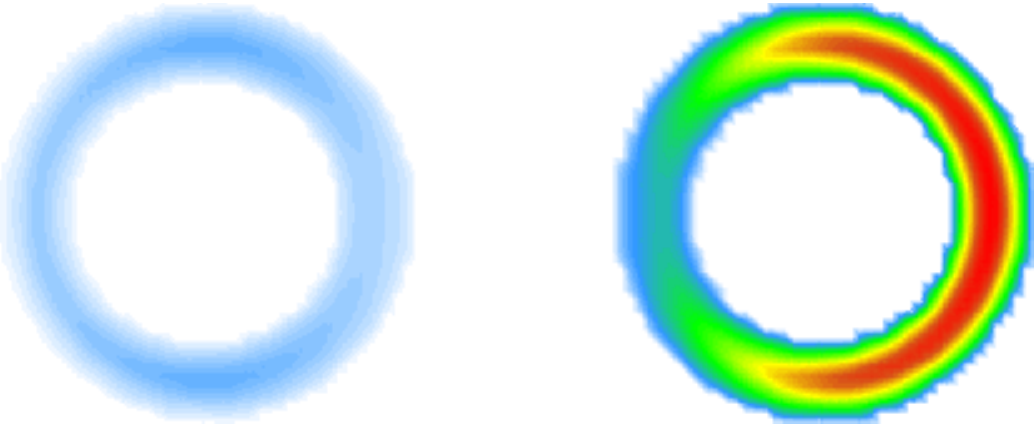} \\
  & & \\
 \includegraphics[width=2.75cm]{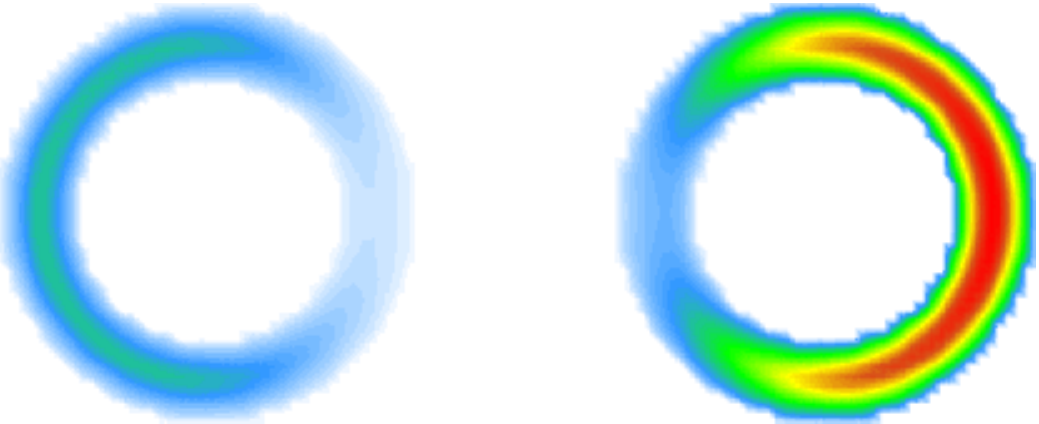} &
 \includegraphics[width=2.75cm]{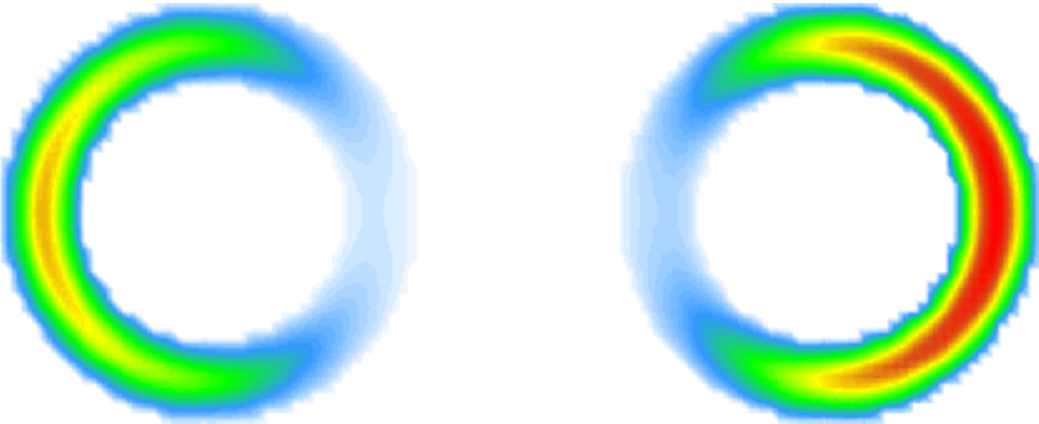} &
 \includegraphics[width=2.75cm]{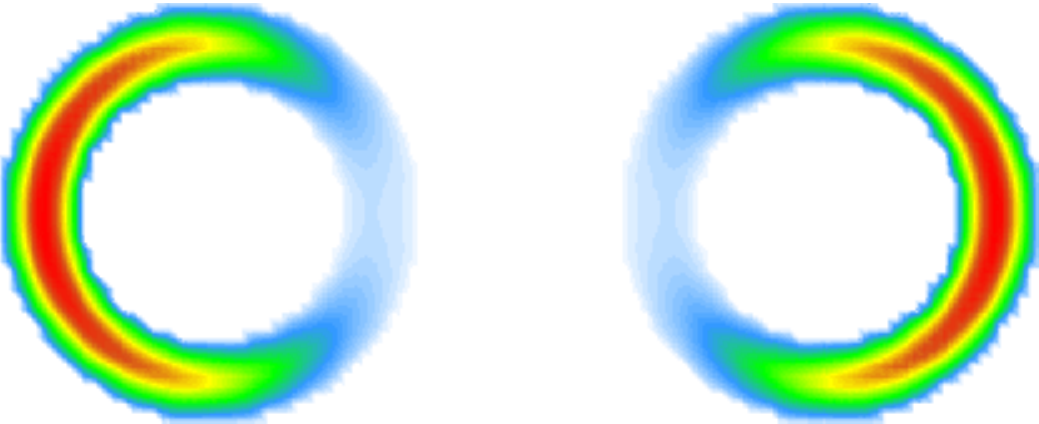} 
\end{tabular} 
\captionsetup{width=0.85\textwidth}
   \caption{Simulated localization of the membrane proteins from its initial position to the outer ring of the torus on a $128^3$ uniform mesh.
$\varepsilon = 0.1$. Time increment $\Delta t = 10^{-3}$. Spontaneous curvatures $C_0^{\rm pro} = 0.5, C_0^{\rm lip} =-0.1$, 
and sampling moments are $t=0, 0.1, 0.25, 0.5, 1.0, 5.0$. Color is scaled by the maximum concentration in each plot.} 
   \label{fig:protein_sdiff_1} 
\end{figure}

To demonstrate the curvature preference of protein localization we consider in the domain $\Omega =(-4,4)^3$  a torus because it has regions with positive and negative mean curvatures where the proteins may populate or not  depending on their spontaneous curvature. The torus surface is given by 
\begin{equation}
(R - \sqrt{x^2 + y^2})^2 + z^2 = r^2, 
\end{equation}
where $R$ and $r$ are the major and minor radii, respectively. Its alternative parametrization 
\begin{equation}
(x,y,z) = ((R+r\cos \theta) \cos \varphi, (R+r\cos\theta) \sin \varphi, r \sin \varphi)
\end{equation}
can be handy when computing the curvature. Here $0 \le \theta \le 2 \pi$ is the angle made from the surface around the center of the tube, known as the poloidal angle, and $0 \le \varphi \le 2 \pi$ is the angle made from the surface to the positive $x$-axis (projected on the $xy$-plane), known as the toroidal angle. When $R>r$, one gets the so-called ring torus. Here we choose $R=2$  and $r=1.1$. The phase field function $\phi$ is set as the signed distance function  with this torus surface. We consider only one species of diffusion proteins and one species of lipids. The saturation condition (\ref{eqn:saturation}) then indicates that we only need to model the distribution of proteins only. The membrane proteins are initially concentrated near the highest point of the positive $y$-axis, smoothly distributed along the surface, and because of the adoption of phase field function which expands the transportation domain from the surface to a small neighborhood in the vicinity of the surface, smoothly distributed from the surface to the bulk:
\begin{equation}
\rho = \rho_0 e^{-\sqrt{x^2+(y-R)^2 + z^2}} e^{-2\left (r - \sqrt{(x-c_x)^2 + (y - c_y)^2 + z^2} \right) },
\end{equation}
where $r = \sqrt{x^2 + y^2 + z^2}$ and $(c_x,c_y,0)$ is the center of the torus tube on the same plane of which locates the point $(x,y,z)$.  The scaling constant $\rho_0$ is chosen such that the maximum of the concentration is $1$ on the torus surface. 
\begin{figure}[!htb]
 \centering
 \begin{tabular}{lll}
 \includegraphics[width=3.125cm]{Torus_N128_t0B}   &  
 \includegraphics[width=3.125cm]{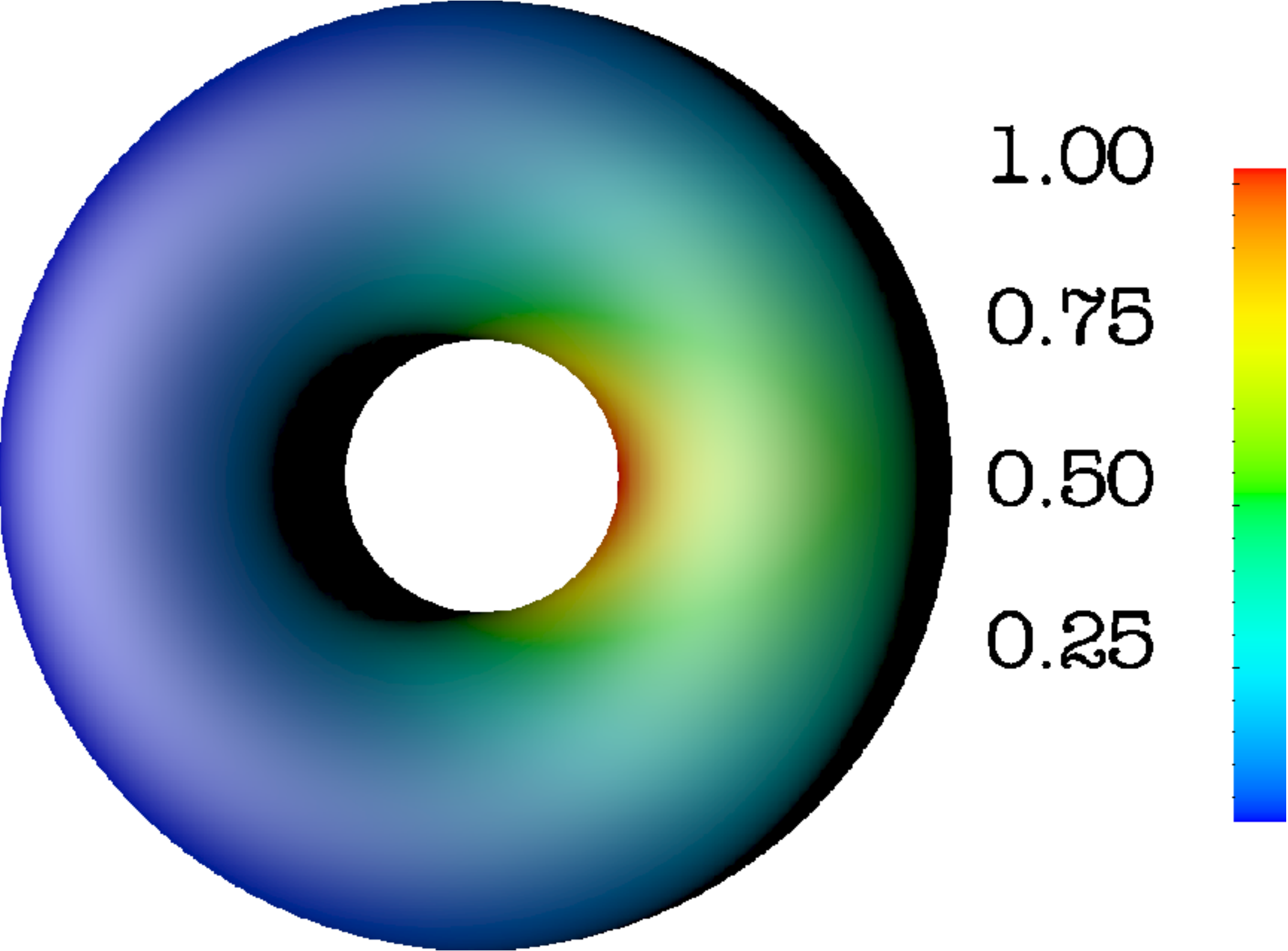} &
 \includegraphics[width=3.125cm]{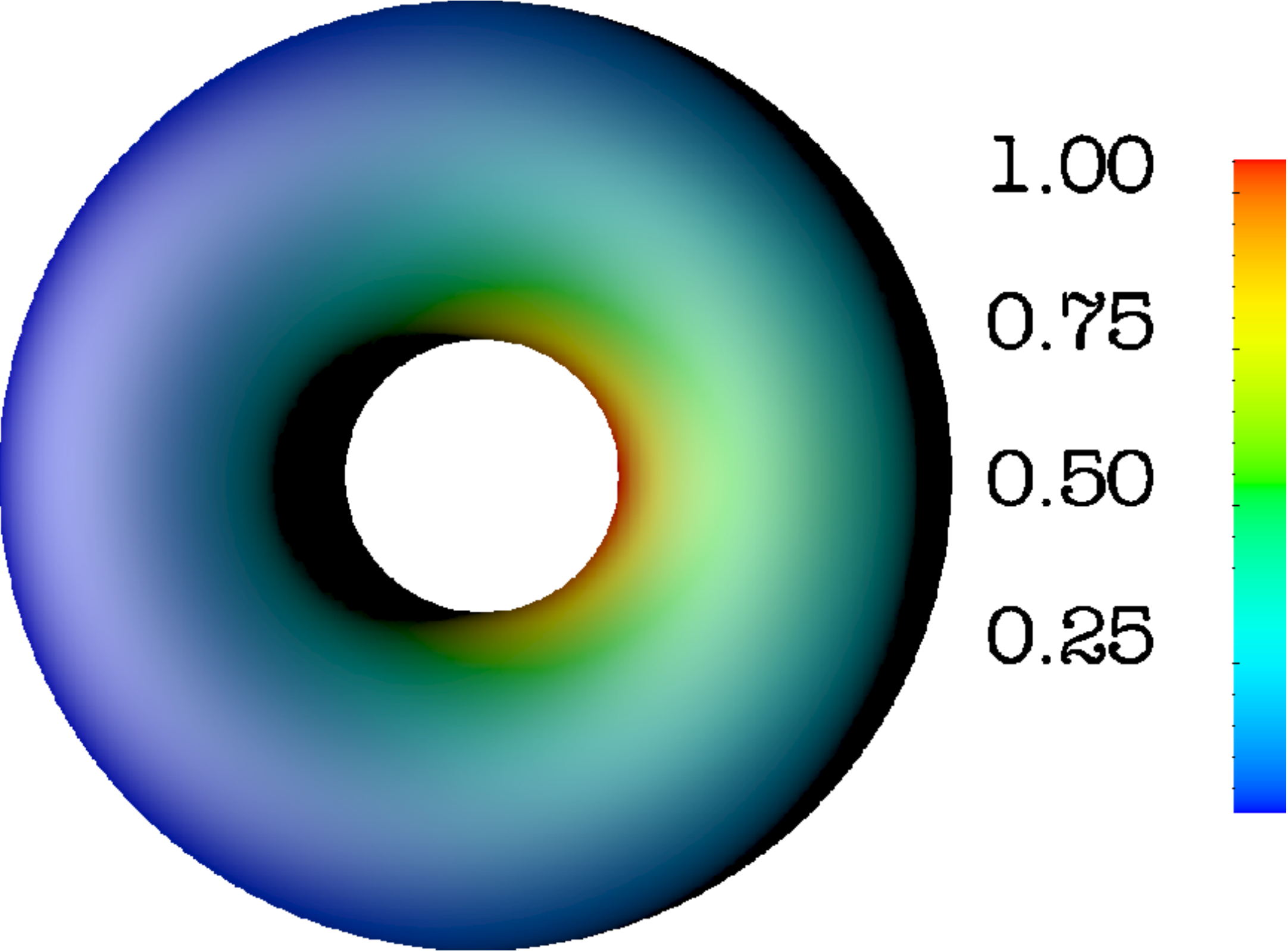} \\ 
 \includegraphics[width=3.125cm]{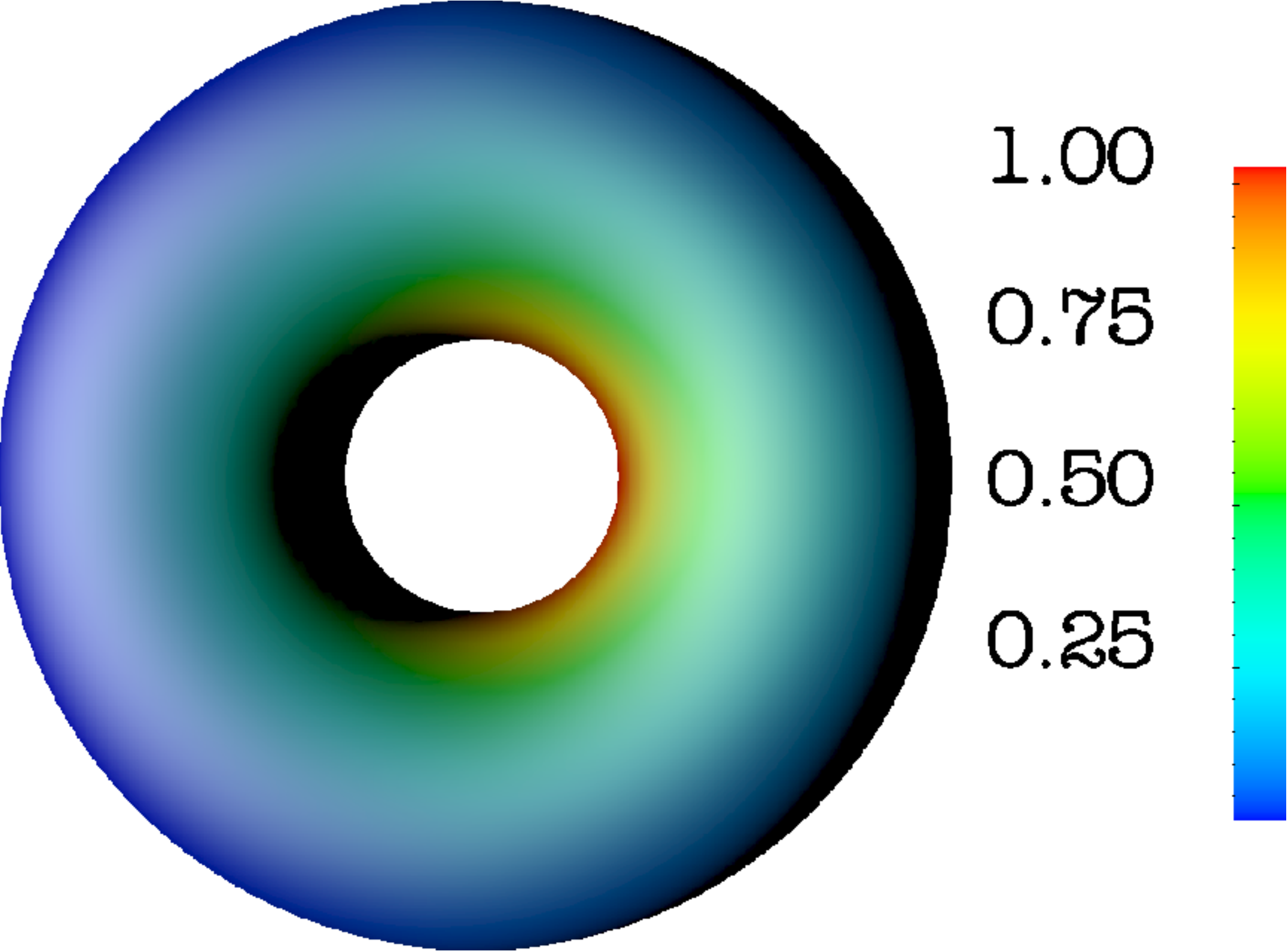} &
 \includegraphics[width=3.125cm]{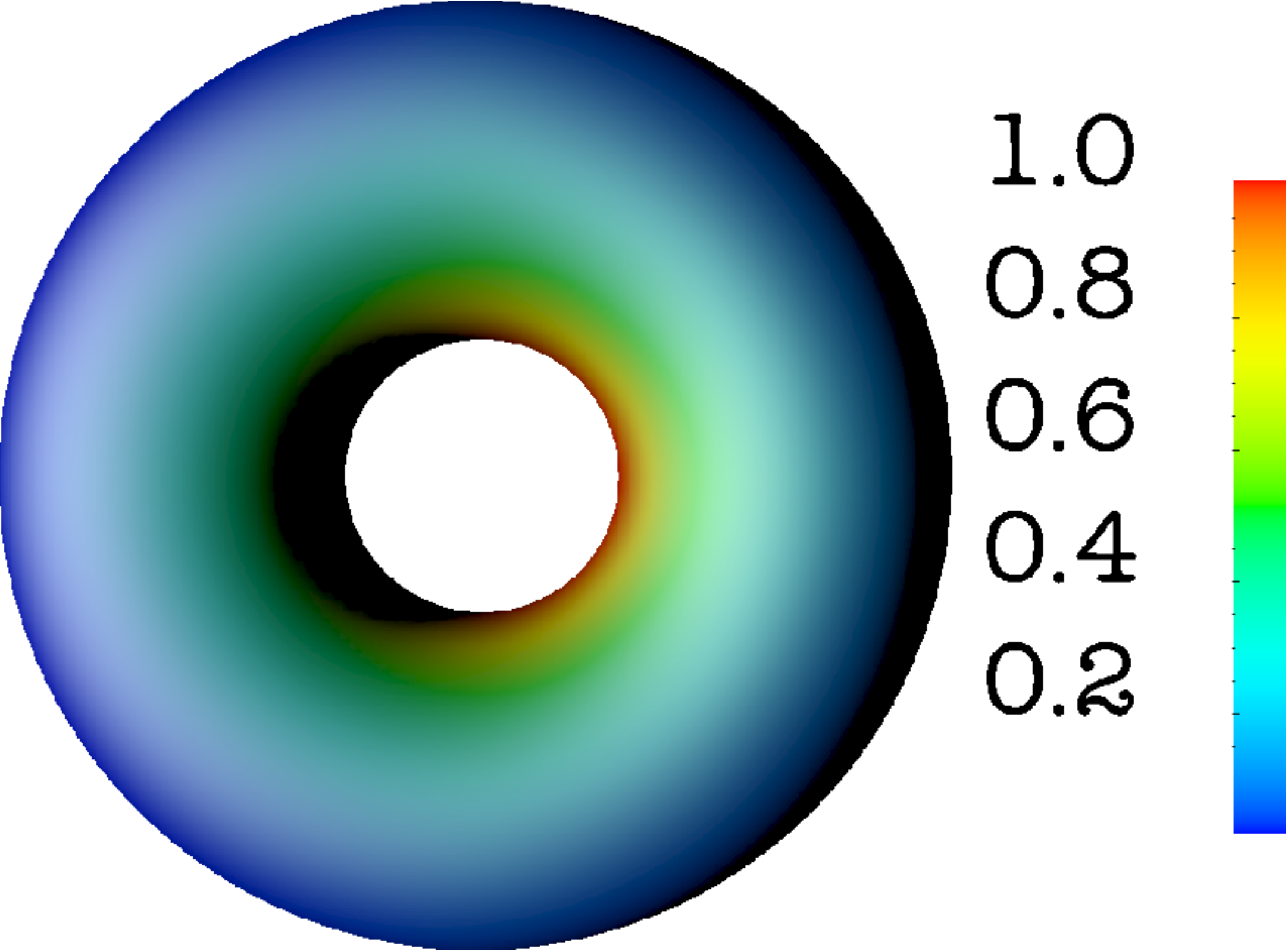} &
 \includegraphics[width=3.125cm]{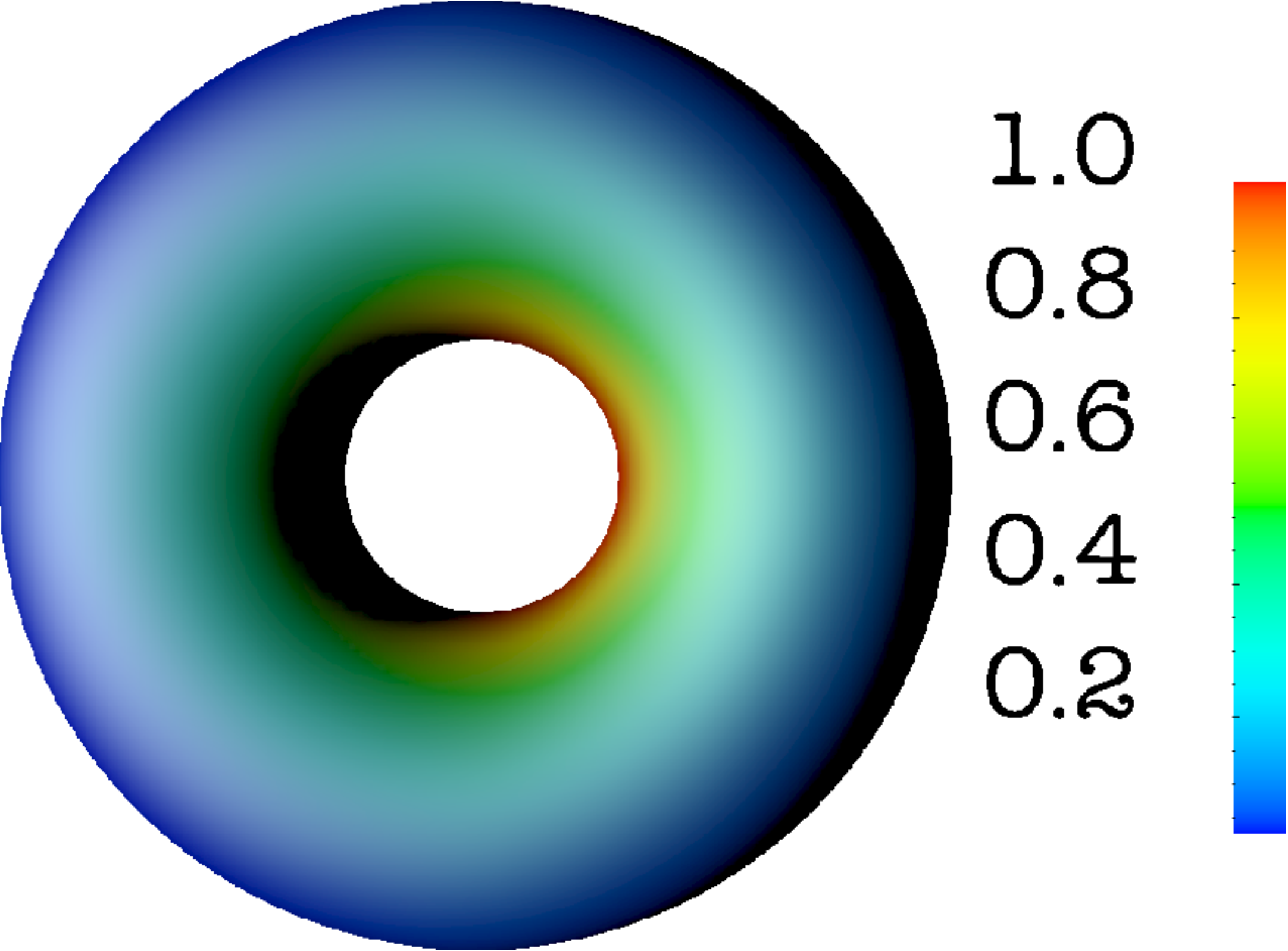} \\
  & & \\
 \includegraphics[width=2.75cm]{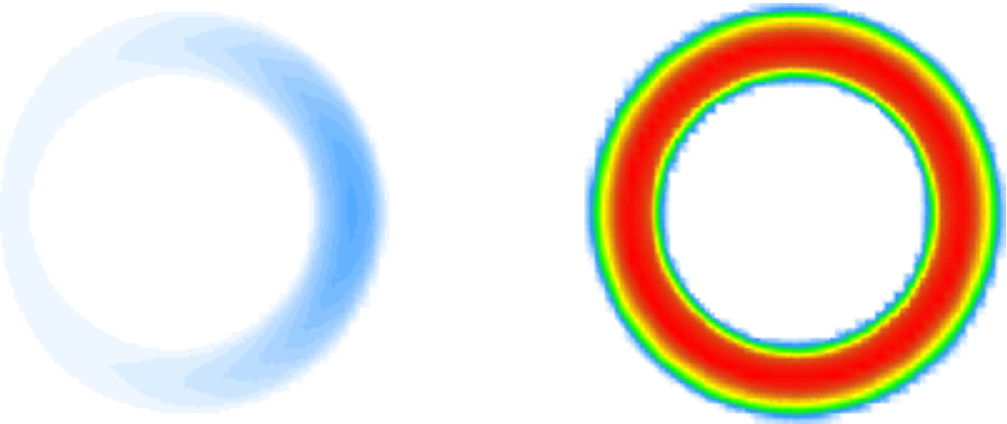}   &  
 \includegraphics[width=2.75cm]{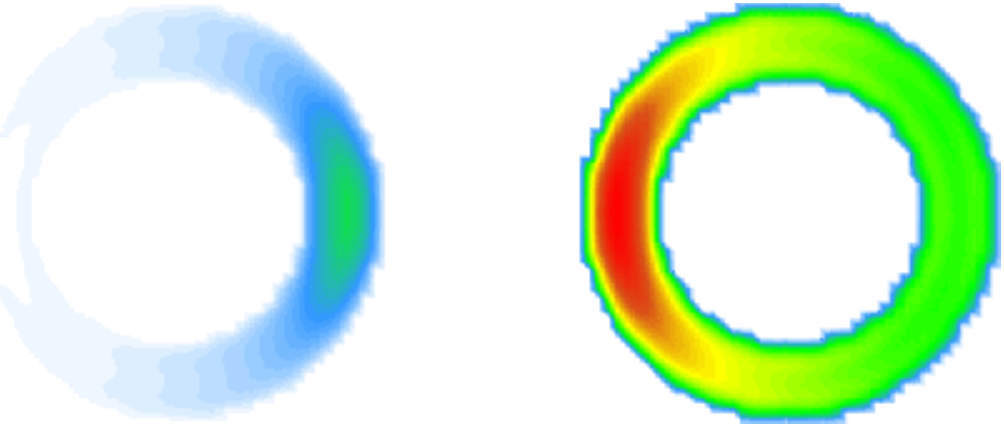} &
 \includegraphics[width=2.75cm]{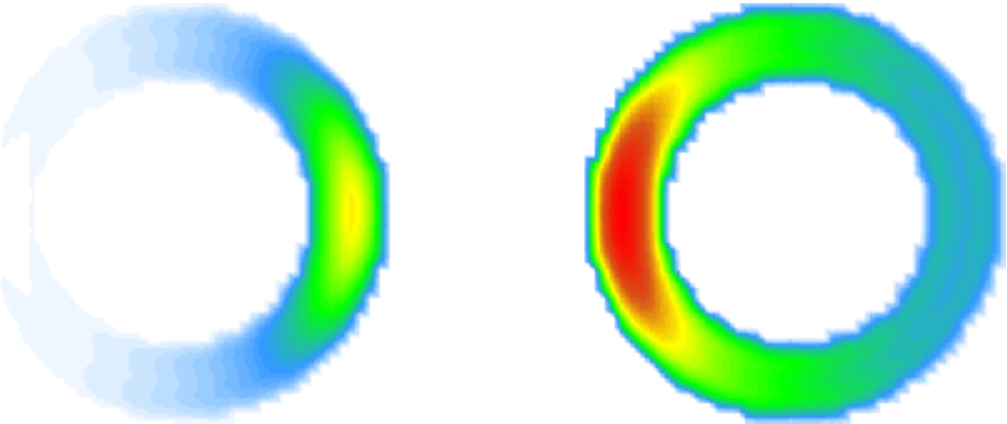} \\ 
  & & \\
 \includegraphics[width=2.75cm]{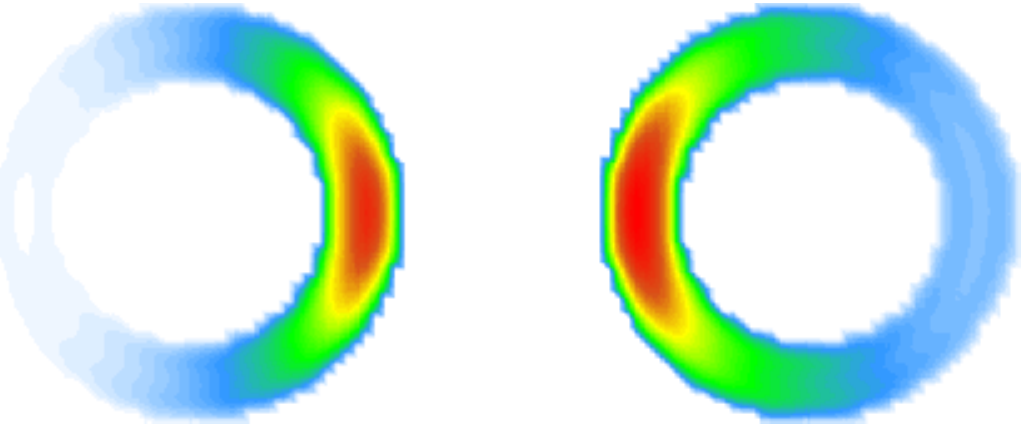} &
 \includegraphics[width=2.75cm]{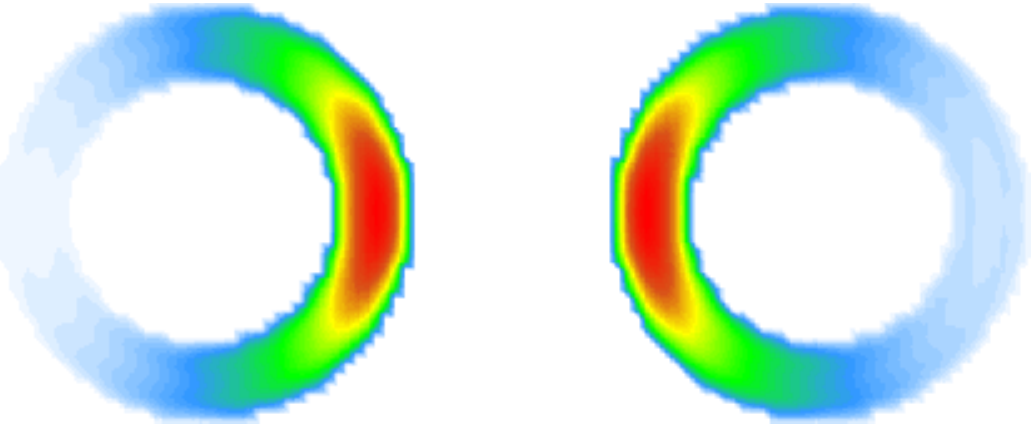} &
 \includegraphics[width=2.75cm]{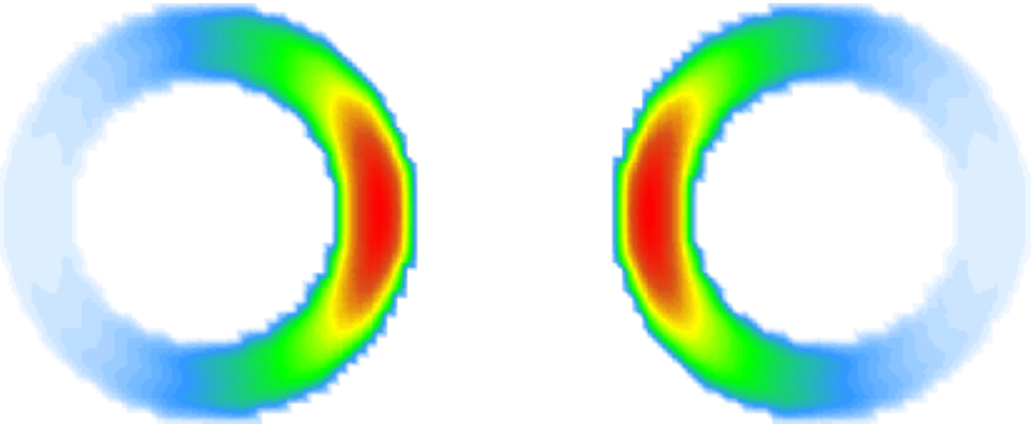} 
 \end{tabular}
\captionsetup{width=0.85\textwidth}
   \caption{Simulated localization of the membrane proteins from its initial position to the outer ring of the torus on a $128^3$ uniform mesh. $\varepsilon = 0.1$. Time increment is $\Delta t = 10^{-3}$. Spontaneous curvatures are $C_0^{\rm pro} =-0.1 $ and $  C_0^{\rm lip} =0.5$, 
and sampling moments are $t=0, 0.1, 0.25, 0.5, 1.0, 5.0$. Color is scaled by the maximum concentration in each plot.} 
  \label{fig:protein_sdiff_2} 
\end{figure}

We first set the spontaneous curvature of membrane proteins and lipids to be $C_0^{\rm pro} = 0.5, C_0^{\rm lip} = -0.1$, respectively. Notice that the mean curvature of a torus is given by
\begin{equation}
H_{\rm torus} = \frac{R + 2r \cos \theta}{2r(R + r \cos \theta)},
\end{equation}
which gives a mean curvature $H_{\rm torus} \approx 0.6158$ for the chosen values of $R,r$ at the outer ring of the torus where $\theta =0$ and $H_{\rm torus} \approx -0.1$ at the inner ring of the torus where $\theta = \pi$. With this first choice of $C_0^{\rm pro},C_0^{\rm lip}$ we expect that the membrane proteins will populate near the outer ring where the mean curvature is close to the specified spontaneous curvature of membrane proteins. Our expectation is verified by Fig.~\ref{fig:protein_sdiff_1}, where the plots of the concentrations  of the membrane proteins on the membrane $\phi=0$ and the cross section $y=0$ at six sampling moments show the transportation of membrane proteins from its initial position to the outer ring of the torus.

In the second simulation we start with same initial condition as in the first simulation but switch the spontaneous curvatures  to  $C_0^{\rm pro} =-0.1 $  and $  C_0^{\rm lip} = 0.5$. It is expected that the membrane proteins will finally populate at the inner ring of the torus,
and this is verified by the snapshots of concentrations in Fig.~\ref{fig:protein_sdiff_2}. 

These two computational simulations demonstrate the successful modeling of the curvature driven membrane protein localization using the drift-diffusion equation (\ref{eqn:membrane_protein_eq_2}). Full version of Eq. (\ref{eqn:membrane_protein_eq}) can  also be considered to include the effects of finite sizes of effects of lipids and proteins, and multiple species of lipids. Our choice of small time increment ($\Delta = 10^{-3}$) is restricted by the stability of the implicit-explicit splitting method used for integrating the nonlinear equation. We expect the development of more efficient numerical methods for the integration of the equation, in particular when it is to be coupled with the dynamic phase field function $\phi$, in that case a membrane velocity shall be added to Eq. (\ref{eqn:membrane_protein_eq_2}) to make it an advection-drift-diffusion equation. Such coupling  reveals the positive feedback of membrane curvature accumulation to membrane protein localization. On the other hand, the number of major membrane proteins involved in the membrane fusion, budding, endocytosis, or exocytosis is a not constant over the entire time course because there is continuous intracellular protein transport. Proteins may be recruited from the solution to membrane at specific regions  of the membrane and meanwhile they are released from the membrane to the solution \cite{RothmanJ1994a,SollnerT1993a}. The model presented here can be extended by adding a reaction term that models the  dynamic exchange of membrane proteins between the membrane and the solution. Indeed, it is shown that some membrane budding proteins such as influenza virus hemagglutinin (HA) and neuraminidase (NA) are associated with raft-like microdomains, while some are not \cite{LeserG2005a}. An integration of the curvature driven localization and local clustering within the microdomains will help elucidate the competing  or collaborative effects of these membrane proteins in the same biophysical process.

\section{Conclusions}
Energetic variational principle constitutes a tangible link between multiscale theory and the experimental observation of biomolecular structure, function,  and dynamics,  aided by computational simulations. Although the applications of variational principle have been well established for research in various areas of mechanics, classical and modern physics, and material sciences, novel insights are offered by this principle when it is applied  to the biomolecular systems.   Among the progresses achieved in recent years, a significant step forward has been made using the geometry of the molecular interface to parametrize the total energy \cite{WeiG2012a,ChenZ2010a,ChenZ2011a,LiB2011a,LiB2009a,Melissa_PhDThesis,Mikucki_PhDThesis,HyonY2014a}.  This unified representation allows the investigators to focus on the identification of energies that characterize various molecular interactions at multiple spatial and temporal scales. The flexibility of the analytical and computational framework of the variational principle ensures that the critical states and dynamics of the biomolecular system can be tracked with confidence by evolving the total energy. Furthermore, by introducing a phase field function we can implicitly define and track the molecular interface which may subject to large deformation and topological change. The three topics presented here demonstrated the desirable flexibilities of formulating the total energy, of parametrizing the energy using phase field function, and of simulating the equilibrium state and dynamics of the system though the numerical solutions of the nonlinear partial differential equations (PDEs) for the geometric flow of the total energy.  

The geometrically parametrized total energy obtained by the energetic variational principles entails a rich body of features for mathematical and numerical analysis, including the stability of its critical points, the coarsening dynamics, the solution periodicity, and the conservative discretization of the resulting PDEs, while most of them remain open as long as the applications to biomolecular problems are concerned. More broad usefulness of the methodology outlined in the present three topics are expected to be established in chemistry, biophysics, and medicine through interdisciplinary research and collaboration.

\bibliographystyle{plain}

\end{document}